\documentclass[11pt]{amsart}

\usepackage{amssymb}

\input xy

\xyoption{all}

\theoremstyle{plain}

\newtheorem{theorem}{Theorem}[section]

\newtheorem{corollary}[theorem]{Corollary}

\newtheorem{lemma}[theorem]{Lemma}

\newtheorem{proposition}[theorem]{Proposition}

\newtheorem{fact}[theorem]{Fact}

\newtheorem{claim}[theorem]{Claim}

\theoremstyle{definition} 

\newtheorem{definition}[theorem]{Definition}

\newtheorem{example}[theorem]{Example}

\newtheorem{remark}[theorem]{Remark}

\theoremstyle{remark}

\newtheorem{hypothesis}[theorem]{Hypothesis}

\newcommand{\cf}{\operatorname{cf}}

\newcommand{\Min}{\operatorname{Min}}

\renewcommand{\phi}{\varphi}

\newcommand{\initial}\lessdot

\newcommand{\id}{\operatorname{id}}

\newcommand{\N}{\mathcal N}

\def\?{?\vadjust

{\vbox to 0pt{\vskip-7pt\hbox to 1.1\hsize{\hfill\huge ?!}}}}

\newcommand{\be}{\begin{enumerate}}
\newcommand{\ee}{\end{enumerate}}
\newcommand{\case}{\emph}
\newcommand{\step}{\emph}
\newcommand{\discussion}{\subsubsection*}
\usepackage{enumerate}
\usepackage{verbatim}
\renewcommand{\epsilon}{\varepsilon}

 \def\nfork{\setbox0\hbox{$\bigcup$}%
 \setbox1=\hbox to \wd0{\hfil\vrule width 0.7pt depth 2pt height 7.5pt\hfil}%
 \wd1=0cm\relax\box1\box0}
 \def\dnf{\mathop{\nfork}\limits}

\begin{document}
\title{Good Frames With A Weak Stability}

\thanks{This work is supported by Number 875 on Shelah's list of publications. It is a part of the PHD thesis of the first author with the
supervisor of Boaz Tsaban}

\author{Adi Jarden}
\email[Adi Jarden]{jardenadi@gmail.com}
\address{Department of Mathematics and computer science.\\ Bar-Ilan University \\ Ramatgan 52900, Israel}

\author{Saharon Shelah}
\email[Saharon Shelah]{shelah@math.huji.ac.il}
\address{Institute of Mathematics\\ Hebrew University\\ Jerusalem, Israel\\ and Rutgers University\\
 New Brunswick\\ NJ, U.S.A}


\maketitle

\begin{abstract}
Let $K$ be an abstract elementary class of models. Assume that
there are less than the maximal number of models in
$K_{\lambda^{+n}}$ (namely models in $K$ of power $\lambda^{+n}$)
for all $n$. We provide conditions on $K_\lambda$, that imply the
existence of a model in $K_{\lambda^{+n}}$ for all $n$. We do this
by providing sufficiently strong conditions on $K_\lambda$, that
they are inherited by a properly chosen subclass of
$K_{\lambda^+}$.

\end{abstract}

\tableofcontents

\section{Introduction}

The book classification theory, \cite{shc}, of elementary classes,
i.e. classes of first order theories, presents properties of
theories, which are so called ``dividing lines'' and investigates
them. When such a property is satisfied, the theory is low, i.e.
we can prove structure theorems, such as:
\begin{enumerate}
\item The fundamental theorem of finitely generated Abelian
groups. \item Artin–Wedderburn Theorem on semi-simple rings. \item
If $V$ is a vector space, then it has a basis B, and $V$ is the
direct sum of the subspaces $span\{b\}$ where $b \in B$.
\end{enumerate}
(we do not assert that these results follow from the model
theoretic analysis, but they merely illustrate the meaning of
`structure'). But when such a property is not satisfied, we have
\emph{non-structure}, namely there is a witness that the theory is
complicated, and there are no structure theorems. This witness can
be the existence of many models in the same power.

 There has been much work on classification of
elementary classes, and some work on other classes of models.

The main topic in the new book, (\cite{shh}), is \emph{abstract
elementary classes} (In short a.e.c.). There are two additional
books which deal with a.e.c.s (\cite{babook} and \cite{grbook}).

From the viewpoint of the algebraist, model theory of first order
theories is somewhat close to universal algebra. But he prefers
focusing on the structures, rather than on sentences and formulas.
Our context, abstract elementary classes, is closer to universal
algebra, as our definitions do not mention sentences or formulas.

As superstability is one of the better dividing lines for first
order theories, it is natural to generalize this notion to
a.e.c.s. A reasonable generalization is that of the existence of a
good $\lambda$-frame, (see Definition \ref{2.1a} on page
\pageref{2.1a}), introduced in \cite{sh600}. In \cite{sh600} we
assume existence of a good $\lambda$-frame and either get a
non-structure property (in $\lambda^{++}$, at least where
$2^\lambda<2^{\lambda^+}<2^{\lambda^{++}}$) or derive a good
$\lambda^+$-frame from it. Our paper generalizes \cite{sh600},
weakening the assumption of a good $\lambda$-frame, or more
specifically weakening the basic stability assumption.

\subsection{The required knowledge}
We assume basic knowledge in set theory (ordinals, cardinals,
closed unbounded subsets and stationary subsets). In model theory,
we just assume the reader is familiar with notions, every student
in algebra knows (theory, model=structure, isomorphism and
embedding). Especially we do not assume the reader is familiar
with formula and elementary substructure, as here we do not deal
with those notions (except in one example). Of course, we do not
assume the reader has read any paper in abstract elementary
classes, and if the reader prefers to translate a model as a
group, he will not lose the main ideas. We sometimes refer to
another paper, for the following four tasks:
\begin{enumerate}
\item To convince the reader that an assumption is reasonable,
i.e.
to prove that we can conclude something from its negation. \item
To give examples. \item To compare it with \cite{sh600}. \item To
point out its continuations.
\end{enumerate}

The best way to read this paper is to read it until its end,
before reading any reference.

\subsection{The assumptions}
When we write a hypothesis, we assume it until we write another
hypothesis, but usually we recall the hypothesis in the beginning
the following section. When we write `but we do not use local
character', the reader may wonder why? The answer is that we want
to apply theorems we prove here, in papers, in which local
character is not assumed (for example \cite{jrsh2}). For the same
reason, in Hypothesis \ref{3.0} we assume weak assumptions.

\subsection{Notations}
 We use the letters $m,n,k,l$ for natural numbers or
integer numbers, $\alpha,\beta,\gamma,i,j,\epsilon,\zeta$ for
ordinal numbers, $\delta$ for a limit ordinal number,
$\kappa,\lambda,\mu$ for cardinal numbers, $p,q$ for types, $P$
for a set of types, $K$ for a class of models, $\frak{k}$
,$\frak{s},U$ for specific uses.

\begin{definition}[Abstract Elementary Classes]\label{1.1}
\mbox{} 
\begin{enumerate}
\item  Let $K$ be a class of models for a fixed vocabulary and let
$\preceq=\preceq$ be a 2-place relation on $K$. The pair
$(K,\preceq)$ is an \emph{a.e.c.} if the following axioms are
satisfied:
\begin{enumerate}
\item  $K,\preceq$ are closed under isomorphisms. In other words,
if $M_1 \in K$, $M_0 \preceq M_1$ and $f:M_1 \rightarrow N_1$ is
an isomorphism then $N_1 \in K$ and $f[M_0]\preceq N_1$. \item
$\preceq$ is a partial order and it is included in the inclusion
relation. \item If $\langle M_\alpha:\alpha<\delta\rangle$ is a
continuous $\preceq$-increasing sequence, then $$M_0 \preceq
\bigcup \{M_\alpha:\alpha<\delta\} \in K.$$ \item Smoothness: If
$\langle M_\alpha:\alpha<\delta \rangle$ is a continuous
$\preceq$-increasing sequence, and for every $\alpha<\delta,\
M_\alpha\preceq N$, then
$$\bigcup \{M_\alpha:\alpha<\delta\} \preceq N.$$
\item If $M_0 \subseteq M_1 \subseteq M_2$ and $M_0 \preceq M_2
\wedge M_1 \preceq M_2$, then $M_0 \preceq M_1$. \item  There is a
Lowenheim Skolem Tarski number, $LST(K,\preceq)$, which is the
first cardinal $\lambda$, such that for every model $N \in K$ and
a subset $A$ of it, there is a model $M \in K$ such that $A
\subseteq M \preceq N$ and the cardinality of $M$ is $\leq
\lambda+|A|$.
\end{enumerate}
\item  $(K,\preceq)$ is an \emph{a.e.c. in $\lambda$} if: The
cardinality of every model in $K$ is $\lambda$, and it satisfies
axioms  a,b,d,e of a.e.c. (Definition \ref{1.1}.1), and axiom
\ref{1.1}.1.c for sequences $\langle M_\alpha:\alpha<\delta
\rangle$ with $\delta<\lambda^+$.
\end{enumerate}
\end{definition}

\begin{remark}
\mbox{}
\begin{enumerate} \item If $K$ is a class of models for a
fixed vocabulary, then $(K,\subseteq)$ satisfies axioms b,d,e of
a.e.c. (Definition \ref{1.1}.1). \item Suppose $(K,\preceq)$ is an
a.e.c.. If $(K,\subseteq)$ satisfies axiom \ref{1.1}.1.c, then
$(K,\subseteq)$ is an a.e.c.. \item If $(K,\preceq)$ is an a.e.c.
and $K' \subseteq K$ then $(K',\preceq \restriction K')$ satisfies
axioms b,d,e of a.e.c. (Definition \ref{1.1}.1).
\end{enumerate}
\end{remark}

We give some simple examples of a.e.c.s. One can see more examples
in \cite{gr21}.

\begin{example}
Let $T$ be a first order theory. Denote $K=:\{M:M \models T\}$.
Define $M\preceq N$ if $M$ is an elementary submodel of $N$.
$(K,\preceq)$ is an a.e.c..
\end{example}

\begin{example}
Let $T$ be a first order theory with $\Pi_2$ axioms, namely axioms
of the form $\forall x \exists y \phi(x,y)$ [For example $(\forall
x,y) (x+y=y+x)$ is OK, as it is equivalent to the $\Pi_2$ axiom
$(\forall x,y) \exists z (x+y=y+x)$]. Denote $K=:\{M:M \models
T\}$. Then $(K,\subseteq)$ is an a.e.c..
\end{example}

\begin{example}
The class of \emph{locally-finite groups} (the subgroup generated
by every finite subset of the group is finite) with the relation
$\subseteq$ is an a.e.c..
\end{example}

\begin{example}
Let $K$ be the class of groups. Let $\preceq=:\{(M,N):M,N$ are
groups, and $M$ is a \emph{pure} subgroup of $N\}$ ($M$ is a pure
subgroup of $N$ if and only if $N \models (\exists y)ry=m$ implies
$M \models (\exists y)ry=m$ for every integer $r$ and every $m \in
M$). $(K,\preceq)$ is an a.e.c..
\end{example}

\begin{example}
The class of models that are isomorphic to $(\mathbb{N},<)$ with
the relation $\subseteq$ is \emph{not} an a.e.c., as it does not
satisfy axiom \ref{1.1}.1.c: $\bigcup
\{\{-n,-n+1,-n+2..0,1,2...\}:0 \leq n \}$ is isomorphic to
$(\mathbb{Z},<)$ although $\{-n,-n+1,-n+2..0,1,2...\}$ is
isomorphic to $(\mathbb{N},<)$.

But the class of models that are isomorphic to $(\mathbb{N},0,<)$
with the relation $\subseteq$ is an a.e.c., (the relation
$\subseteq$ in this case is actually the equality, and this a.e.c.
has just one model).
\end{example}

\begin{example}
The class of \emph{Banach spaces} with the relation $\subseteq$ is
\emph{not} an a.e.c., as it does not satisfy axiom \ref{1.1}.1.c.
\end{example}

\begin{example}
The class of \emph{sets} (i.e. models without relations or
functions) of cardinality less than $\kappa$, where $\aleph_0 \leq
\kappa$ and the relation is $\subseteq$, is \emph{not} an a.e.c.,
as it does not satisfy axiom \ref{1.1}.1.c.

The class of sets with the relation $\preceq=\{(M,N):M \subseteq
N$ and $||N-M||>\kappa\}$ where $\aleph_0 \leq \kappa$, is not an
a.e.c., as it does \emph{not} satisfy smoothness (axiom
\ref{1.1}.1.d).
\end{example}

\begin{definition}
We say $M\prec N$ when $M \preceq N$ and $M \neq N$.
\end{definition}

\begin{definition}\label {1.2}
\mbox{} 
$K_\lambda=:\{M \in K:||M||=\lambda\},\ K_{<\lambda}=\{M\in
K:||M||<\lambda\}$, etc.
\end{definition}

By the following proposition we can replace the increasing
continuous sequence in axioms c,d in Definition \ref{1.1} by a
directed order.
\begin{proposition}
Let $(K,\preceq)$ be an a.e.c., $I$ be a \emph{directed} order and
suppose that for $s,t \in I$ we have $M_s \in K$ and $s \leq_I t
\Rightarrow M_s \preceq M_t$. Then:
\begin{enumerate}
\item $M_0 \preceq \bigcup \{M_s:s \in I \} \in K.$ \item If for
every $s \in I,\ M_s \preceq N \in K$, then $\bigcup \{M_s:s \in I
\} \preceq N.$
\end{enumerate}
\end{proposition}

\begin{proof}
We prove the two items of the proposition simultaneously, by
induction on $|I|$. For finite $I$, there is nothing to prove, so
assume $I$ is infinite. There is an increasing continuous sequence
of subsets of $I$, $\langle I_\alpha:\alpha<|I| \rangle$, such
that $|I_\alpha|<|I|$. Denote $M_{I_\alpha}:=\bigcup\{M_s:s \in
I_\alpha \}$ and $M_I:=\bigcup \{M_s:s \in I \}$. If
$\alpha<\beta<|I|$ then by item (1) of the induction hypothesis,
$s \in I_\alpha \Rightarrow M_s \preceq M_{I_\alpha}$. But as
$I_\alpha \subseteq I_\beta$, $s \in I_\beta$, so $M_s \preceq
M_{I_\beta}$. So by item (2) of the induction hypothesis,
$M_{I_\alpha} \preceq M_{I_\beta}$. Hence the sequence $\langle
M_{I_\alpha}:\alpha<|I| \rangle$ is increasing. But it is also
continuous, as the sequence $\langle I_\alpha:\alpha<|I| \rangle$
is continuous. So by axiom c of Definition \ref{1.1} $M_{I_\alpha}
\preceq M_I \in K$. So as $\preceq$ is transitive and $M_s \preceq
M_{I_\alpha}$ for $s \in I_\alpha$, we have $M_s \preceq M_I \in
K$. Hence we have proved item (1) of the proposition for the
cardinality $|I|$. Now we prove item (2) of the proposition for
$|I|$. If for every $s \in I,\ M_s \preceq N \in K$, then by item
(2) of the induction hypothesis, for $\alpha<|I|$, we have
$M_{I_\alpha} \preceq N \in K$, hence we can apply axiom (d) of
Definition \ref{1.1} for the increasing continuous sequence
$\langle M_{I_\alpha}:\alpha<|I| \rangle$, so $\bigcup
\{M_{I_\alpha}:\alpha<|I| \} \preceq N$. But $M_I=\bigcup
\{M_{I_\alpha}:\alpha<|I| \}$.
\end{proof}

\begin{definition} \label{up}
$(K,\preceq)^{up}:=(K^{up},\preceq^{up})$ where we define:
\begin{enumerate}
\item $K^{up}$ is the class of models with the vocabulary of $K$,
such that there are a directed order $I$, and a set of models
$\{M_s:s \in I\}$ such that: $M=\bigcup \{M_s:s \in I \}$ and $s
\leq_I t \Rightarrow M_s \preceq M_t$. \item For $M,N \in K^{up},\
M \preceq^{up}N$ iff there are directed orders $I,J$ and sets of
models $\{M_s:s \in I \},\ \{N_t:t \in J \}$ respectively such
that: $M=\bigcup \{M_s:s \in I \},\ N=\bigcup \{N_t:t \in J \},\ I
\subseteq J,\ s \leq_J t \Rightarrow N_s \preceq N_t,\ s \leq_I t
\Rightarrow M_s \preceq M_t \preceq N_t$.
\end{enumerate}
\end{definition}

\begin{proposition}
If
\begin{enumerate}
\item $(K_1,\preceq_1),(K_2,\preceq_2)$ are a.e.c.s in $\lambda$.
\item $K_1 \subseteq K_2$. \item $\preceq_2 \restriction K_1$ is
$\preceq_1$.
\end{enumerate}
Then $K_1^{up} \subseteq K_2^{up}$ and $(\preceq_2)^{up}
\restriction K_1^{up}$ is $(\preceq_1)^{up}$.
\end{proposition}

\begin{proof}
Easy.
\end{proof}

\begin{fact}[Lemma 1.23 in \cite{sh600}] \label{1.14}
Let $(K,\preceq)$ be an a.e.c. in $\lambda$. Then
\begin{enumerate}
\item $(K,\preceq)^{up}$ is an a.e.c.. \item $(K^{up})_\lambda=K$.
\item $\preceq^{up} \restriction K$ is $\preceq$. \item
$LST(K,\preceq)^{up}=\lambda$.
\end{enumerate}
\end{fact}

\begin{definition}
\label{1.2b} \mbox{}
\begin{enumerate}
\item Let $M,N$ be models in $K$, $f$ is an injection of $M$ to
$N$. We say that $f$ is a $\preceq$-\emph{embedding} and write
$f:M \to N$, or shortly $f$ is an \emph{embedding} (if $\preceq$
is clear from the context), when $f$ is an injection with domain
$M$ and $Im(f)\preceq N$. \item A function $f:B \to C$ is
\emph{over} $A$, if $A \subseteq B \bigcap C$ and $x \in A
\Rightarrow f(x)=x$.
\end{enumerate}
\end{definition}

\begin{definition}\label{K^3}
\mbox{}
\begin{enumerate}
\item $K^3_{K,\preceq}=:\{(M,N,a):M,N \in K,\ M \preceq N,\ a \in
N\}$. When the class $(K,\preceq)$ is clear from the context we
omit it writing $K^3$. \item $K^3_\lambda:=\{(M,N,a):M,N \in
K_\lambda,\ M \preceq N,\ a \in N\}$.
\end{enumerate}
\end{definition}

\begin{definition}
\label{1.3} \mbox{}
\begin{enumerate}
\item $E^*_{K,\preceq}$ is the following relation on
$K^3_{K,\preceq}$: $(M_0,N_0,a_0)E^* \allowbreak (M_1,N_1,a_1)$
iff $M_1=M_0$ and for some $N_2 \in K_\lambda$ with $N_1 \preceq
N_2$ there is an embedding $f:N_0 \to N_2$ over $M_0$ with
$f(a_0)=a_1$. \item $E_{K,\preceq}$ is the closure of
$E^*_{K,\preceq}$ under transitivity, i.e. the closure to an
equivalence relation.
\end{enumerate}
When $(K,\preceq)$ is clear from the context we omit it writing
$E^*,E$.
\end{definition}

\begin{definition}
\label{1.4} \mbox{}
\begin{enumerate} \item We say that $(K_\lambda,\preceq \restriction K_\lambda)$
has \emph{amalgamation} when: For every $M_0,M_1,M_2$ in
$K_\lambda$, such that $n<3 \Rightarrow M_0 \preceq M_n$, there
are $f_1,f_2,M_3$ such that: $f_n:M_n \to M_3$ is an embedding
over $M_0$, i.e. the diagram below commutes. In such a case we say
that $M_3$ is an amalgam of $M_1,M_2$ over $M_0$.

\[ \xymatrix{\ar @{} [dr]
M_1
\ar[r]^{f_1}  &M_3 \\
M_0 \ar[u]^{\id} \ar[r]_{\id} & M_2 \ar[u]_{f_2} }
\] \item

we say that $K_\lambda$ has \emph{joint embedding} when: If
$M_1,M_2 \in K_\lambda$, then there are $f_1,f_2,M_3$ such that
for $n=1,2$ $f_n:M_n \to M_3$ is an embedding and $M_3 \in
K_\lambda$. \item A model $M$ in $K_\lambda$ is \emph{superlimit}
when:
\begin{enumerate}
\item If $\langle M_\alpha:\alpha \leq \delta \rangle$ is an
increasing continuous sequence of models in $K_\lambda$,
$\delta<\lambda^+$ and $\alpha<\delta \Rightarrow M_\alpha \cong
M$, then $M_\delta \cong M$. \item $M$ is $\preceq$-universal.
\item $M$ is not $\preceq$-maximal.
\end{enumerate}
\item $M \in K$ is $\preceq$-\emph{maximal} if there is no $N \in
K$ with $M\prec N$.
\end{enumerate}
\end{definition}

\begin{proposition} \label{1.5}
\mbox{}
\begin{enumerate} \item For every $M,N_0,N_1 \in K_\lambda$, $a \in N_0-M$ and $b \in N_1-M$, $(M,N_0,a)E^* \allowbreak (M,N_1,b)$ iff there
is an amalgamation $N,f_0,f_1$ of $N_0,N_1$ over $M$ such that
$f_0(a)=f_1(b)$. \item $E^*$ is a reflexive, symmetric relation.
\item If $(K_\lambda,\preceq \restriction K_\lambda)$ has
amalgamation, then $E^*_\lambda$ is an equivalence relation.
\end{enumerate}
\end{proposition}
\begin{proof}
Easy.
\end{proof}

\begin{definition} \label{definition of a
type} \mbox{}
\begin{enumerate}
\item For every $(M,N,a) \in K^3$ let $tp_{K,\preceq}(a,M,N)$, the
\emph{type} of $a$ in $N$ over $M$, be the equivalence class of
$(M,N,a)$ under $E_{K,\preceq}$  When the class $(K,\preceq)$ is
clear from the context we omit it, writing $tp(a,M,N)$ (In other
texts, it is called `$ga-tp(a/M,N)$'). \item For every $M \in K$,
$S(M):=\{tp(a,M,N):(M,N,a) \in K^3\}$. \item If $p=tp(a,M_1,N)$
and $M_0 \preceq M_1$, then we define $p\restriction
M_0=tp(a,M_0,N)$,
\end{enumerate}

\end{definition}

\begin{remark}
By the definitions of $E,E^*$ it is easy to check that
$p\restriction M_0$ does not depend on the representative of $p$.
\end{remark}

\begin{proposition}
For every $M,N,N^+ \in K$ and $a \in N-M$ with $M \bigcup \{a\}
\subseteq N \preceq N^+$, $tp(a,M,N)=tp(a,M,N^+)$.
\end{proposition}

\begin{proof}
Easy.
\end{proof}

\begin{definition} \label{1.7}
Suppose  $M \preceq N$.
\begin{enumerate}
\item For $p \in S(M)$, we say that $N$ \emph{realizes} $p$ if for
some $a \in N$ $p=tp(a,M,N)$. \item For $P \subseteq S(M)$, we say
that $N$ \emph{realizes} $P$ if $N$ realizes every type in P.
\item For $p \in S(M)$ and $a \in N-M$, we say that $a$
\emph{realizes} p, when $p=tp(a,M,N)$.
\end{enumerate}
\end{definition}

\begin{proposition} \label{1.8}
Let $M,M_0 \in K_\lambda,\ M_0 \preceq M$. Suppose
$(K_\lambda,\preceq \restriction K_\lambda)$ has amalgamation and
$LST(K,\preceq) \leq \lambda$. Let $P$ be a set of types over
$M_0,\ |P| \leq \lambda$. Then there is a model $N$ in $K_\lambda$
such that $M\preceq N$ and $N$ realizes P.
\end{proposition}
\begin{proof} Easy.
\end{proof}

\begin{definition} \label{1.9}
Let $M,N \in K$. $M$ is said to be \emph{full} over $N$ when $M$
satisfies $S(N)$. $M$ is said to be \emph{saturated} in
$\lambda^+$ over $\lambda$, when for every model $N \in
K_\lambda$, if $N \preceq M$ then $M$ is full over $N$.
\end{definition}

\begin{remark}
This is the reasonable sense of saturated model we can use in our
context, as we do not want to assume anything about
$K_{<\lambda}$, especially not stability and not amalgamation, (so
a saturated model in $\lambda^+$ over $\lambda$ may not be full
over a model $N \in K_{<\lambda},\ N \preceq M$), see the
following example from \cite{bks}.
\end{remark}

\begin{example}
Let $\tau$ contain infinitely many unary predicates $P_n$ and one
binary predicate $E$. Define a first order theory $T$ such that
$P_{n+1}(x) \Rightarrow P_n(x)$, $E$ is an equivalence relation
with two classes, which are each represented be exactly one point
in $P_n-P_{n+1}$ for each $n$. Now let $K$ be the class of models
in $T$, that omit the type of two inequivalent points that satisfy
all the $P_n$. Then a model $M \in K$ is determined up to
isomorphism by $\mu(M):=|\{x \in M:(\forall n)P_n(x)\}|$. So $K$
is categorical in every uncountable powers, but has $\aleph_0$
countable models (none of them is finite). Now let $\preceq$ be
the relation of being submodel. Then $(K,\preceq)$ is an a.e.c.
with $L.S.T.(K,\preceq)=\aleph_0$. Let $M_0,M_1,M_2 \in K$ be such
that $\mu(M)=0,\mu(M_1)=\mu(M_2)=1$ and $M_1,M_2$ are not
isomorphic over $M_0$. Then there is no amalgamation of $M_1,M_2$
over $M_0$. Now if $\lambda>\aleph_0$ then every model $M \in
K_{\lambda^+}$ is saturated (over $\lambda$). But it is not
saturated over $\aleph_0$, since it can not realize
$tp(a_1,M_0,M_1),tp(a_2,M_0,M_2)$, (where $a_n$ is the unique
element of $M_n-M_0$ of course).
\end{example}

\begin{definition} \label{1.10}
Let $M$ be a model in $K_{\lambda^+}$. $M$ is said to be
\emph{homogenous} in $\lambda^+$ over $\lambda$ if for every
$N_1,N_2 \in K_\lambda$ with $N_1 \preceq M \wedge N_1 \preceq
N_2$, there is a $\preceq$-embedding $f:N_2 \to M$ over $N_1$.
\end{definition}

\begin{definition}
A \emph{representation} of a model $M$ is an $\preceq$-increasing
continuous sequence $\langle M_\alpha:\alpha<||M|| \rangle$ of
models with union $M$, such that $||M_\alpha||<||M||$ for each
$\alpha$ and if $||M||=\lambda^+$ then $||M_\alpha||=\lambda$ for
each $\alpha$.
\end{definition}

The following proposition is a version of Fodor's lemma (there is
no mathematical reason to choose this version, but we think that
it is comfortable).
\begin{proposition} \label{1.11}
There are no $\langle M_\alpha:\alpha \in \lambda^+ \rangle,\
\langle N_\alpha:\alpha \in \lambda^+ \rangle,\ \langle
f_\alpha:\alpha \in \lambda^+ \rangle,S$ such that the following
conditions are satisfied:
\begin{enumerate}
\item The sequences $\langle M_\alpha:\alpha \in \lambda^+
\rangle,\ \langle N_\alpha:\alpha \in \lambda^+ \rangle$ are
$\preceq$-increasing continuous sequences of models in
$K_\lambda$. \item For every $\alpha<\lambda^+$ $f_\alpha:M_\alpha
\to N_\alpha$ is a $\preceq$-embedding. \item $\langle
f_\alpha:\alpha \in \lambda^+ \rangle$ is an increasing continuous
sequence. \item $S$ is a stationary subset of $\lambda^+$. \item
For every $\alpha \in S$,
there is $a \in M_{\alpha+1}-M_\alpha$ 
$M_{\lambda^+}-M_{\alpha}$)
such that $f_{\alpha+1}(a) \in
N_\alpha$.
\end{enumerate}
\end{proposition}

\begin{proof}
Suppose there are such sequences. Denote $M=\bigcup
\{f_\alpha[M_\alpha]:\alpha \in \lambda^+\}$. By clauses 4,5
$||M||=K_{\lambda^+}$. $\langle f_\alpha[M_\alpha]:\alpha \in
\lambda^+ \rangle$, $\langle N_\alpha \bigcap M:\alpha \in
\lambda^+ \rangle$ are representations of $M$. So they are equal
on a club of $\lambda^+$. Hence there is $\alpha \in S$ such that
$f_\alpha[M_\alpha]=N_\alpha \bigcap M$. Hence $f_\alpha[M_\alpha]
\subseteq N_\alpha \bigcap f_{\alpha+1}[M_{\alpha+1}] \subseteq
N_\alpha \bigcap M = f_\alpha[M_\alpha]$ and so this is an
equivalences chain. Especially $f_{\alpha+1}[\allowbreak
M_{\alpha+1}] \bigcap N_\alpha = f_\alpha[M_\alpha]$, in
contradiction to condition 5.
\end{proof}

\begin{proposition} [saturation = model homogeneity]\label{1.12}
Let $(K,\preceq)$ be an a.e.c. such that $K_\lambda$ has
amalgamation, and $LST(K,\preceq) \leq \lambda$. Let $M$ be a
model in $K_{\lambda^+}$. Then $M$ is saturated in $\lambda^+$
over $\lambda$ iff $M$ is a homogenous model in $\lambda^+$ over
$\lambda$.
\end{proposition}
\begin{proof}
One direction is trivial, so let us prove the other direction.
Suppose $M^*_1$ is saturated in $\lambda^+$ over $\lambda$,
${N_0,N_1} \subseteq K_\lambda,\ N_0 \preceq N_1,\ N_0 \preceq
M^*_1$ and there is no embedding of $N_1$ to $M^*_1$ over $N_0$.
Construct by induction on $\alpha \in \lambda^+$ a triple
$(N_{0,\alpha},N_{1,\alpha},f_\alpha)$ such that:
\begin{enumerate}
\item For $n<2$ $\langle N_{n,\alpha}: \alpha \in \lambda^+
\rangle$ is a $\preceq$-increasing continuous sequence of models
in $K_\lambda$. \item $N_{0,0}=N_0,\ N_{1,0}=N_1,\
f_0=id\restriction N_0$. \item For $\alpha \in \lambda^+,\
N_{0,\alpha} \preceq M^*_1$. \item $\langle f_\alpha:\alpha \in
\lambda^+ \rangle$ is an increasing continuous sequence. \item
$f_\alpha:N_{0,\alpha} \to N_{1,\alpha}$ is an embedding. \item
For every $\alpha \in \lambda^+$ there is $a \in
N_{0,\alpha+1}-N_{0,\alpha}$ such that $f_{\alpha+1}(a) \in
N_{1,\alpha}$.
\end{enumerate}

Why can we carry out the construction?\\
for $\alpha=0$ see 2. For $\alpha$ limit, take unions. Suppose we
have chosen $N_{0,\alpha},N_{1,\alpha},f_\alpha$, how will we
choose $N_{0,\alpha+1},\ N_{1,\alpha+1},\ f_{\alpha+1}$?
$f_\alpha[N_{0,\alpha}] \neq N_{1,\alpha}$ (otherwise
$f^{-1}_\alpha\restriction N_1$ is an embedding of $N_1$ to
$M^*_1$ over $N_0$, in contradiction to our assumption). Hence
there is $c \in N_{1,\alpha}-f_\alpha[N_{0,\alpha}]$. As $M^*_1$
is saturated in $\lambda^+$ over $\lambda$, there is $a \in M^*_1$
such that
$tp(a,N_{0,\alpha},M^*_1)=f^{-1}_\alpha(tp(c,f_\alpha[N_{0,\alpha}],N_{1,\alpha})$.
Now $LST(K,\preceq) \leq \lambda$ so there is $N_{0,\alpha+1} \in
K_\lambda$, such that $N_{0,\alpha} \bigcup \{a\} \subseteq
N_{0,\alpha+1} \preceq M^*_1$. So by axiom e of a.e.c.
$N_{0,\alpha} \preceq N_{0,\alpha+1}$. Hence
$f_\alpha(tp(a,N_{0,\alpha},N_{0,\alpha+1}))=tp(c,f_\alpha[N_{0,\alpha}],N_{1,\alpha})$.
By the definition of type and having amalgamation, for some model
$N_{1,\alpha+1}$ and an embedding $f_{1,\alpha+1}$ the following
hold: $N_{1,\alpha} \preceq N_{1,\alpha+1},\ f_{1,\alpha+1}(a)=c$
and $f_\alpha \subseteq f_{\alpha+1}:N_{0,\alpha+1} \to
N_{1,\alpha+1}$. Hence
we can carry out the construction.\\
Now the conditions on the existence of the sequences $\langle
N_{0,\alpha}:\alpha \in \lambda^+ \rangle, \langle
N_{1,\alpha}:\alpha \in \lambda^+ \rangle, \langle f_\alpha:\alpha
\in \lambda^+ \rangle$ contradict Proposition \ref{1.11}
(requirement 5 in Proposition \ref{1.11} is satisfied by
requirement 6 in the construction here).
\end{proof}

Now we prove the uniqueness of the saturated model, although we do
not know its existence.
\begin{theorem}[the uniqueness of the saturated model]\label{1.13}
Suppose $(K_\lambda,\preceq \restriction K_\lambda)$ has the
amalgamation property and $LST(K,\preceq) \leq \lambda$.
\begin{enumerate}
\item Let $N \in K_\lambda$ and for $n=1,2$ $N \preceq M_n$ and
$M_n$ is saturated in $\lambda^+$ over $\lambda$. Then $M_1,\ M_2$
are isomorphic over $N$. \item If $M_1,\ M_2$ are saturated in
$\lambda^+$ over $\lambda$ and $(K_\lambda,\preceq \restriction
K_\lambda)$ has the joint embedding property then $M_1,\ M_2$ are
isomorphic.
\end{enumerate}
\end{theorem}

\begin{proof}
(1) We use the hence and forth method. For $n=1,2$ Let $\langle
a_{n,\alpha}:\alpha \in \lambda^+ \rangle$ be an enumeration of
$M_n$ without repetitions. We choose by induction on $\alpha \in
\lambda^+$ a triple $(N_{1,\alpha},N_{2,\alpha},f_\alpha)$ such
that:
\begin{enumerate}[(a)]
\item $N_{n,0}=N,\ f_0=id$. \item $N_{n,\alpha} \preceq M_n$.
\item The sequence $\langle N_{n,\alpha}:\alpha \in \lambda^+
\rangle$ is an increasing continuous sequence of models in
$K_\lambda$. \item $\langle f_\alpha:\alpha \in \lambda^+ \rangle$
is increasing and continuous. \item $f_\alpha:N_{1,\alpha} \to
N_{2,\alpha}$ is an embedding. \item $a_{n,\alpha} \in
N_{n,2\alpha+n}$.
\end{enumerate}
Why can we carry out the construction?\\
For $\alpha=0$ see a. Let $\alpha$ be a limit ordinal. For $n=1,2$
Define $N_{n,\alpha}= \bigcup \{N_{n,\beta}: \beta<\alpha \},\
f_\alpha=\bigcup \{f_\beta:\beta<\alpha \}$. By axiom c of a.e.c.
(i.e. the closure under increasing continuous sequences) for
$n=1,2$ $\beta<\alpha \Rightarrow N_{n,\beta} \preceq
N_{n,\alpha}$ and by axiom \ref{1.1}.1.d (smoothness)
$N_{n,\alpha} \preceq M_n$. So there is no problem in the limit
case. Suppose we have defined $N_{1,\alpha},\ N_{2,\alpha},\
f_\alpha$. Suppose $\alpha=2\beta$. As $LST(K,\preceq) \leq
\lambda$, there is a model $N_{1,\alpha+1} \in K_\lambda$ such
that $N_{1,\alpha} \bigcup\{a_{1,\beta}\} \subseteq N_{1,\alpha+1}
\preceq M_1$. By the induction hypothesis (b) $N_{1,\alpha}
\preceq M_1$. Now by axiom \ref{1.1}.1.c (closure under increasing
continuous sequences) $N_{1,\alpha} \preceq N_{1,\alpha+1}$. Let
$f^+_\alpha$ be an injection with domain $N_{1,\alpha+1}$ such
that $f_\alpha \subseteq f_{1,\alpha+1}$. Actually it is an
isomorphism of its domain to its range. The relation $\preceq$ is
closed under isomorphisms, so $N_{2,\alpha}=f_\alpha[N_{1,\alpha}]
\preceq f^+_\alpha[N_{1,\alpha+1}]$. $M_2$ is saturated in
$\lambda^+$ over $\lambda$ and so by Lemma \ref{1.12} it is model
homogenous in $\lambda^+$ over $\lambda$. So there is an embedding
$g:f^+_\alpha[N_{1,\alpha+1}] \to M_2$ over $N_{2,\alpha}$. Define
$f_{\alpha+1}=:g \circ f^+_\alpha,\
N_{2,\alpha+1}=:f_{\alpha+1}[N_{1,\alpha+1}]$. $f_\alpha \subseteq
f_{\alpha+1}$ and so (d) is satisfied. Requirement a is not
relevant for the successor case. (b) is satisfied for n=1 by the
definition of $N_{n,\alpha+1}$ and for n=2 as g is
$\preceq-$embedding. (c) is satisfied for n=1 by the construction
and for n=2 as $\preceq$ respects isomorphisms. (e) is satisfied
by the definition of $f_{\alpha+1}$. (f) is relevant only for n=1.
Hence we can carry out the construction in the $\alpha+1$ step for
$\alpha$ even. The case $\alpha$ is an odd number is symmetric, so
we have to change $a,b$.
Hence we can carry out the construction.\\
Now by (b),(f) $\bigcup \{N_{n,\alpha}:\alpha \in \lambda^+
\}=M_n$ . Define $f=\bigcup \{f_\alpha:\alpha \in \lambda^+ \}$.
By (e) $f:M_1 \to M_2$ is an isomorphism. By (a),(d) this
isomorphism is over $N$.\\
(2) For $n=1,2$ As $LST(K,\preceq) \leq \lambda$ there is $N_n
\preceq M_n$ in $K_\lambda$. $K_\lambda$ has the joint embedding
property and so there is a model $N$ and embeddings $f_n:N_n \to
N$. Let $f^+_n$ an injection with domain $M_n$ such that $f_n
\subseteq f^+_n$. By Lemma \ref{1.12} for $n=1,2$ there is an
embedding $g_N:N \to f^+_n[M_n]$ over $f_n[N_n]$. Now $f=g_1 \circ
g^{-1}_2$ is an isomorphism and so there is an injection $g^+$
with domain $f^+_2[M_2]$ such that $g \subseteq g^+$. By the
definition of $g_2,\ g_2[N] \preceq f^+_2[M_2]$ and so as
$\preceq$ respects isomorphisms, $g_1[N]=g[g_2[N]] \preceq
g^+[f_2[M_2]]$. By item a $f^+_1[M_1],\ g^+[f^+_2[M_2]]$ are
isomorphic over $g_1[N]$. Hence $M_1,M_2$ are isomorphic.
\end{proof}
\section{Non-forking frames}
\discussion{The plan} Suppose we know something about $K_\lambda$,
especially that there is no $\preceq$-maximal model. Can we say
something about $K_{\lambda^{+n}}$? At least we want to prove that
$K_{\lambda^{+n}} \neq \emptyset$. It is easy to prove that
$K_{\lambda^+} \neq \emptyset$ [How? We choose $M_\alpha$ by
induction on $\alpha<\lambda^+$ such that $M_\alpha \prec
M_{\alpha+1}$ and if $\alpha$ is limit we define
$M_\alpha:=\bigcup \{M_\beta:\beta<\alpha\}$ (by axiom
\ref{1.1}.1.c $M_\alpha \in K$). In the end $M_{\lambda^+} \in
K_{\lambda^+}$]. What about $K_{\lambda^{+2}}$? The main topic in
this paper is semi-good frames. If there is a semi-good
$\lambda$-frame, then by Proposition \ref{3.3}.2 there is no
$\preceq$-maximal model in $K_{\lambda^+}$. So $K_{\lambda^{++}}
\neq \emptyset$. Moreover, Theorem \ref{11.1}.1 says that if
$\frak{s}$ is a semi-good $\lambda$-frame with some additional
assumptions and $\lambda$ satisfies specific set-theoretic
assumptions, then there is a good $\lambda^+$-frame
$\frak{s}^+=(K^+,\preceq^+,S^{bs,+},\dnf^+)$, such that $K^+
\subseteq K$ and the relation $\preceq ^+ \restriction K^+$ is
included in the relation $\preceq \restriction K^+$. So
$K_{\lambda^{+3}} \neq \emptyset$ and so on. Thus we prove
$K_{\lambda^{+n}} \neq \emptyset$ by induction on $n<\omega$,
(assuming reasonable assumptions).

Definition \ref{definition good frame} is an axiomatization of the
non-forking relation in superstable elementary class. If we
subtract axiom \ref{definition good frame}.3.c., we get the basic
properties of the non-forking relation in $(K_\lambda,\preceq
\restriction K_\lambda)$ where $(K,\preceq)$ is stable in
$\lambda$.

Sometimes we do not find a natural independence relation on all
the types. So first we extend the notion of an a.e.c. in $\lambda$
by adding a new function $S^{bs}$ which assigns a collection of
basic types to each model in $K_\lambda$, and then we add an
independence relation $\dnf$ on basic types.

It is reasonable to assume \emph{categoricity} in some cardinality
$\lambda$ for some reasons:
\begin{enumerate}
\item In Example \ref{the main example} $K$ is categorical in
$\lambda$. \item If $K$ is not categorical in any cardinality,
then we know $\{\lambda:K$ is categorical in $\lambda\}$, it is
the empty set.
 \item If there is a superlimit model in $K_\lambda$,
then we can reduce $(K_\lambda, \preceq \restriction K_\lambda)$
to the models which are isomorphic to it, and therefore obtain
categoricity in $\lambda$ (see section 1 in \cite{sh600}). However
this case requires stability.
\end{enumerate}

We do not assume \emph{amalgamation}, but we assume amalgamation
in $(K_\lambda, \preceq \restriction K_\lambda)$. This is a
reasonable assumption because it is proved in \cite{sh88r} that if
an a.e.c. is categorical in $\lambda$ and amalgamation fails in
$\lambda$ then under plausible set theoretic assumptions there are
$2^{\lambda^+}$ models in $K_{\lambda^+}$.

\begin{definition} \label{2.1a} \label{definition good frame}
\label{definition of a good frame}
$\frak{s}=(K,\preceq,S^{bs},\dnf)$ is a \emph{good $\lambda$-frame} if:\\
(1)
\begin{enumerate}[(a)]
\item $(K,\preceq)$ is an a.e.c.. \item
$LST(K,\preceq) \leq \lambda$. \item $(K_\lambda, \preceq \restriction K_\lambda)$ 
has joint
embedding. \item $(K_\lambda, \preceq \restriction K_\lambda)$ has amalgamation. \item There is no $\preceq$-maximal model in $K_\lambda$.\\
\end{enumerate}
(2) $S^{bs}$ is a function with domain $K_\lambda$, which
satisfies the following axioms:
\begin{enumerate}[(a)]
\item $S^{bs}(M) \subseteq S^{na}(M)=:\{tp(a,M,N):M\prec N \in
K_\lambda,\ a\in N-M\}$. \item It respects isomorphisms. \item
Density of the basic types: If $M,N \in K_\lambda$ and $M \prec
N$, then there is $a \in N-M$ such that $tp(a,M,N) \in S^{bs}(M)$.
\item Basic stability: For every $M \in K_\lambda$, the
cardinality of $S^{bs}(M)$ is $\leq \lambda$.
\end{enumerate}
(3) the relation $\dnf$ satisfies the following axioms:
\begin{enumerate}[(a)]
\item $\dnf$ is a subset of $\{(M_0,M_1,a,M_3):M_0,M_1,M_3 \in K,\
||M_0||=||M_1||=\lambda,\ M_0 \preceq M_1 \preceq M_3 ,\ a \in
M_3-M_1\ and \ n<2 \Rightarrow tp(a,M_n,M_3) \in S^{bs}(M_n)\}$.
\item Monotonicity: If $M_0 \preceq M^*_0 \preceq M^*_1 \preceq
M_1 \preceq M_3 \preceq M_3^*,\ M^*_1\bigcup \{a\} \subseteq
M^{**}_3 \preceq M^*_3$, then $\dnf(M_0,M_1,a,M_3) \Rightarrow
\dnf(M^*_0,M^*_1,a,M^{**}_3)$. See Remark \ref{dnf not as a sign}.
\item Local character: For every limit ordinal $\delta<\lambda^+$
if $\langle M_\alpha:\alpha \leq \delta \rangle$ is an increasing
continuous sequence of models in $K_\lambda$, and
$tp(a,M_\delta,M_{\delta+1}) \in S^{bs}(M_\delta)$, then there is
$\alpha<\delta$ such that $tp(a,M_\delta,\allowbreak
M_{\delta+1})$ does not fork over $M_\alpha$. \item Uniqueness of
the non-forking extension: If $p,q \in S^{bs}(N)$ do not fork over
$M$, and $p\restriction M=q\restriction M$, then $p=q$. \item
Symmetry: If $M_0 \preceq M_1 \preceq M_3,\ a_1 \in M_1,\
tp(a_1,M_0,M_3) \in S^{bs}(M_0)$, and $tp(a_2,M_1,M_3)$ does not
fork over $M_0$, \emph{then} there are $M_2,M^*_3$ such that $a_2
\in M_2,\ M_0 \preceq M_2 \preceq M^*_3,\ M_3 \preceq M^*_3$, and
$tp(a_1,M_2,M^*_3)$ does not fork over $M_0$. \item Existence of
non-forking extension: If $p \in S^{bs}(M)$ and $M\prec N$, then
there is a type $q \in S^{bs}(N)$ such that $q$ does not fork over
$M$ and $q\restriction M=p$. \item Continuity: Let $\langle
M_\alpha:\alpha \leq \delta \rangle$ be an increasing continuous
sequence. Let $p \in S(M_\delta)$. If for every $\alpha \in
\delta,\ p\restriction M_\alpha$ does not fork over $M_0$, then $p
\in S^{bs}(M_\delta)$ and does not fork over $M_0$.
\end{enumerate}
\end{definition}

\begin{remark}\label{dnf not as a sign}
If $\dnf(M_0,M_1,a,M_3)$ and $tp(b,M_1,M_3^{**})=tp(a,M_1,M_3)$
then by Definition \ref{definition good frame}.3.b (the
monotonicity axiom) $\dnf(M_0,M_1,a,M_3)$. Therefore we can say
``$p$ does not fork over $M_0$'' instead of $\dnf(M_0,M_1,a,M_3)$.
\end{remark}

While in \cite{sh600} we study good frames, so basic stability is
assumed, here we assume basic \emph{weak} stability so the
following definition is central:
\begin{definition} \label{2.1b}
$\frak{s}=((K^\frak{s},\preceq^\frak{s},S^{bs,\frak{s}},\dnf^\frak{s})=(K,\preceq,S^{bs},\dnf)$
is a \emph{semi}-good $\lambda$-frame, if $\frak{s}$ satisfies the
axioms of a good $\lambda$-frame except that instead of assuming
basic stability, we assume that $\frak{s}$ has basic weak
stability, namely for every $M \in K_\lambda$ $S^{bs}(M)$ has
cardinality at most $\lambda^+$.

$\frak{s}$ is said to be a \emph{semi-good frame} if it is a
semi-good $\lambda$-frame for some $\lambda$.
\end{definition}

\begin{remark}\label{local character alomost implies continuity}
If for each $M \in K_\lambda$ $S^{bs}(M)=\{tp(a,M,N):M \prec N, a
\in N-M\}$ then the continuity axiom is an easy consequence of the
local character.
\end{remark}

Can we define in our context independence, orthogonality and more
things like in superstable theories? The answer is: See
\cite{sh705} (mainly sections 5,6) and \cite{jrsi3}.

Now we give examples of good frames and examples of semi-good
frames.
\begin{example}\label{example superstable regular}
Let $T$ be a superstable first order theory and let $\lambda$ be a
cardinal $\geq |T|+\aleph_0$ such that $T$ is stable in $\lambda$.
Let $K_{T,\lambda}$ be the class of models of $T$ of cardinality
at least $\lambda$. Let $\preceq$ denote the relation of being an
elementary submodel. Let $S^{bs}(M)$ be the set of regular types
over $M$. Let $\dnf$ be as usual. Then by Claim 3.1 on page 52 in
\cite{sh600} (or see \cite{sh91})
$(K_{T,\lambda},\preceq,S^{bs},\dnf)$ is a good $\lambda$-frame.
\end{example}

\begin{example}
The same as Example \ref{example superstable regular} but the
basic types are the non-algebraic types (see \cite{hs89}).
\end{example}

\begin{example}[the main example] \label{the main example}
An example of a semi-good $\lambda$-frame which appears in section
3 of \cite{sh600} and is based on \cite{sh88r}: Let $(K,\preceq)$
be an a.e.c. with a countable vocabulary,
$LST(K,\preceq)=\aleph_0$, $(K,\preceq)$ is $PC_{\aleph_0}$ (i.e.
$K$ is the class of reduced models to a smaller language, of some
countable elementary class, which omit a countable set of types,
and the relation $\preceq$ is defined similarly), it has an
intermediate number of non-isomorphic models of cardinality
$\aleph_1$ and $2^{\aleph_0}<2^{\aleph_1}$. Then we can derive a
semi-good $\aleph_0$-frame from it.

How? For each $M \in K_{\aleph_0}$ define $K_M=\{N \in K:N
\equiv_{L_{\infty,\omega}}M\},\ \preceq_M=\{(N_1,N_2):N_1,N_2 \in
K_M,\ N_1 \preceq N_2,\ and \ N_1\preceq_{L_{\infty,\omega}}N_2
\}$. There is a model $M \in K_{\aleph_0}$ such that
$(K_M)_{\aleph_1} \neq \emptyset$. Fix such an $M$. For every $N
\in K_M$ with $||N||=\aleph_0$ define
$S^{bs}(N)=\{tp(a,N,N^*):N\prec_M N^* \in K_M,\ a \in N^*-N \}$.
Define $\dnf:=\{(M_0,M_1,a,M_2):M_0,M_1,M_2 \in K_M,\
||M_0||=||M_1||=\aleph_0,\ M_0 \preceq_M M_1 \prec_M M_2$ and
$tp(a,M_1,M_2)$ is definable over some finite subset $A$ of $M_0$
in the following sense: For every type $q \in S^{bs}(M_1)$ if `$q
\restriction A=tp(a,A,M_2)$' [more precisely for some
$b,M_2^*,M_3,f$, $M_1 \preceq M_2^*$, $tp(b,M_1,M_2^*)=q$ $M_2
\preceq M_3$ and $f$ is an isomorphism of $M_2^*$ to $M_3$ over
$A$ with $f(b)=a$] then $q=tp(a,M_1,M_2)$.
By the proof of Theorem 3.4 on page 54 in \cite{sh600}
$\frak{s}:=(K_M,\preceq_M,S^{bs},\dnf)$ is a semi-good
$\aleph_0$-frame [Why? We assumed here assumptions
$(\alpha),(\beta),(\gamma)$ of Theorem 3.4. So by Theorem 3.4.1
for some $M \in K_{\aleph_0}$ we have $(\delta^-),(\epsilon)$ too.
So if $(\delta)$ (namely stability) holds then by Theorem 3.4.2
$\frak{s}$ is a good $\aleph_0$-frame. Here we have $(\delta^-)$
(namely weak stability) only. In the beginning of the proof of
item 2 (on page 56) it is written `we assume $(\delta^-)$ instead
of $(\delta)$'. By the continuation of the proof, we see that
$\frak{s}$ is a semi-good $\aleph_0$-frame]
\end{example}

\begin{definition}\label{2.2}
Let $p_0 \in S(M_0), p_1 \in S(M_1)$. We say that $p_0,p_1$
\emph{conjugate} if for some $a_0,M_0^+,a_1,M_1^+,f$ the following
hold:
\begin{enumerate}
\item For $n=0,1$, $tp(a_n,M_n,M_n^+)=p_n$. \item $f:M_0^+ \to
M_1^+$ is an isomorphism. \item $f \restriction M_0:M_0 \to M_1$
is an isomorphism. \item $f(a_0)=a_1$.
\end{enumerate}
\end{definition}

\begin{proposition}\label{2.3}
Let $p_0 \in S(M_0), p_1 \in S(M_1)$ and assume that $p_0,p_1$
conjugate.
\begin{enumerate}
\item If for $n=1,2$, $tp(a_n,M_n,M_n^+)=p_n$ and there is an
embedding $f:M_0^+ \to M_1^+$ such that $f \restriction M_0:M_0
\to M_1$ is an isomorphism and $f(a_0)=a_1$, then the types
$p_0,p_1$ conjugate. \item Suppose that $a_0,M_0^+,a_1,M_1^+,f$
are as in Definition \ref{2.2}, i.e. they witness that $p_0,p_1$
conjugate. If $tp(a_1^*,M_1,M_1^*)=tp(a_1,M_1,M_1^+)$ then for
some $M_1^{**}$ and $f^*$ such that:
\begin{enumerate}
\item $M_1^* \preceq M_1^{**}$. \item $f^*:M_0^+ \to M_1^{**}$ is
an embedding. \item $a_0,M_0^+,a_1^*,f^*[M_0^+],f^*$ witness that
$p_0,p_1$ conjugate.
\end{enumerate}
\item Assume that ($p_0,p_1$ conjugate and) the types
$p_1,p_2$ conjugate. Then $p_0,p_2$ conjugate.
\end{enumerate}
\end{proposition}

\begin{proof}
\mbox{}
\begin{enumerate}
\item Substitute $f[M_0^+]$ instead of $M_1^+$ in Definition
\ref{2.2}. \item Since $tp(a_1^*,M_1,M_1^*)=tp(a_1,M_1,M_1^+)$,
there is an amalgamation $(id_{M_1^*},g,M_1^{**})$ of
$M_1^*,M_1^+$ over $M_1$ such that $g(a_1)=a_1^*$. Now define
$f^*:=g \circ f$. So $f^*(a_0)=g(f(a_0))=g(a_1)=a_1^*$.
\item Suppose that $a_1^*,M_1^*,a_2,M_2^*,g$ witness that
$p_1,p_2$ conjugate. Since the types $p_0,p_1$ conjugate, there
are witnesses for this. So by 
Proposition \ref{2.3}.2 for some $a_0,M_0^{**},f^*$ the following
hold:
\begin{displaymath}
\xymatrix{& M_1^{**} \ar[r]^{g^+} & M_2^+\\
M_0^+ \ar[ru]^{f^*} \ar[r]^{id} & M_1^* \ar[u]^{id} \ar[r]^{g} &
M_2^* \ar[u] ^{id}\\
M_0 \ar[u]^{id} \ar[r]^{f^* \restriction M_0} & M_1 \ar[u]^{id}
\ar[r]^{g \restriction M_1} & M_2 \ar[u]^{id}}
\end{displaymath}
\begin{enumerate}
\item $M_1^* \preceq M_1^{**}$. \item $f^*:M_0^+ \to M_1^{**}$ is
an embedding. \item $a_0,M_0^+,a_1^*,f^*[M_0^+],f^*$ witness that
$p_0,p_1$ conjugate.
\end{enumerate}
Since $(K_\lambda, \preceq \restriction K_\lambda)$ has
amalgamation, for some $g^+,M_2^+$ the following hold:
\begin{enumerate}
\item $g^+:M_1^{**} \to M_2^+$ is an embedding. \item $M_2^*
\preceq M_2^+$. \item $g^+ \restriction M_1^*=g$.
\end{enumerate}
Now define $f:=g^+ \circ f^*$. So $a_0,M_0^+,a_2,f[M_0^+],f$
witness that the types $p_0,p_2$ conjugate. Why? We verify the
conditions in Definition \ref{2.2}:
\begin{enumerate}[(1)]
\item $tp(a_0,M_0,M_0^+)=p_0$ because
$a_0,M_0^+,a_1^*,f^*[M_0^+],f^*$ witness that $p_0,p_1$ conjugate.
$tp(a_2,M_2,M_2^+)=p_2$ because $a_1^*,M_1^*,a_2,M_2^*,g$ witness
that $p_1,p_2$ conjugate. \item $f:M_0^+ \to f^*[M_0^+]$ is an
isomorphism. \item  $f \restriction M_0:M_0 \to M_2$ is or course
an embedding, but why is it onto? Take $z \in M_2$. Since $g
\restriction M_1$ is an isomorphism (as $a_1^*,M_1^*,a_2,M_2^*,g$
witness that $p_1,p_2$ conjugate), there is $y \in M_1$ such that
$g(y)=z$. Since $f^* \restriction M_0$ is an isomorphism (as
$a_0,M_0^+,a_1^*,f^*[M_0^+],f^*$ witness that $p_0,p_1$
conjugate), there is $x \in M_0$ such that $g(x)=z$. Therefore
$g^+(f^*(x))=g(f^*(x))=g(y)=z$. \item
$f(a_0)=g^+(f^*(a_0))=g^+(a_1^*)=a_2$.
\end{enumerate}
\end{enumerate}
\end{proof}

\begin{definition}\label{f(p)}
Let $p=tp(a,M,N)$. Let $f$ be a bijection with domain $M$. Define
$f(p)=tp(f(a),f[M],f^+[N])$, where $f^+$ is an extension of $f$
(and the relations and functions on $f^+[N]$ are defined such that
$f^+:N \to f^+[N]$ is an isomorphism).
\end{definition}

\begin{proposition}
The definition of $f(p)$ in Definition \ref{f(p)} does not depend
on the representative $(M,N,a) \in p$.
\end{proposition}

\begin{proof}
By Proposition \ref{2.3}.2.
\end{proof}



\begin{definition} \label{2.4}
Let $\frak{s}$ be a semi-good $\lambda$-frame. We say that
$\frak{s}$ \emph{has conjugation} when: $K_\lambda$ is categorical
and if $M_1,M_2 \in K_\lambda$, $M_1 \preceq M_2$ and $p_2 \in
S^{bs}(M_2)$ is the non-forking extension of $p_1 \in
S^{bs}(M_1)$, then the types $p_1,p_2$ conjugate.
\end{definition}

\begin{proposition}
The semi-good frame
in Example \ref{the main example} has conjugation.
\end{proposition}

\begin{proof}
Assume that $||M_0||=||M_1||=\aleph_0$, $M_0 \prec_M M_1$ and $p
\in S^{bs}(M_1)$ is definable over some finite subset $A$ of
$M_0$. We have to prove that the types $p,p_0:=p \restriction M_0$
conjugate. Since $M_0 \prec_{L_{\infty,\omega}}M_1$, there is an
isomorphism $f:M_0 \to M_1$ over $A$. So $f(p)$ does not fork over
$f[M_0]$. But $p$ does not fork over $M_0$ too. $p \restriction
A=(p \restriction M_0) \restriction A=p_0 \restriction A=f(p_0)
\restriction A$. Since $f(p_0),p$ are definable over $A$,
$f(p_0)=p$.
\end{proof}

\begin{proposition}[versions of extension]
\label{a version of extension}  If for $n<3$ $M_n \in K_\lambda$,
$M_0 \preceq M_n$, and $tp(a,M_1,M_0) \in S^{bs}(M_0)$ then:
\begin{enumerate}
\item There are $M_3,f$ such that:
\begin{enumerate}
\item $M_2 \preceq M_3$. \item $f:M_1 \to M_3$ is an embedding
over $M_0$. \item $tp(f(a),M_2,M_3)$ does not fork over $M_0$.
\end{enumerate}
\item There are $M_3,f$ such that:
\begin{enumerate}
\item $M_1 \preceq M_3$. \item $f:M_2 \to M_3$ is an embedding
over $M_0$. \item $tp(a,f[M_2],M_3)$ does not fork over $M_0$.
\end{enumerate}
\end{enumerate}
\end{proposition}
\begin{proof}
Easily by Definition \ref{2.1a}.2.
\end{proof}

\begin{proposition}[The transitivity proposition] \label{2.7} \label{s is
transitive}\label{transitivity}
Suppose $\frak{s}$ is a semi-good $\lambda$-frame. Then: If $M_0
\preceq M_1 \preceq M_2,\ p \in S^{bs}(M_2)$ does not fork over
$M_1$ and $p\restriction M_1$ does not fork over $M_0$, then $p$
does not fork over $M_0$.
\end{proposition}
\begin{proof} Suppose $M_0\prec M_1\prec M_2,\ n<3 \Rightarrow M_n \in
K_\lambda,\ p_2 \in S^{bs}(M_2)$ does not fork over $M_1$ and
$p_2\restriction M_1$ does not fork over $M_0$. For $n<2$ define
$p_n=p_2\restriction M_n$. By axiom f there is a type $q_2 \in
S^{bs}(M_2)$ such that $q_2\restriction M_0=p_0$ and $q_2$ does
not fork over $M_0$. Define $q_1=q_2\restriction M_1$. By axiom b
(monotonicity) $q_1$ does not fork over $M_0$. So by axiom d
(uniqueness) $q_1=p_1$. Using again axiom e, we get $q_2=p_2$, as
they do not fork over $M_1$. By the definition of $q_2$ it does
not fork over $M_0$.
\end{proof}

\begin{proposition} \label{axiom i is redundant}
Suppose
\begin{enumerate}
\item $\frak{s}$ satisfies the axioms of a semi-good
$\lambda$-frame. \item $n<3 \Rightarrow M_0 \preceq M_n$. \item
For $n=1,2$, $a_n \in M_n-M_0\ and \ tp(a_n,M_0,M_n) \allowbreak
\in S^{bs}(M_0)$.
\end{enumerate}
Then there is an amalgamation $M_3,f_1,f_2$ of $M_1,M_2$ over
$M_0$ such that for $n=1,2$ $tp(f_n(a_n),f_{3-n}[M_{3-n}],M_3)$
does not fork over $M_0$.
 \end{proposition}
\begin{proof} Suppose for $n=1,2$ $M_0\prec M_n \wedge tp(a_n,M_0,M_n)
\in S^{bs}(M_0)$. By Proposition \ref {a version of extension}.1
there are $N_1,f_1$ such that:
\begin{enumerate}
\item $M_1 \preceq N_1$. \item $f_1:M_2 \to N_1$ is an embedding
over $M_0$. \item $tp(f_1(a_2),M_1,N_1)$ does not fork over $M_0$.
\end{enumerate}

\begin{displaymath}
\xymatrix{M_2 \ar[r]^{f_1} & N_1 \ar[rr]^{id} && N_2 \ar[r]^{id}
\ar[d]^{f_2} & N_3\\
&& N_2^* \ar[r]^{id} \ar[ru]^{id} & f_2[N_2] \ar[ru]^{id}\\
M_0 \ar[r]^{id} \ar[uu]^{id} \ar[rru]^{id} & M_1 \ar[uu]^{id} }
\end{displaymath}

By axiom f (the symmetry axiom), there are a model $N_2$, $N_1
\preceq N_2 \in K_\lambda$ and a model $N_2^* \in K_\lambda$ such
that: $M_0 \bigcup \{f_1(a_2)\} \subseteq N_2^* \preceq N_2$ and
$tp(a_1,N_2^*,N_2)$ does not fork over $M_0$.\\
By Proposition \ref{a version of extension}.2 (substituting
$N_2^*,N_2,N_2,a_1$ which appear here instead of $M_0,M_1,M_2,a$
there) there are $N_3,f_2$ such that:
\begin{enumerate}
\item $N_2 \preceq N_3$. \item $f_2:N_2 \to N_3$ is an embedding
over $N_2^*$. \item $tp(a_1,f_2[N_2],N_3)$ does not fork over
$N_2^*$.
\end{enumerate}
So by Proposition \ref{s is transitive} (on page \pageref{s is
transitive}), $tp(a_1,f_2[N_2],N_3)$ does not fork over $M_0$. So
as $M_0 \preceq f_2 \circ f_1[M_2] \preceq f_2[N_2]$ by axiom b
(monotonicity) $tp(a_1,f_2 \circ f_1[M_2],N_3)$ does not fork over
$M_0$. As $f_2 \restriction N_2^*=id_{N_2^*}$,
$f_2(f_1(a_1))=f_1(a_1)$.
\end{proof}

\begin{theorem} \label{4.0}
Suppose $\frak{s}$ satisfies conditions 1 and 2 of a semi-good
$\lambda$-frame (so actually the relation $\dnf$ is not relevant).
\begin{enumerate}
 \item Suppose:
\begin{enumerate}
\item $\langle  M_\alpha:\alpha \leq \lambda^+ \rangle$ is an
increasing continuous sequence of models in $K_\lambda$. \item
There is a stationary set $S \subseteq \lambda^+$ such that for
every $\alpha \in S$ and every model $N$, with $M_\alpha \prec N$
there is a type $p \in S(M_\alpha)$ which is realized in
$M_{\lambda^+}$ and in $N$.
\end{enumerate}
Then $M_{\lambda^+}$ is full over $M_0$ and is saturated in
$\lambda^+$ over $\lambda$. \item Suppose:
\begin{enumerate}
   \item $\langle  M_\alpha:\alpha \leq \lambda^+ \rangle$ is an increasing continuous
sequence of models in $K_\lambda$.
   \item For every $\alpha \in \lambda^+$ and every $p \in
S^{bs}(M_\alpha)$ there is $\beta \in (\alpha,\lambda^+)$ such
that $p$ is realized in $M_\beta$.
\end{enumerate}
Then $M_{\lambda^+}$ is full over $M_0$ and $M_{\lambda^+}$ is
saturated in $\lambda^+$ over $\lambda$. \item There is a model in
$K_{\lambda^+}$ which is saturated in $\lambda^+$ over $\lambda$.
\item $M \in K_\lambda \Rightarrow |S(M)| \leq \lambda^+$.
\end{enumerate}
\end{theorem}

\begin{proof}
Obviously $1 \Rightarrow 2$ and $3 \Rightarrow 4$. To show $2
\Rightarrow 3$, we construct a chain satisfying the hypotheses of
2. Let $cd$ be an injection from $\lambda^+ \times \lambda^+$ onto
$\lambda^+$. Define by induction on $\alpha<\lambda^+$ $M_\alpha$
and $\langle p_{\alpha,\beta}:\beta<\lambda^+ \rangle$ such that:
\begin{enumerate}
\item $\langle  M_\alpha:\alpha<\lambda^+ \rangle$ is an
increasing continuous sequence of models in $K_\lambda$. \item
$\{p_{\alpha,\beta}:\beta<\lambda^+\}=S^{bs}(M_\alpha)$. \item
$M_{\alpha+1}$ realizes $p_{\gamma,\beta}$, where we denote:
$A_\alpha:= \{cd(\gamma,\beta):\gamma \leq \alpha,\
p_{\gamma,\beta}$ is not realized in $M_{\alpha}$\},
$\epsilon_\alpha=\Min(A_\alpha)$ and
$(\gamma,\beta)=cd^{-1}(\epsilon_\alpha)$.
\end{enumerate}
We argue that $M_{\lambda^+}:=\bigcup
\{M_\alpha:\alpha<\lambda^+\}$ is saturated in $\lambda^+$ over
$\lambda$. By 2 it is sufficient to prove that for every $\alpha
\in \lambda^+$ and every $p \in S^{bs}(M_\alpha)$ there is $\beta
\in (\alpha,\lambda^+)$ such that $p$ is realized in $M_\beta$.
Towards a contradiction choose $\alpha^*$ so that $p \in
S^{bs}(M_{\alpha^*})$ is not realized in $M_{\lambda^+}$. There is
$\beta<\lambda^+$ such that $p=p_{\alpha^*,\beta}$. Denote
$\epsilon:=cd(\alpha^*,\beta)$. For every $\alpha \geq \alpha^*$
$\epsilon \in A_\alpha$, so $A_\alpha$ is nonempty and
$\epsilon_\alpha$ is defined. But $\epsilon_\alpha \neq \epsilon$,
(as otherwise $p$ is realized in $M_{\alpha+1}$), so
$\epsilon_\alpha<\epsilon$. The function $f:[\alpha^*,\lambda^+)
\rightarrow \epsilon$, $f(\alpha)=\epsilon_\alpha$ is injection
which is impossible.

It remains to prove item 1. Fix $N$, with $M_0 \prec N$. It is
sufficient to prove that there is an embedding of $N$ to
$M_{\lambda^+}$ over $M_0$. We choose
$(\alpha_\epsilon,N_\epsilon,f_\epsilon)$ by induction on
$\epsilon<\lambda^+$ such that:

\begin{displaymath}
\xymatrix{N \ar[r]^{id} & N_\epsilon \ar[r]^{id} &
N_{\epsilon+1}\\
M_0 \ar[u]^{f_0} \ar[r]^{id} & M_{\alpha_\epsilon}
\ar[u]^{f_\epsilon} \ar[r]^{id} & M_{\alpha_{\epsilon+1}}
\ar[u]^{f_{\epsilon+1}}}
\end{displaymath}

\begin{enumerate}
\item $\langle \alpha_\epsilon:\epsilon<\lambda^+ \rangle$ is an
increasing continuous sequence of ordinals in $\lambda^+$. \item
The sequence $\langle N_\epsilon:\epsilon<\lambda^+ \rangle$ is
increasing and continuous. \item $\langle
f_\epsilon:\epsilon<\lambda^+ \rangle$ is increasing continuous.
\item $N_0:=N$, $\alpha_0:=0$ and $f_0=id_{M_0}$. \item
$f_\epsilon:M_{\alpha_\epsilon} \to N_\epsilon$ is an embedding.
\item For every $\alpha \in S$ there is $a \in
M_{\alpha_{\epsilon+1}}-M_{\alpha_\epsilon}$ such that
$f_{\epsilon+1}(a) \in N_\epsilon$.
\end{enumerate} By Proposition \ref{1.11}
we cannot carry out this construction. Where will we get stuck?
For $\epsilon=0$ or limit we will not get stuck. Suppose we have
defined $(\alpha_\zeta,N_\zeta,f_\zeta)$ for $\zeta \leq
\epsilon$. If $f_\epsilon[M_{\alpha_\epsilon}]=N_\epsilon$ then
$f_{\epsilon}^{-1}\restriction N$ is an embedding of $N$ into
$M_{\lambda^+}$ over $M_0$, hence we are finished. So without loss
of generality $f_\epsilon[M_{\alpha_\epsilon}] \neq N_\epsilon$.
If $\alpha_\epsilon$ $\notin S$ then we define
$\alpha_{\epsilon+1}:=\alpha_\epsilon+1$ and use the amalgamation
in $(K_\lambda, \preceq \restriction K_\lambda)$ to find
$N_{\epsilon+1},f_{\epsilon+1}$ as needed.

Suppose $\alpha_\epsilon \in S$. By the theorem's assumption,
there is a type $p \in S(M_{\alpha_\epsilon})$ such that $p$ is
realized in $M_{\lambda^+}$ and $f_\epsilon(p)$ is realized in
$N_\epsilon$. Define $\alpha_{\epsilon+1}:=\Min\{\alpha \in
\lambda^+:p$ is realized in $M_\alpha$ \}. Take $a \in
M_{\alpha_{\epsilon+1}}$ such that
$tp(a,M_{\alpha_\epsilon},M_{\alpha_{\epsilon+1}})=p$ and take $b
\in N_\epsilon$ such that
$tp(b,f_{\epsilon}(M_{\alpha_\epsilon}),N_\epsilon)=f_\epsilon(p)$.
Then
$f_\epsilon(tp(a,M_{\alpha_\epsilon},M_{\lambda^+}))=tp(b,M_{\alpha_\epsilon},N_\epsilon)$.
By the definition of type (Definition \ref{definition of a type}.1
on page \pageref{definition of a type}), there are
$N_{\epsilon+1},f_{\epsilon+1}$ with $N_{\epsilon} \preceq
N_{\epsilon+1},\ f_{\epsilon+1}$ is an embedding of
$M_{\alpha_{\epsilon+1}}$ into $N_{\epsilon+1}$, $f_\epsilon
\subseteq f_{\epsilon+1}$ and $f_{\epsilon+1}(a)=b$.

Since the hypotheses of 5 applies to any cofinal segment of the
sequence $\langle M_\alpha:\alpha<\lambda^+ \rangle$ and any
submodel of size $\lambda$ lies in some $M_\alpha$ we conclude
that $M_{\lambda^+}$ is saturated in $\lambda^+$ over $\lambda$.
\end{proof}

\subsection{non-forking with greater models}
Now we extend our non-forking notion to include models of
cardinality greater than $\lambda$.

\begin{definition}\label{preparation for forking for big models}
\label{2.9} $\dnf^{\geq \lambda}$ is the class of quadruples
$(M_0,a,M_1,M_2)$ such that:
\begin{enumerate}
\item $M_0 \preceq M_1 \preceq M_2$. \item $\lambda \leq ||M_0||$.
\item For some model $N_0 \in K_\lambda$ with $N_0 \preceq M_0$
for each model $N \in K_\lambda$, $N_0 \preceq N \preceq M_1
\Rightarrow \dnf(N_0,a,N,M_2)$.
\end{enumerate}
\end{definition}

\begin{remark}
If $\dnf(M_0,a,M_1,M_2)$ then $\dnf^{\geq
\lambda}(M_0,a,M_1,M_2)$.
\end{remark}

\begin{definition}\label{forking for big models}
Let $M_0 \preceq M_1$ and $p \in S(M_1)$. We say that $p$ does not
fork over $M_0$, when for some triple $(M_1,M_2,a) \in p$ we have
$\dnf^{\leq \lambda}(M_0,a,M_1,M_2)$.
\end{definition}

\begin{remark}
We can replace the `for some' in Definition \ref{forking for big
models} by `for each'.
\end{remark}

\begin{definition}\label{basic for big models}
Let $M \in K_{>\lambda},\ p \in S(M)$. $p$ is said to be
\emph{basic} when there is $N \in K_\lambda$ such that $N \preceq
M$ and $p$ does not fork over $N$. For every $M \in K_{>\lambda},\
S^{bs}_{> \lambda}(M)$ is the set of basic types over $M$.
Sometimes we write $S^{bs}_{\geq \lambda}(M)$, meaning $S^{bs}(M)$
or $S^{bs}_{> \lambda}(M)$ (the unique difference is the
cardinality of $M$).
\end{definition}

Now we present a weak version of local character, which is
peripheral in the continuation of this paper.
\begin{definition}\label{weak local character for <^* sequences}
Let $\frak{s}$ be a semi good frame except local character.
$\frak{s}$ is said to satisfy \emph{weak local character for
$\prec^*$-increasing sequences} when: If $ \alpha^*< \lambda^+$
and $\langle M_\alpha:\alpha \leq \alpha^*+1 \rangle$ is an
$\prec^*-$increasing continuous sequence of models, then for some
element $a \in M_{\alpha^*+1}-M_{\alpha^*}$ and ordinal
$\alpha<\alpha^*$, $tp(a,M_{\alpha^*},M_{\alpha^*+1})$ does not
fork over $M_\alpha$.
\end{definition}

\begin{definition}\label{weak local character for fast sequences}
Let $\frak{s}$ be a semi good $\lambda$-frame except local
character. $\frak{s}$ is said to satisfy \emph{weak local
character for fast sequences} when for some relation $\prec^*$ the
following hold:
\begin{enumerate}
\item $\prec^*$ is a relation on $K_\lambda$. \item If $M_0
\prec^* M_1$ then $M_0 \preceq M_1$. \item If $M_0 \prec^*M_1
\preceq M_2 \in K_\lambda$ then $M_0 \prec^* M_2$. \item
$\frak{s}$ satisfies weak local character for $\prec^*$-increasing
sequences. \item If $M_0 \in K_\lambda$, $M_0 \prec M_2 \in
K_{\lambda^+}$, then there is a model $M_1 \in K_\lambda$ such
that: $M_0 \prec^* M_1 \preceq M_2$.
\end{enumerate}
\end{definition}

\begin{remark}
If $\frak{s}$ is a semi good $\lambda$-frame and $\prec^*$ is a
relation on $K_\lambda$ such that $M \prec^* N \Rightarrow M
\preceq N$ then $\frak{s}$ satisfies weak local character for
$\prec^*$-increasing sequences.
\end{remark}


The following theorem asserts that a non-forking relation in
$(K_\lambda, \preceq \restriction K_\lambda)$ can be lifted to
$K_{\geq \lambda}$ with many properties preserved.
\begin{theorem} \label{2.10}\label{lifting}
Let $\frak{s}$ be a semi-good $\lambda$-frame, except local
character.
\begin{enumerate}
\item Density: If $\frak{s}$ satisfies weak local character for
fast sequences and $M\prec N,\ M \in K_{\geq \lambda}$
 then there is $a
\in N-M$ such that $tp(a,M,N) \in S^{bs}_{\geq \lambda}(M)$. \item
Monotonicity: Suppose $M_0 \preceq M_1 \preceq M_2,\ n<3
\Rightarrow M_n \in K_{\geq \lambda},\
||M_2||\allowbreak>\lambda$. If $p \in S^{bs}_{\geq \lambda}(M_2)$
does not fork over $M_0$, then
\begin{enumerate}
\item $p$ does not fork over $M_1$. \item $p\restriction M_1$ does
not fork over $M_0$.
\end{enumerate}
\item Transitivity: Suppose $M_0,M_1,M_2 \in K_{\geq \lambda}$ and
$M_0 \preceq M_1 \preceq M_2$. If $p \in S^{bs}_{\geq
\lambda}(M_2)$ does not fork over $M_1$, and $p\restriction M_1$
does not fork over $M_0$, then $p$ does not fork over $M_0$. \item
About local character: Let $\delta$ be a limit ordinal. Suppose
$\frak{s}$ satisfies local character or $\lambda^+ \leq
cf(\delta)$. If $\langle M_\alpha:\alpha \leq \delta \rangle$ is
an increasing continuous sequence of models in $K_{>\lambda}$, and
$p \in S^{bs}_{>\lambda}(M_\delta)$ then there is $\alpha<\delta$
such that $p$ does not fork over $M_\alpha$. \item Continuity:
Suppose $\langle M_\alpha:\alpha \leq \delta+1 \rangle$ is an
increasing continuous sequence of models in $K_{\geq\lambda}$. Let
$c \in M_{\delta+1}-M_\delta$. Denote
$p_\alpha=tp(c,M_\alpha,M_{\delta+1})$. If for every
$\alpha<\delta,\ p_\alpha$ does not fork over $M_0$, then
$p_\delta$ does not fork over $M_0$.
\end{enumerate}
\end{theorem}
\begin{proof}
\mbox{}
(1) Density: Suppose $M\prec N$.\\
\case{Case 1:} $||M||=\lambda$. Choose $a \in N-M$.
$LST(K,\preceq) \leq \lambda$ and so there is $N^*\prec N$ such
that: $||N^*||=\lambda$ and $M \bigcup \{a\} \subseteq N^*$. By
axiom e of a.e.c $M \preceq N^*$ But $a \in N^*-M$ and so $M\prec
N^*$. By the existence axiom in $\frak{s}$, there is $c \in N^*-M$
such that $tp(c,M,N^*)$ is basic. So $tp(c,M,N) \in
S^{bs}(M)$.\\
\case{Case 2:} $||M||>\lambda$. We choose $M_n,N_n$ by induction
on $n<\omega$ such that:

\begin{displaymath}
\xymatrix{c \in N_n \ar[rr]^{id} && N_n \ar[r]^{id} & N_\omega \ar[r]^{id} & N\\
M_n \ar[u]^{id} \ar[r]^{id} & M_{n,c} \ar[r]^{id} & M_{n+1}
\ar[u]^{id} \ar[r]^{id} & N_\omega \ar[r]^{id} \ar[u]^{id} & M
\ar[u]^{id} }
\end{displaymath}

\begin{enumerate}[(a)] \item $\langle N_n:n \leq \omega \rangle$ is an
$\prec-$increasing continuous sequence of models in $K_\lambda$.
\item $\langle M_n:n \leq \omega \rangle$ is an
$\prec^*-$increasing continuous sequence of models in $K_\lambda$.
\item $M_n\prec M$ (see the end of Definition \ref{2.1a}). \item
$N_n\prec N$. \item $N_0 \nsubseteq M$. \item For every $c \in
N_n,\ M_{n,c} \subseteq M_{n+1}$ where we choose $M_{n,c} \in
K_\lambda$ such that: If $tp(c,M_n,N_n) \in S^{bs}(M_n)$ but does
fork over $M_n$ then $M_{n,c}$ is a witness for this, namely
$M_n\prec M_{n,c} \prec M$ and $tp(c,M_{n,c},N)$ forks over $M_n$.
Otherwise $M_{n,c}=M_n$.
\end{enumerate}

The construction is of course possible.

Now we define $M_\omega:=\bigcup \{M_n:n<\omega\}$ and
$N_\omega:=\bigcup \{N_n:n<\omega\}$. By axiom \ref{1.1}.1.d
(smoothness) $M_\omega \preceq N_\omega$. By the local character
for $\prec^*$-increasing sequences, for some element $c \in
N_\omega-M_\omega$ and there is $n<\omega$ such that
$tp(c,M_\omega,N_\omega) \in S^{bs}(M_\omega)$ does not fork over
$M_n$. By the monotonicity without loss of generality $c \in N_n$.
We will prove that $tp(c,M,N)$ does not fork over $M_\omega$. Take
$M^*$ with $M_\omega \prec M^{*}\prec M$. By way of contradiction
suppose $tp(c,M^{*},N)$ forks over $M_\omega$. By the monotonicity
in $\frak{s}$ (axiom b), $tp(c,M^{*},N)$ forks over $M_n$. So by
the definition of $M_{n,c},\ tp(c,M_{n,c},N)$ forks over $M_n$.
Hence by axiom b (monotonicity) $tp(c,M_\omega,N)$ forks over
$M_n$, a contradiction.

(2) Monotonicity: We use the same witness.

(3) Transitivity:

\begin{displaymath}
\xymatrix{& N \ar[r]^{id} & N^{**} \ar[r]^{id} & M_2 & p\\
N_1 \ar[rr]^{id} && N^* \ar[r]^{id} \ar[u]^{id} & M_1
\ar[u]^{id}\\
& N_0 \ar[rr]^{id} \ar[uu]^{id} \ar[ru]^{id} && M_0 \ar[u]^{id} }
\end{displaymath}

Suppose $M_0\prec M_1\prec M_2$, $p \in S^{bs}(M_2)$ does not fork
over $M_1$ and $p\restriction M_1$ does not fork over $M_0$. We
can find $N_0\prec M_0$ such that $N_0$ witnesses that
$p\restriction M_1$ does not fork over $M_0$. We will prove that
$N_0$ witnesses that $p$ does not fork over $M_0$. Let $N \in
K_\lambda$ be such that $N_0\prec N\prec M_2$. We have to prove
that $p\restriction N $ does not fork over $N_0$. First we find a
model $N_1$ that witnesses that $p$ does not fork over $M_1$. As
$LST(K,\preceq) \leq \lambda$ there is $N^* \in K_\lambda$ such
that $N_0 \bigcup N_1 \subseteq N^* \preceq M_1$ and there is
$N^{**} \in K_\lambda$ such that $N^* \bigcup N \subseteq N^{**}
\preceq M_2$. As $N_1$ witnesses that $p$ does not fork over
$M_1$, $p\restriction N ^{**}$ does not fork over $N_1$. By the
Definition \ref{2.1a}.3.b (monotonicity), $p\restriction N^{**}$
does not fork over $N^*$. $N_0$ witnesses that $p\restriction M_1$
does not fork over $M_0$, so $p\restriction N ^*$ does not fork
over $N_0$. By the transitivity proposition (Proposition \ref{s is
transitive}), $p\restriction N^{**}$ does not fork over $N_0$. So
by Definition \ref{2.1a}.3.b (monotonicity), $p\restriction N $
does not fork over $N_0$.

(4) About local character: Let $\langle M_\alpha:\alpha<\delta
\rangle$ be an increasing continuous sequence of models in
$K_{>\lambda}$. Let $p \in S^{bs}_{>\lambda}(M_\delta)$ and $N^*$
a witness for this, i.e. $p$ does not fork over $N^* \in
K_\lambda$. Let $\langle \alpha(\epsilon):\epsilon \leq cf(\delta)
\rangle$ be an increasing continuous sequence of ordinals with
$\alpha(cf(\delta))=\delta$.

\case{Case a:} $\lambda^+ \leq cf(\delta)$.  By cardinality
considerations, there is $\epsilon<cf(\delta)$ such that: $N^*
\subseteq M_{\alpha(\epsilon)}$. By axiom \ref{1.1}.1.e $N^*
\preceq M_{\alpha(\epsilon)}$. As $N^*$ witnesses that the type
$p$ is basic, by Definition \ref{2.9} $N^*$ witnesses that $p$
does not fork over $M_{\alpha(\epsilon)}$.

\case{Case b:} $\frak{s}$ satisfies local character and
$cf(\delta) \leq \lambda$. Using $LST(K,\preceq) \leq \lambda$ and
smoothness, we can choose $N_{\alpha(\epsilon)}$ by induction on
$\epsilon \leq cf(\delta)$ such that:

\begin{displaymath}
\xymatrix{N^* \ar[r]^{id} & N_\delta \ar[r]^{id} & M_\delta & p\\
 & N_{\alpha(\epsilon)} \ar[r]^{id} \ar[u]^{id} &
 M_{\alpha(\epsilon)}
\ar[u]^{id} }
\end{displaymath}

\begin{enumerate}[(a)]  \item $N_{\alpha(\epsilon)} \in
K_\lambda$. \item $\langle N_{\alpha(\epsilon)}:\epsilon \leq
cf(\delta) \rangle$ is an increasing continuous sequence. \item
$M_{\alpha(\epsilon)} \bigcap N^* \subseteq N_{\alpha(\epsilon)}
\preceq M_{\alpha(\epsilon)}$.
\end{enumerate}

By axiom \ref{1.1}.1.e $N^* \preceq N_\delta \preceq M_\delta$.
Since $p$ does not fork over $N^*$, by monotonicity (Theorem
\ref{lifting}.2) $p$ does not fork over $N_\delta$. By local
character, for some $\epsilon<cf(\delta)$, $p\restriction
N_\delta$ does not fork over $N_{\alpha(\epsilon)}$. By
transitivity (Theorem \ref{lifting}.3), $p$ does not fork over
$N_{\alpha(\epsilon)}$. By monotonicity (Theorem \ref{lifting}.2),
$p$ does not fork over $M_{\alpha(\epsilon)}.$

(5) Continuity: For every $\alpha \in \delta$ denote
$p_\alpha:=p\restriction M_\alpha$. $p_7$ does not fork over
$M_0$. So for some $N_0 \in K_\lambda$, $N_0 \preceq M_0$ and
$p_7$ does not fork over $N_0$. By monotonicity (Theorem
\ref{lifting}.2) and transitivity (Theorem \ref{lifting}.2) for
every $\alpha<\delta$ $p_\alpha$ does not fork over $N_0$. We will
prove that $p$ does not fork over $N_0$. Take $N^* \in K_\lambda$
with $N_0 \preceq N^* \preceq M_\delta$. We have to prove that
$p\restriction N^*$ does not fork over $N_0$. Let $\langle
\alpha(\epsilon):\epsilon \leq cf(\delta) \rangle$ be an
increasing continuous sequence of ordinals with
$\alpha(cf(\delta))=\delta$.

\case{Case a:} $\lambda^+ \leq cf(\delta)$. By cardinality
considerations there is $\epsilon<cf(\delta)$ such that $N^*
\subseteq M_{\alpha(\epsilon)}$. But $M_{\alpha(\epsilon)} \preceq
M_\delta$ and $N^* \preceq M_\delta$, so by axiom \ref{1.1}.1.e
$N^* \preceq M_{\alpha(\epsilon)}$. Since $p_{\alpha(\epsilon)}$
does not fork over $N_0$, by monotonicity (Theorem
\ref{lifting}.2) $p \restriction N^*$ does not fork over $N_0$.

\case{Case b:} $cf(\delta) \leq \lambda^+$. We choose
$N_{\alpha(\epsilon)}$ by induction on $\epsilon \in
(0,cf(\delta)]$ such that:
\begin{enumerate}[(a)]
\item The sequence $\langle N_{\alpha(\epsilon)}:\epsilon \leq
cf(\delta) \rangle$ is increasing continuous. \item $\epsilon \leq
cf(\delta) \Rightarrow N^* \bigcap M_{\alpha(\epsilon)} \subseteq
N_{\alpha(\epsilon)} \preceq M_{\alpha(\epsilon)}$. \item
$N_{\alpha(\epsilon)} \in K_\lambda$.
\end{enumerate}

For every $\epsilon<\cf(\delta)$, $p_{\alpha(\epsilon)}$ does not
fork over $N_0$, so $p \restriction N_{\alpha(\epsilon)}$ does not
fork over $N_0$. So by Definition \ref{2.1a}.3.g (continuity) (in
$\frak{s}$), $p\restriction N _\delta$ does not fork over $N_0$.
$N^* \subseteq N_\delta$, hence by axiom \ref{1.1}.1.e $N^*
\preceq N_\delta$. Therefore by Definition \ref{2.1a}.3.b
(monotonicity), $p\restriction N^*$ does not fork over $N_0$.

\end{proof}

\section{The decomposition and amalgamation method}
\discussion{Discussion} In section 2 we defined an extension of
the non-forking notion to cardinals bigger than $\lambda$. But we
did not prove all the good frame axioms. The purpose from here
until the end of the paper is to construct a good
$\lambda^+$-frame, which is derived from the semi good
$\lambda$-frame. In a sense, the main problem is that amalgamation
in $(K_\lambda, \preceq \restriction K_\lambda)$ does not imply
amalgamation in $(K_{\lambda^+},\preceq \restriction
K_{\lambda^+})$. Suppose for $n<3$ $M_n \in K_{\lambda^+}$, $M_0
\preceq M_n$ and we want to amalgamate $M_1,M_2$ over $M_0$. We
take representation of $M_0,M_1,M_2$. We want to amalgamate
$M_1,M_2$ by amalgamating their representations. For this goal, we
will find in section 5, a relation of ``a canonical amalgamation''
or ``a non-forking amalgamation''. Sections 3,4 are preparations
for section 5. If the reader wants to know the plan of the other
sections now, he may see the discussion at the beginning of
section 10.
 \discussion{The decomposition and amalgamation method} Suppose
 for $n=1,2$ $M_0\preceq M_n$
and we want to prove that there is an amalgamation of $M_1,M_2$
over $M_0$ which satisfies specific properties (for example
disjointness or uniqueness, see below). Sometimes there is a
property of triples, $K^{3,*} \subseteq K^3$ such that if
$(M_0,M_1,a) \in K^{3,*}$ and $(M_0,M_1,a) \preceq (M_2,M_3,a)$
then the amalgamation $M_3$ satisfies the required property. A
classic example of this property in the context of fields is
`$M_1$ is the algebraic closure of $M_0 \bigcup \{a\}$. What can
we do, if there is no $a \in M_1-M_0$ such that $(M_0,M_1,a) \in
K^{3,*}$? Theorem \ref{3.7} says under some assumptions that we
can decompose an extension of $M_1$ over $M_0$ by triples in
$K^{3,*}$. By Proposition \ref{3.3}.2 we may amalgamate $M_2$ with
the decomposition we have obtained. \discussion{Applications of
the decomposition and amalgamation method}

\begin{enumerate}
\item By Proposition \ref{3.3}(2) there is no $\preceq$-maximal
model in $K_{\lambda^+}$. \item By Proposition \ref{3.12} the
uniqueness triples are dense with respect to $\preceq_{bs}$ (see
Definition \ref{3.1}.2). It enables to prove Theorem \ref{3.13}
(the disjoint amalgamation existence), by the decomposition and
disjoint method. \item By
assumption \ref{5.0} the uniqueness triples are dense with respect
to $\preceq_{bs}$. The density enables to prove Theorem \ref{the
existence theorem for NF} (the exitance theorem for $NF$). \item
Using again assumption \ref{5.0}, we prove Proposition
\ref{opposite} . But for this, we have to prove Proposition
\ref{3.4}, a generalization of \ref{3.3}, which says that we can
amalgamate two sequences of models, not just a model and a
sequence.
\end{enumerate}

\begin{hypothesis} \label{3.0}
$\frak{s}$ is a semi good $\lambda$-frame, except basic weak
stability and local character.
\end{hypothesis}

\subsection{The a.e.c. $(K^{3,bs},\preceq_{bs})$ and amalgamations}

\begin{definition} \label{3.1}  \mbox{} \begin{enumerate}
\item $K^{3,bs}=:\{(M,N,a):M,N \in K_\lambda,\ a \in N-M$ and
$tp(a,M,N) \in S^{bs}(M)\}$. \item $\preceq_{bs}$ is the relation
on $K^{3,bs}$ defined by: $(M,N,a) \preceq_{bs} (M*,N^*,a^*)$ iff
$M\preceq M^*,\ N\preceq N^*,\ a^*=a$
 and $tp(a,M^*,N^*)$ does not fork over $M$.
\item The sequence $\langle (M_\alpha,N_\alpha,a):\alpha<\theta
\rangle$ is said to be $\preceq_{bs}$-increasing continuous if
$\alpha <\theta\Rightarrow (M_\alpha,N_\alpha,a)\preceq_{bs}
(M_{\alpha+1},N_{\alpha+1},a)$ and the sequences $\langle
(M_\alpha:\alpha<\theta \rangle,\ \langle N_\alpha:\alpha<\theta
\rangle$ are continuous (and clearly also increasing).
\end{enumerate}
\end{definition}

\begin{proposition} \label{3.2}
$(K^{3,bs},\preceq_{bs})$ is an a.e.c. in $\lambda$ and it has no
$\preceq_{bs}$-maximal model.

\end{proposition}

\begin{proof}
First we note that $K^{3,bs}$ is not the empty set [There is $M
\in K_\lambda$, and as it has no $\preceq$-maximal model, there is
$N \in K_\lambda$ with $M\prec N$. Now by Definition
\ref{2.1a}.3.f, there is $a \in N-M$ such that $tp(M,N,a) \in
S^{bs}(M)$]. Why is axiom c of a.e.c. (defintion \ref{1.1}.1.c)
satisfied? Suppose $\delta<\lambda^+$ and $\langle
(M_\alpha,N_\alpha,a):\alpha<\delta \rangle$ is increasing and
continuous. Denote $M=\bigcup \{M_\alpha:\alpha<\delta\},\
N=\bigcup \{N_\alpha:\alpha<\delta\}$. By axiom c of a.e.c., $M,N
\in K_\lambda,\ \alpha<\delta \Rightarrow M_\alpha \preceq M,\
N_\alpha \preceq N$. By the definition of $\preceq_{bs}$ for every
$\alpha<\delta$, $tp(a,M_\alpha,N_\alpha)$ does not fork over
$M_0$. So by Definition \ref{2.1a}.3.g (continuity), $tp(a,M,N)$
is basic and does not fork over $M_0$. By the smoothness, $M
\preceq N$. By axiom c of a.e.c. $M_0 \preceq M$ and $N_0 \preceq
N$. So $(M_0,N_0,a) \preceq _{bs}(M,N,a) \in K^{3,bs}$. Why is the
smoothness satisfied? Suppose $\langle
(M_\alpha,N_\alpha,a):\alpha\leq \delta+1 \rangle$ is continuous
and for $\alpha<\beta \leq \delta+1$, we have $\alpha \neq \delta
\Rightarrow (M_\alpha,N_\alpha,a)
\preceq_{bs}(M_\beta,N_\beta,a)$. So $\delta \neq \alpha < \beta
\leq \delta+1 \Rightarrow M_\alpha \preceq M_\beta$. But by the
continuity of the sequence $\langle (M_\alpha,N_\alpha,a):\alpha
\leq \delta+1 \rangle$ we have $M_\delta=\bigcup
\{M_\alpha:\alpha<\delta\}$. So by the smoothness of
$(K,\preceq)$, $M_\delta \preceq M_{\delta+1}$. In a similar way
$N_\delta \preceq N_{\delta+1}$. $(M_0,N_0,a) \preceq_{bs}
(M_{\delta+1},N_{\delta+1},a)$, so by the definition,
$tp(a,M_{\delta+1},N_{\delta+1})$ does not fork over $M_0$.
Therefore by Definition \ref{2.1a}.3.b (monotonicity),
$tp(a,M_{\delta+1},N_{\delta+1})$ does not fork over $M_\delta$.
Why does $(K^{3,bs},\preceq_{bs})$ satisfy axiom \ref{1.1}.1.e?
Suppose $(M_0,N_0,a) \subseteq (M_1,N_1,a) \preceq (M_2,N_2,a),\
(M_0,N_0,a) \preceq_{bs} (M_2,N_2,a)$. By the definition of
$\preceq_{bs}$ we have $M_0 \subseteq M_1 \preceq M_2$ and $M_0
\preceq M_2$. Hence by axiom \ref{1.1}.1.e we have $M_0 \preceq
M_1$. In a similar way $N_0 \preceq N_1$. By the definition of
$\preceq_{bs}$, $tp(a,M_2,N_2)$ does not fork over $M_0$. By
\ref{definition good frame}.3.b (monotonicity), $tp(a,M_1,N_1)$
does not fork over $M_0$. So $(M_0,N_0,a) \preceq_{bs}
(M_1,N_1,a)$. Why is there no maximal element in
$(K^{3,bs},\preceq_{bs})$? Let $(M_0,N_0,a) \in K^{3,bs}$. In
$K_\lambda$ there is no $\preceq$-maximal element, and so there is
$M_0\prec M_1^* \in K_\lambda$. By Proposition \ref{axiom i is
redundant} there is $N_0 \preceq N_1 \in K_\lambda$ and there is
an embedding $f:M_1^* \to N_1$ such that $tp(a,M_1,N_1)$ does not
fork over $M_0$ where $M_1:=f[M_1^*]$. Hence $(M_0,N_0,a)
\preceq_{bs} (M_1,N_1,a)$.
\end{proof}

\begin{proposition} \label{3.3}

\begin{displaymath}
\xymatrix{
N_0 \ar[r]^{id} & N_1 \ar[r]^{id} & N_2 \ar[rr]^{id} && N_\alpha \ar[r]^{id} & N_{\alpha+1} \ar[rr]^{id} && N_\theta\\
M_0 \ar[r]^{id} \ar[u]^{id} & M_1 \ar[r]^{id} \ar[u]^{id} & M_2
\ar[rr]^{id} \ar[u]^{id} && M_\alpha \ar[r]^{id} \ar[u]^{id} &
M_{\alpha+1} \ar[rr]^{id} \ar[u]^{id} && M_\theta \ar[u]^{id}}
\end{displaymath}

\begin{enumerate}
\item Let $\langle M_\alpha:\alpha \leq \theta \rangle$ be an
increasing continuous sequence of models in $K_\lambda$. Let $N
\in K_\lambda$ with $M_0\prec N$, and for $\alpha<\theta$, let
$a_\alpha \in M_{\alpha+1}-M_\alpha,\
(M_\alpha,M_{\alpha+1},a_\alpha) \in K^{3,bs}$ and $b \in N-M_0,\
(M_0,N,b) \in K^{3,bs}$. Then there are $f,\langle N_\alpha:\alpha
\leq \theta \rangle$ such that:
\begin{enumerate}
\item $f$ is an isomorphism of $N$ to $N_0$ over $M_0$. \item
$\langle N_\alpha:\alpha \leq \theta \rangle$ is an increasing
continuous sequence. \item $M_\alpha \preceq N_\alpha$. \item
$tp(a_\alpha,N_\alpha,N_{\alpha+1})$ does not fork over
$M_\alpha$. \item $tp(f(b),M_\alpha,N_\alpha)$ does not fork over
$M_0$.
\end{enumerate}
\item $K_{\lambda^+}\neq \emptyset$, and it has no
$\preceq$-maximal model. \item There is a model in $K$ of
cardinality $\lambda^{+2}$.
\end{enumerate}
\end{proposition}

\begin{proof}
\mbox{} (1) First we explain the idea of the proof. If we `fix'
the models in the sequence $\langle M_\alpha:\alpha \leq \theta
\rangle$, then we will `change' $N$ $\theta$ times. So in limit
steps we will be in a problem. The solution is to fix $N$, and
`change' the sequence $\langle M_\alpha:\alpha \leq \theta
\rangle$. At the end of the proof we `return the sequence to its
place'.

\emph{The proof itself:} We choose $(N_\alpha^*,f_\alpha)$ by
induction on $\alpha$ such that $(*)_\alpha$ holds where
$(*)_\alpha$ is:
\begin{enumerate}[(i)]
\item $\alpha \leq \theta \Rightarrow N_\alpha^* \in K_\lambda$.
\item $(N_0^*,f_0)=(N,id_{M_0})$. \item The sequence $\langle
N_\alpha^*:\alpha \leq \theta \rangle$ is increasing and
continuous. \item For every $\alpha \leq \theta$, the function
$f_\alpha$ is an embedding of $M_\alpha$ to $N_\alpha^*$. \item
The sequence $\langle f_\alpha:\alpha \leq \theta \rangle$ is
increasing and continuous. \item For every $\alpha<\theta$
$tp(f_\alpha(a_\alpha),N_\alpha^*,N_{\alpha+1}^*)$ does not fork
over $f_\alpha[M_\alpha]$. \item For every $\alpha \leq \theta$
$tp(b,f_\alpha[M_\alpha],N_\alpha^*)$ does not fork over $M_0$.
\end{enumerate}

Now $f_\theta:M_\theta \to N_\theta^*$ is an embedding. Extend
$f_\theta^{-1}$ to a function with domain $N_\theta^*$ and define
$f:=g \restriction N$.
Define $N_\alpha:=g[N_\alpha^*]$.


(2) $K_{\lambda^+} \neq \emptyset$, as we can choose an increasing
continuous sequence of models in $K_\lambda$, $\langle
M_\alpha:\alpha<\lambda^+ \rangle$, and so its union is a model in
$K_{\lambda^+}$ [As there is no $\preceq$-maximal model in
$K_\lambda$ and in limit step use axiom \ref{1.1}.1.c].

Why is there no maximal model in $K_{\lambda^+}$? Let $M \in
K_{\lambda^+}$. Let $\langle M_\alpha:\alpha<\lambda^+ \rangle$ be
a representation of $M$. By the Definition \ref{2.1a}.3.f
(existence, on page \pageref{2.1a}), for every $\alpha \in
\lambda^+$ there is an element $a_\alpha \in
M_{\alpha+1}-M_\alpha$ (we do not use $a_\alpha$, but as we have
written it in 1, for shortness, we have to write it here). As in
$K_\lambda$ there is no maximal model, there is a model $N$ such
that $M_0\prec N \in K_\lambda$ and without loss of generality $N
\bigcap M=M_0$. By Definition \ref{definition good frame}.2.c (the
density of basic types), there is $b \in N-M_0$
such that $tp(b,M_0,N)$ is basic. Now by 
Proposition \ref{3.3}.1, there is an increasing continuous
sequence $\langle N_\alpha:\alpha<\lambda^+ \rangle$ and $f$ such
that $f:N \to N_0$ is an isomorphism over $M_0$ and for $\alpha
\in \lambda^+$ we have $M_\alpha \preceq N_\alpha$ and
$tp(f(b),M_\alpha,N_\alpha)$ does not fork over $M_0$. So by
Definition \ref{2.1a}, (on page \pageref{2.1a}), $f(b)$ does not
belong to $M_\alpha$ for $\alpha \in \lambda^+$. So $f(b)$ does
not belong to $M$. But it belongs to $N_{\lambda^+}$, so $M \neq
N_{\lambda^+}$, and for this we defined $b$. But it is easy to see
that $M \subseteq N_{\lambda^+}$ and $N_{\lambda^+} \in
K_{\lambda^+}$. By the smoothness (i.e. Definition \ref{1.1}.1.d
on page \pageref{1.1}) $M \preceq N_{\lambda^+}$. So $M$ is not a
maximal model.

(3) We construct a strictly increasing continuous sequence of
models in $K_{\lambda^+}$, $\langle M_\alpha:\alpha<\lambda^{+2}
\rangle$. So its union is a model in $K_{\lambda^{+2}}$. As by 2
there is no maximal model in $K_{\lambda^+}$, there is no problem
to choose this sequence.
\end{proof}

\begin{proposition}[a rectangle which amalgamate two
sequences] \label{3.4}
For $x=a,b$ let $\langle
M_{x,\alpha}:\alpha<\theta^x \rangle$ be an increasing continuous
sequence of models in $K_\lambda$ such that $M_{a,0}=M_{b,0}$ and
let $\langle d_{x,\alpha}:\alpha<\theta^x \rangle$ be a sequence
such that $d_{x,\alpha} \in M_{x,\alpha+1}-M_{x,\alpha}$, and the
type $tp(d_{x,\alpha},M_{x,\alpha},M_{x,\alpha+1})$ is basic.
Denote $\alpha^*=\theta^a,\ \beta^*=\theta^b$. Then there are a
``rectangle of models'' $\{M_{\alpha,\beta}:\alpha<\alpha^*,\
\beta<\beta^*\}$ and a sequence $\langle f_\beta:\beta<\beta^*
\rangle$ such that:
\begin{enumerate}
\item $(\alpha<\alpha^* \wedge \beta<\beta^*) \Rightarrow
M_{\alpha,\beta} \in K_\lambda$. \item $f_\beta:M_{b,\beta}\to
M_{0,\beta}$ is an isomorphism. \item $M_{\alpha,0}=M_{a,\alpha}$.
\item $f_0$ is the identity on $M_{a,0}=M_{b,0}$. \item $\langle
f_\beta:\beta<\beta^* \rangle$ is increasing and continuous. \item
For every $\alpha,\beta$ which satisfies $\alpha+1<\alpha^*$ and
$\beta<\beta^*$, the type
$tp(d_{a,\alpha},M_{\alpha,\beta},M_{\alpha+1,\beta})$ does not
fork over $M_{\alpha,0}$. \item For every $\alpha,\beta$ which
satisfies $\alpha<\alpha^*$ and $\beta+1<\beta^*$, the type
$tp(d_{b,\beta},M_{\alpha,\beta},M_{\alpha,\beta+1})$ does not
fork over $M_{0,\beta}$. \item If $\bigcup
\{Im(f_\beta):\beta<\beta^* \} \bigcap \bigcup
\{M_{a,\alpha}:\alpha<\alpha^* \} = \bigcup
\{M_{b,\beta}:\beta<\beta^* \} \bigcap \bigcup
\{M_{a,\alpha}:\alpha<\alpha^* \} = M_{a,0}$, then $(\forall \beta
\in \beta^*)f_\beta=id\restriction M_{b,\beta}$. \item For all
$\alpha(1)<\alpha^*$ the sequence $\langle
M_{\alpha(1),\beta}:\beta<\beta^* \rangle$ is increasing and
continuous. \item For all $\beta(1)<\beta^*$ the sequence $\langle
M_{\alpha,\beta(1)}:\alpha<\alpha^* \rangle$ is increasing and
continuous.
\end{enumerate}
\end{proposition}

\begin{displaymath}
\xymatrix{
d_{a,\alpha} \in M_{\alpha+1,0}=M_{a,\alpha+1} \ar[r]^{id} & M_{\alpha+1,\beta} \ar[r]^{id} & M_{\alpha+1,\beta+1}\\
M_{\alpha,0}=M_{a,\alpha}  \ar[r]^{id} \ar[u]^{id} & M_{\alpha,\beta} \ar[r]^{id} \ar[u]^{id} & M_{\alpha,\beta+1} \ar[u]^{id} \\
M_{0,0}=M_{a,0}=M_{b,0} \ar[r]^{id} \ar[u]^{id} &
M_{0,\beta}=f_\beta[M_{b,\beta}] \ar[r]^{id} \ar[u]^{id}&
M_{0,\beta+1}=f_{\beta+1}[M_{b,\beta+1}] \ar[u]^{id}}
\end{displaymath}

\begin{proof} We define by induction on $\beta<\beta^*$ $f_\beta,\{M_{\alpha,\beta}:\alpha<\alpha^*\}$ such that
 the conditions 1-6 and 8,9 are satisfied. For $\beta=0$ see 3,4. For
$\beta$ a limit ordinal, we define $f_\beta=\bigcup
\{f_\gamma:\gamma<\beta \},\ M_{\alpha,\beta}=\bigcup
\{M_{\alpha,\gamma}:\gamma<\beta \}$. Why does 6 satisfy, i.e. why
for every $\alpha$, does
$tp(d_{a,\alpha},M_{\alpha,\beta},M_{\alpha+1,\beta})$ not fork
over $M_{\alpha,0}$? By the induction hypothesis 6 is satisfied
for every $\gamma<\beta$, i.e.
$tp(d_{a,\alpha},M_{\alpha,\gamma},M_{\alpha+1,\gamma})
=tp(d_{a,\alpha},M_{\alpha,\gamma},M_{\alpha+1,\gamma})$ does not
fork over $M_{0,\gamma}$. By Definition \ref{2.1a}.3.b
(monotonicity) and Definition \ref{2.1a}.3.g (continuity)
$tp(d_{a,\alpha},M_{\alpha,\beta},M_{\alpha+1,\beta})$ does not
fork over $M_{\alpha,0}$. So condition 6 is satisfied. For
$\beta=\gamma+1$ use Proposition \ref{3.3}.1. So we can carry out
the induction. Now without loss of generality condition 7 is
satisfied too.
\end{proof}

\subsection{Decomposition}

When we speak about $tp(a,M,N)$ the order of $N$ is peripheral.
Now we consider classes $K^{3,*}$ of triples $(M,N,a)$ where the
order of $N$ is very important. For example $N$ is the algebraic
closure of $M \bigcup \{a\}$, where $(K,\preceq)$ is the class of
fields with the partial order of being sub-field.

\begin{definition} \label{3.5}
Let $K^{3,*} \subseteq K^{3,bs}$ be closed under isomorphisms.
\begin{enumerate}
\item $K^{3,*}$ is \emph{dense with respect to $\preceq_{bs}$} if
for every $(M,N,a) \in K^{3,bs}$ there is $(M^*,N^*,a^*) \in
K^{3,*}$ such that $(M,N,a) \preceq_{bs}(M^*,N^*,a^*)$. \item
$K^{3,*}$ \emph{has existence} if for every $(M,N,a) \in K^{3,bs}$
there are $N^*,a^*$ such that $tp(a^*,M,N^*)=tp(a,M,N)$ and
$(M,N^*,a^*) \in K^{3,*}$. In other words if $p \in S^{bs}(M)$
then $p \bigcap K^{3,*} \neq \emptyset$.
\end{enumerate}
\end{definition}

\begin{definition} \label{3.6}
Let $K^{3,*} \subseteq K^{3,bs}$ be closed under isomorphisms. We
say that $M^*$ is \emph{decomposable by $K^{3,*}$} over $M$, if
there is a sequence $ \langle
d_\epsilon,N_\epsilon:\epsilon<\alpha^* \rangle^\frown \langle
N_{\alpha^*} \rangle$ such that:
\begin{enumerate}
\item $\epsilon<\alpha^* \Rightarrow N_\epsilon \in K_\lambda$.
\item $ \langle N_\epsilon:\epsilon \preceq \alpha^* \rangle$ is
increasing and continuous. \item $N_0=M$. \item
$N_{1,\alpha^*}=M^*$. \item
$(N_\epsilon,N_{\epsilon+1},d_\epsilon) \in K^{3,*}$.
\end{enumerate}
In such a case we say that the sequence $ \langle
d_\epsilon,N_\epsilon:\epsilon<\alpha^* \rangle ^\frown \langle
N_{\alpha^*} \rangle$ is a decomposition of $M^*$ over $M$ by
$K^{3,*}$. The main case is $K^{3,*}=K^{3,uq}$ (which we have not
defined yet), and in such a case we may omit it.
\end{definition}

\begin{theorem}[the extensions decomposition theorem]
\label{3.7} \label{the extensions decomposition theorem} Let
$K^{3,*} \subseteq K^{3,bs}$ be closed under isomorphisms.
\begin{enumerate}
\item Suppose $\frak{s}$ has conjugation. If $K^{3,*}$ is dense
with respect to $\preceq_{bs}$ then it has existence. \item
Suppose $K^{3,*}$ has existence. If $N \in K_\lambda$ and
$p=tp(a,M,N) \in S^{bs}(M)$ then there are $N^*,N^+$ such that
$(M,N^*,a) \in K^{3,*} \bigcap p,\ N \preceq N^+,\ N^* \preceq
N^+$. \item Suppose $K^{3,*}$ has existence and $M\prec N$. Then
there is $M^* \succeq N$ such that $M^*$ is decomposable over $M$
by $K^{3,*}$. Moreover, letting $a \in N-M$, $tp(a,M,N)$ is basic,
we can choose $d_0=a$, where $d_0$ is the element which appears in
Definition \ref{3.6}.
\end{enumerate}
\end{theorem}
\begin{proof}
\mbox{} (1) Suppose $p=tp(M,N,a) \in S^{bs}(M)$. As $K^{3,*}$ is
dense with respect to $\preceq_{bs}$, there are $M^*,N^*,b$ with
$(M,N,a) \preceq_{bs}(M^*,N^*,b)$. As $\frak{s}$ has conjugation,
$p^*=:tp(M^*,N^*,b)$ conjugate to $p$. $K^{3,*}$ is closed under
isomorphisms and so $p \bigcap K^{3,*} \neq \emptyset$.\\ (2)
$K^{3,*}$ has existence and so there are $b,N^*$ such that:
$tp(b,M,N^*)=p,\ (M,N^*,b) \in K^{3,*}$. By the definition of a
type (i.e. the definition of equivalence between triples in a
type), there are a model $N^+$, $N \preceq N^+$ and an embedding
$f:N^* \to N^+$ over $M$ such that $f(b)=a$. Denote
$N^{**}=f[N^*]$. Now as $K^{3,*}$ respects isomorphisms,
$(M,N^{**},a) \in K^{3,*}$. $M \preceq N^{**} \preceq N^+$.\\ (3)
Assume toward a contradiction that $M\prec N$ and there is no
$M^*$ as required. We try to construct
$M_\alpha,a_\alpha,N_\alpha$ by induction on $\alpha \in
\lambda^+$ such that (see the diagram below):
\begin{enumerate}[(a)]
\item $M_0=M,\ N_0=N$. \item $(M_\alpha,M_{\alpha+1},d_\alpha) \in
K^{3,*}$. \item $M_\alpha \preceq N_\alpha$. \item For every
$\alpha \in \lambda^+,\ d_\alpha \in M_{\alpha+1} \bigcap
N_\alpha-M_\alpha$. \item The sequence $\langle
M_\alpha:\alpha<\lambda^+ \rangle$ is increasing and continuous.
\item The sequence $\langle N_\alpha:\alpha<\lambda^+ \rangle$ is
increasing and continuous.
\end{enumerate}

\begin{displaymath}
\xymatrix{
N_0 \ar[r]^{id} & N_1 \ar[rr]^{id} && N_\alpha\\
M_0 \ar[r]^{id} \ar[u]^{id} & M_1 \ar[rr]^{id} \ar[u]^{id} &&
M_\alpha \ar[r]^{id} \ar[u]^{id} &
 M_{\alpha+1} \ni a_\alpha}
\end{displaymath}

We cannot succeed because otherwise substituting the sequences
$\langle M_\alpha:\alpha \in \lambda^+ \rangle,\ \langle
N_\alpha:\alpha \in \lambda^+ \rangle,\ \langle
id_{M_\alpha}:\alpha \in \lambda^+ \rangle$ in Proposition
\ref{1.11} we get a contradiction. So where will we get stuck? For
$\alpha=0$ there is no problem. For $\alpha$ limit take unions. 3
is satisfied by (smoothness) (Definition \ref{1.1}.1.d). What will
we do for $\alpha+1$, (assuming we have defined
$(M_\alpha,N_\alpha,d_\alpha$)? If $N_\alpha=M_\alpha$ then
$N_\alpha$ is decomposable over $M$ by $K^{3,*}$ and the proof has
reached to its end. Otherwise by the existence of the basic types
(\ref{definition good frame}), there is $d_\alpha \in
N_\alpha-M_\alpha$ such that $(M_\alpha,N_\alpha,d_\alpha) \in
K^{3,bs}$ (and for the ``more over'' take $d_0=a$ if $\alpha=0$).
By assumption $K^{3,*}$ has existence, so there are
$d_\alpha^*,M_{\alpha+1}^*$ such that:
$(M_\alpha,M_{\alpha+1}^*,d_\alpha^*) \in K^{3,*},\
tp(d_\alpha^*,M_\alpha,M_{\alpha+1}^*)=tp(d_\alpha,M_\alpha,N_\alpha)$.
By the definition of a type, there are $N_{\alpha+1}$, $N_\alpha
\preceq N_{\alpha+1}$ and an embedding $f:M_{\alpha+1}^* \to
N_{\alpha+1}$ over $M_\alpha$ such that $f(d_\alpha^*)=d_\alpha$.
Denote $M_{\alpha+1}=Im(f)$. We have $N_\alpha \preceq
N_{\alpha+1},\ M_{\alpha+1} \preceq N_{\alpha+1}$ and
$(M_\alpha,M_{\alpha+1},d_\alpha) \in K^{3,*}$. So 2,3,4 are
guaranteed.
\end{proof}

\begin{proposition}[existence of decomposition over two models]
\label{existence of decomposition over two models} \label{3.8} If
$n<2 \Rightarrow M_n \preceq N$ then there is $M^*$ such that: $N
\preceq M^*$ and $M^*$ is decomposable over $M_0$ and over $M_1$.
\end{proposition}
\begin{proof}
Choose an increasing continuous sequence $\langle  M_n:2 \preceq n
\leq \omega \rangle$ such that:
\begin{enumerate}
\item $N \preceq M_2$. \item For every $n \in \omega$, $M_{n+2}$
is decomposable over $M_n$.
\end{enumerate}
The construction is possible by Theorem \ref{3.7}. Now by the
following proposition $M_\omega$ is decomposable over $M_0$ and
$M_1$.
\end{proof}

\begin{proposition} [the decomposable extensions transitivity]
\label{3.9} Let $\langle  M_\epsilon:\epsilon \leq \alpha^*
\rangle$ be an increasing continuous sequence of models, such that
for every $\epsilon<\alpha^*,\ M_{\epsilon+1}$ is decomposable
over $M_\epsilon$. Then $M_{\alpha^*}$ is decomposable over $M_0$.
\end{proposition}
\begin{proof}
Easy.
\end{proof}

\subsection{A disjoint amalgamation}

The following goal is to prove the existence of a disjoint
amalgamation. For this we are going to prove the density of the
reduced triples.

\begin{definition} \label{definition
of disjoint amalgamation} The amalgamation $f_1,f_2,M_3$ of
$M_1,M_2$ over $M_0$ is said to be \emph{disjoint} when $f_1[M_1]
\bigcap f_2[M_2]=M_0$.
\end{definition}

\begin{definition} \label{3.11}
The triple $(M,N,a) \in K^{3,bs}_\lambda$ is \emph{reduced} if
$(M,N,a) \preceq_{bs} (M^*,N^*,a) \Rightarrow M^* \bigcap N = M$.
We define $K^{3,r}:=\{(M,N,a) \in K^{3,bs}:(M,N,a)$ is reduced\}.
\end{definition}

\begin{proposition} \label{3.12}
The reduced triples are dense with respect to $\preceq_{bs}$: For
every $(M,N,a) \in K^{3,bs}_\lambda$ there is a reduced triple
$(M^*,N^*,a)$ which is $\preceq_{bs}$-bigger than it.
\end{proposition}
\begin{proof}
Suppose towards contradiction that over $(M,N,a)$ there is no
reduced triple. We will construct models $M_\alpha,N_\alpha$ by
induction on $\alpha<\lambda^+$ such that:
\begin{enumerate}[(i)]
\item $(M_0,N_0,a)=(M,N,a)$. \item For every $\alpha \in
\lambda^+,\ (M_\alpha,N_\alpha,a)
\preceq_{bs}(M_{\alpha+1},N_{\alpha+1},a)$. \item For every
$\alpha \in \lambda^+,\ M_{\alpha+1} \bigcap N_\alpha \neq
M_\alpha$. \item The sequence $\langle
(M_\alpha,N_\alpha,a):\alpha<\lambda^+ \rangle$ is continuous,
(see Definition \ref{3.1} on page \pageref{3.1}).
\end{enumerate}

Why can we carry out the construction? For $\alpha=0$ see clause
(i) of the construction. For limit $\alpha$ see clause (iv).
Suppose we have defined $\langle M_\beta,N_\beta,a):\beta \leq
\alpha \rangle$. By Proposition \ref{3.2}
$(K^{3,bs},\preceq_{bs})$ is closed under increasing union. So by
clauses (i),(ii),(iv) $(M,N,a) \preceq_{bs}(M_\alpha,N_\alpha,a)$.
So by the assumption $(M_\alpha,N_\alpha,a)$ is not a reduced
triple, i.e. there are $M_{\alpha+1},N_{\alpha+1}$ which satisfies
clauses (ii),(iii). Hence we can carry out this construction.

Now  we have:
\begin{enumerate}
\item The sequences $\langle M_\alpha:\alpha<\lambda^+ \rangle,\
\langle N_\alpha:\alpha<\lambda^+ \rangle$ are increasing (by
clause (ii) and the definition of $\preceq_{bs}$). \item These
sequences are continuous (by clause (iv)). \item For $\alpha \in
\lambda^+,\ M_\alpha \subseteq N_\alpha$ (by the definition of
$K^{3,bs}$). \item For every $\alpha \in \lambda^+,\ M_{\alpha+1}
\bigcap N_\alpha \neq M_\alpha$ (by clause (iii)).
\end{enumerate}
We got a contradiction to Proposition \ref{1.11}.
\end{proof}

\begin{theorem} [The disjoint amalgamation existence theorem]\label{3.13}
Assume that:
\begin{enumerate}
\item $\frak{s}$ has conjugation. \item $M_0,M_1,M_2 \in
K_\lambda$, $M_0 \preceq M_1$ and $M_0 \preceq M_2$.
\end{enumerate}
Then there are $M_3,f$ such that $f:M_2 \to M_3$ is an embedding
over $M_0,\ M_1 \preceq M_3$, and $f[M_2] \bigcap M_1 = M_0$.
Moreover if $a \in M_1-M_0$ and $tp(a,M_0,M_1) \in S^{bs}(M_0)$
then we can add that $tp(a,f[M_2],M_3)$ does not fork over $M_0$.
\end{theorem}

\begin{proof} If $M_1=M_0$ then the theorem is trivial. Otherwise
by the density of basic types (see Definition \ref{definition good
frame}, page \pageref{definition good frame}) there is an element
$a \in M_1-M_0$ such that $tp(a,M_0,M_1) \in S^{bs}(M_0)$. So it
is sufficient to prove the ``moreover''. By Proposition \ref{3.12}
the reduced triples are dense with respect to $\preceq_{bs}$. So
by Theorem \ref{3.7} (the extensions decomposition theorem), as
$\frak{s}$ has conjugation, there is a model $M_1^*$ such that
$M_1 \preceq M_1^*$ and $M_1^*$ is decomposable over $M_1$ by
reduced triples, i.e. there is an increasing continuous sequence
$\langle N_{0,\alpha}:\alpha \leq \delta \rangle$ of models in
$K_\lambda$ such that: $N_{0,0}=M_0,\ M_{0,\delta}=M_1^*$ and
there is a sequence $\langle d_\alpha:\alpha<\delta \rangle$ such
that $(N_{0,\alpha},N_{0,\alpha+1},d_\alpha)$ is a reduced triple
and $d_0=a$. By Proposition \ref{3.3}.1
there is an isomorphism $f$ of $M_2$ over $M_0$ and there is an
increasing continuous sequence $\langle N_{1,\alpha}:\alpha \leq
\delta \rangle$ such that: $N_{0,\alpha} \preceq N_{1,\alpha},\
f[M_2]=N_{1,0}$ and $tp(d_\alpha,N_{1,\alpha},N_{1,\alpha+1})$
does not fork over $N_{0,\alpha}$. So for $\alpha<\delta,\
(N_{0,\alpha},N_{0,\alpha+1},d_\alpha) \allowbreak
\preceq_{bs}(N_{1,\alpha},N_{1,\alpha+1},d_\alpha)$. But the
triple $(N_{0,\alpha},N_{0,\alpha+1},d_\alpha)$ is reduced, so
$N_{1,\alpha} \bigcap \allowbreak N_{0,\alpha+1}=N_{0,\alpha}$.
Hence $N_{1,0} \bigcap N_{0,\delta}=N_{0,0}$ [Why? let $x \in
N_{1,0} \bigcap N_{0,\delta}$. Let $\alpha$ be the first ordinal
such that $x \in N_{0,\alpha}$. $\alpha$ cannot be a limit ordinal
as the sequence is continuous. If $\alpha=\beta+1$ then $x \in
N_{0,\alpha} \bigcap N_{1,\beta}=N_{0,\beta}$, in contradiction to
the definition of $\alpha$ as the first such an ordinal. So we
must have $\alpha=0$, i.e. $x \in N_{0,0}$]. Hence $M_1 \bigcap
f[M_2]=N_{0,0}=N_0$. Denote $M_3=N_{0,\delta}$.
\end{proof}

\section{Uniqueness triples}
\begin{hypothesis}
$\frak{s}$ is a semi-good $\lambda$-frame.
\end{hypothesis}



\begin{definition}\label{equivalent amalgamations}
Suppose
\begin{enumerate} \item
$M_0,M_1,M_2 \in K_\lambda$, $M_0 \preceq M_1 \wedge M_0 \preceq
M_2$. \item For $x=a,b$, $(f_1^x,f_2^x,M_3^x)$ is an amalgamation
of $M_1,M_2$ over $M_0$.
\end{enumerate}
$(f_1^a,f_2^a,M_3^a),(f_1^b,f_2^b,M_3^b)$ are said to be
\emph{equivalent} over $M_0$ if there are $f^a,f^b,M_3$ such that
$f^b \circ f^b_1=f^a \circ f^a_1$ and $f^b \circ f^b_2=f^a \circ
f^a_2$, namely the following diagram commutes:
\begin{displaymath}
\xymatrix{&M_3^b \ar[r]^{f^b} & M_3\\
M_1 \ar[ru]^{f_1^b} \ar[rr]^{f_1^a} && M_3^a \ar[u]_{f^a}\\
M_0 \ar[u]^{id} \ar[r]_{id} & M_2 \ar[ru]_{f_2^a} \ar[uu]^{f_2^b}}
\end{displaymath}
We denote the relation `to be equivalent over $M_0$' between
amalgamations over $M_0$, by $E_{M_0}$.
\end{definition}

\begin{proposition}\label{equivalence proposition}
The relation $E_{M_0}$ is an equivalence relation.
\end{proposition}

\begin{proof}
Assume $(f^a_1,f^a_2,M^a_3)E_{M_0}(f^b_1,f^b_2,M^b_3)$ and
$(f^b_1,f^b_2,M^b_3)E_{M_0}(f^c_1,f^c_2,\allowbreak M^c_3)$. We
have to prove that
$(f^a_1,f^a_2,M^a_3)E_{M_0}(f^c_1,f^c_2,M^c_3)$. Take witnesses
$g_1,g_2,\allowbreak M^{a,b}_3$ for $(f^a_1,f^a_2,M^a_3)E_{M_0}
\allowbreak (f^b_1,f^b_2,M^b_3)$, and witnesses
$g_3,g_4,M^{b,c}_3$ for $(f^b_1,f^b_2,\allowbreak
M^b_3)E_{M_0}(f^c_1,f^c_2,M^c_3)$. Take an amalgamation
$(h_1,h_2,M_3)$ of $M^{a,b}_3$ and $M^{b,c}_3$ over $M^b_3$. Now
we will prove that $h_2 \circ g_4,h_1 \circ g_1,M_3$ witness that
that $(f^a_1,f^a_2,M^a_3)E_{M_0}(f^c_1,f^c_2,M^c_3)$, i.e. to
prove that the following diagram commutes:
\begin{displaymath}
\xymatrix{&M_3^c \ar[r]^{h_2 \circ g_4} & M_3\\
M_1 \ar[ru]^{f_1^c} \ar[rr]^{f_1^a} && M_3^a \ar[u]_{h_1 \circ g_1}\\
M_0 \ar[u]^{id_{M_0}} \ar[r]_{id_{M_0}} & M_2 \ar[ru]_{f_2^a}
\ar[uu]^{f_2^c}}
\end{displaymath}
i.e. to prove that the following diagram commutes:
\begin{displaymath}
\xymatrix{&M_3^c \ar[r]^{g_4} & M^{b,c}_3 \ar[r]^{h_2} & M_3\\
&& M_3^b \ar[u]^{g_3} \ar[r]^{g_2} & M_3^{a,b} \ar[u]^{h_1}\\
M_1 \ar[ruu]^{f_1^c} \ar[rru]^{f_1^b} \ar[rrr]^{f_1^a} &&& M_3^a
\ar[u]^{g_1}\\
M_0 \ar[u]^{id} \ar[r]^{id} & M_2 \ar[uuu]^{f_2^c}
\ar[uur]^{f_2^b} \ar[urr]^{f_2^a}}
\end{displaymath}

$(h_2 \circ g_4) \circ f_2^c=h_2 \circ g_3 \circ f_2^b=h_1 \circ
g_2 \circ f_2^b=(h_1 \circ g_1) \circ f_2^a$ and similarly $(h_2
\circ g_4) \circ f_1^c=(h_1 \circ g_1) \circ f_1^a$.
\end{proof}

\begin{definition} \label{4.1}\label{definition of uniqueness triples} $K^{3,uq}=K^{3,uq}_\frak{s}$ is the class of
triples $(M_0,N_1,a) \in K^{3,bs}$ such that if $M_0 \preceq N_1
\in K_\lambda$ then up to equivalence over $M_0$ there is a unique
amalgamation $(f_1,f_2,N_1)$ of $M_1,N_0$ over $M_0$ such that
$tp(f_1(a),f_2[N_0],N_1)$ does not fork over $M_0$. Equivalently,
if for $n=1,2$ $(M,N,a) \preceq_{bs} (M^*_n,N^*_n,a)$ and $f:M^*_1
\to M^*_2$ is an isomorphism over $M$, then for some $f_1,f_2,N^*$
the following hold: $f_n:N^*_n \to N^*$ is an embedding over $N$,
and $f_1\restriction M^*_1=f_2\restriction M^*_2 \circ f$.A
\emph{uniqueness triple} is a triple in $K^{3,uq}$.
\end{definition}

\begin{proposition} \label{4.2}
\mbox{}
\begin{enumerate}
\item If $p_0,p_1$ are conjugate types and in $p_0$ there is a
uniqueness triple, then also in $p_1$ there is such a triple.
\item If $\frak{s}$ has conjugation, then every uniqueness triple
is reduced.
\end{enumerate}
\end{proposition}

\begin{proof}
\mbox{}
\begin{enumerate}
\item Suppose $p_0=tp(a,M,N),\ (M,N,a) \in K^{3,uq}$. Let $f$ be
an isomorphism with domain $M$, such that $f(p_0)=p_1$.
$K,\preceq$ are closed under isomorphisms, so it is easy to prove
that $(f[M],f^+[N],f^+(a)) \in K^{3,uq}$, where $f \subseteq f^+,\
dom(f^+)=N$. But $(f[M],f^+[N],f^+(a)) \in p_1$. \item Suppose
$(M_0,N_0,a) \in K^{3,uq}$ and $(M_0,N_0,a) \preceq_{bs}
(M_1,N_1,a)$. By Theorem \ref{3.13} (the existence of a disjoint
amalgamation), there are $f,N_2$ such that $f:M_1 \to N_2$ is an
embedding over $M_0,\ N_0 \preceq N_2,\ f[M_1] \bigcap N_0=M_0$
and $tp(a,f[M_1],N_2)$ does not fork over $M_0$. By Definition
\ref{4.1}, there are $f_1,f_2,N^*$ such that: $f_n:N_n \to N^*$
and embedding over $N_0$ and $f_1\restriction M_1=f_2 \circ f$.
Let $x \in M_1-M_0$. Then $x \notin N_0$ [Why? otherwise $f(x) \in
f[M_1]-M_0$, so $f(x) \notin N_0$, so $f_1(x)=f_2(f(x))\notin N_0$
and hence $x \notin N_0$].
\end{enumerate}
\end{proof}

\begin{definition} \label{4.3}
Let $\frak{s}$ be a semi good $\lambda$-frame.
\begin{enumerate}
\item $\frak{s}$ is \emph{weakly successful in the sense of
density}, if $K^{3,uq}$ is dense with respect to $\preceq_{bs}$.
\item $\frak{s}$ is \emph{weakly successful} if $K^{3,uq}$ has
existence.
\end{enumerate}
\end{definition}

\begin{proposition} \label{4.4}
\mbox{}
\begin{enumerate}
\item If $\frak{s}$ is weakly successful in the sense of density
and it has conjugation then it is weakly successful. \item Let
$\frak{s}$ be weakly successful. If $p=tp(a,M,N) \in S^{bs}(M)$,
then there is a model $N^*$ such that $(M,N^*,a) \in K^{3,uq}
\bigcap p$.
\end{enumerate}
\end{proposition}
 \begin{proof}
\mbox{}
\begin{enumerate}
\item Substitute $K^{3,*}:=K^{3,uq}$ in Theorem \ref{3.7}.1 (page
\pageref{3.7}). \item By Theorem \ref{3.7}.2.
\end{enumerate}
 \end{proof}

Now the reader can believe that the assumption that $\frak{s}$ is
weakly successful is reasonable and jump to section 5 or to read
the rest of this section (which is based on \cite{sh838}).

\begin{hypothesis} \label{4.5}
$\frak{s}$ is (a semi-good $\lambda$-frame and) not weakly
successful in the sense of density.
\end{hypothesis}

\discussion{Discussion toward defining nice construction frame:}
Every model $M \in K_{\mu^+}$ can be represented as $\bigcup
\{M_{\beta}:\beta<\mu^+\}$ where each $M_\beta$ is in $K_\mu$ (and
the sequence is increasing and continuous). Now we can represent
each $M_\beta$ as $\bigcup \{M_{\alpha,\beta}:\alpha<\mu\}$ where
each $M_{\alpha,\beta}$ is in $K_{<\mu}$. So we can approximate a
model $M$ in $K_{\mu^+}$ by a ``rectangle''
$\{M_{\alpha,\beta}:\alpha<\mu,\beta<\mu^+\}$ of models in
$K_{<\mu}$, where $\langle M_{\alpha,\beta}:\alpha<\beta \rangle$
is an increasing continuous sequence of models in $K_{<\mu}$,
$\langle \bigcup \{M_{\alpha,\beta}:\alpha<\mu \}:\beta<\mu^+
\rangle$ is an increasing continuous sequence of models in
$K_{\mu}$ and $\bigcup \{M_{\alpha,\beta}:\alpha<\mu,\beta<\mu^+
\}=M$.

Now we want to violate this rectangle. For $n=1,2$ we will define
a relation $FR_n$ such that $(\forall \alpha,\beta)
[(M_{\alpha,\beta},M_{\alpha+1,\beta},\allowbreak
I_{\alpha,\beta}) \in FR_1 \wedge
(M_{\alpha,\beta},M_{\alpha,\beta+1},J_{\alpha,\beta}) \in FR_2$,
where $I_{\alpha,\beta}$ and $J_{\alpha,\beta}$ are witnesses for
the extensions, namely $I_{\alpha,\beta} \subset
M_{\alpha+1,\beta}-M_{\alpha,\beta}$ and $J_{\alpha,\beta} \subset
M_{\alpha,\beta+1}-M_{\alpha,\beta}$. So essentially, $FR_n$ is a
relation on extensions.

We have to violate also the pairs of such pairs, i.e.
$((M_{\alpha,\beta},M_{\alpha+1,\beta}),(M_{\alpha,\beta+1},\allowbreak
M_{\alpha+1,\beta+1}))$. In other words, we have to define
2-dimensional relations $\leq_1,\leq_2$ on $FR_1,FR_2$
respectively.

\begin{definition} \label{4.7}
$U=(\mu,k^u,FR_1,FR_2,\leq_1,\leq_2)$ is a nice construction frame
if:
\begin{enumerate}
\item $\aleph_0<\mu$ is a regular cardinal. \item
$k^U=(K^U,\preceq^U)$ is an a.e.c. in $<\mu$. The vocabulary of
$K^U$ will denoted $\tau^U$. \item For $n=1,2$ $FR^n$ is a class
of triples $(M,N,J)$ such that:
\begin{enumerate}
\item $M,N \in K^U,\ M\preceq^UN,\ J\subseteq N-M$.
    \item For every $M \in K^U$ there are $N,J$ such that: $J \neq \emptyset$
  and $(M,N,J) \in FR_n$.
    \item If $M \preceq^U N,$ then $(M,N,\emptyset) \in FR_n$.
    \end{enumerate}
    \item ``$(FR_n,\leq_n)$ satisfies some axioms of a.e.c. and
disjointness'':
\begin{enumerate}
\item $\leq_n$ is an order relation of $FR_n$. \item The relations
$FR_n,\leq_n$ are closed under isomorphisms. \item If
$(M_{0,0},M_{0,1},J_0) \leq_n (M_{1,0},M_{1,1},J_1)$ then
  $(n_1 \leq n_2<2 \wedge m_1 \leq m_2<2) \Rightarrow M_{n_1,m_1}
  \preceq^U M_{n_2,m_2}$ and $M_{1,0} \bigcap M_{0,1}=M_{0,0}$.

  \item Axiom c of a.e.c.: For every $\delta<\mu$ and an
  $\leq_n$-increasing continuous sequence
  $\langle (M^\alpha,N^\alpha,J^\alpha):\alpha<\delta \rangle$ we have\\
  $(M^0,N^0,J^0) \leq_n(\bigcup \{M^\alpha:\alpha<\delta\},\bigcup
  \{N^\alpha:\alpha<\delta \},\bigcup \{J^\alpha:\alpha<\delta
\  \})$.
\end{enumerate}
\item $U$ has disjoint amalgamation (at first glance one can think
that the disjointness is in the assumption, but it is in the
conclusion, see 4c):\\
If $(M_0,M_1,J_1) \in FR_1$, $(M_0,M_2,J_2) \in FR_2$ and $M_1
\bigcap M_2=M_0$ then there are $M_3,J_1^*,J_2^*$ such that for
$n=1,2$ $M_n \preceq^U M_3$ and $(M_0,M_n,J_n) \leq_n
(M_{3-n},M_3,J_n^*)$.
\end{enumerate}
\end{definition}

A way to force an amalgamation to be disjoint, is to replace the
equality relation by an equivalence one. This is the role of $E$
in the following definition.
\begin{definition} \label{4.8}
Let $U$ be a nice construction frame. Let $(K,\preceq)$ be an
a.e.c. with a vocabulary $\tau$, such that $\tau \subseteq \tau^U$
and there is a 2-place predicate $E \in \tau^U-\tau$ (in the main
case $\tau^U=\tau \bigcup \{E\}$), such that for $M \in K^U$ we
have:
\begin{enumerate}
\item $E^M$ is an equivalence relation. \item If $R$ is a
predicate in $\tau^U$ different from = and $xE^My$ then
$R^M(x_0,... ,\allowbreak x_{i-1},\allowbreak x,x_{i+1},...x_n)$
iff $R^M(x_0,... ,x_{i-1},y,x_{i+1},...x_n)$.
\end{enumerate}
Similarly for function symbols.\\
We write $(K,\preceq)=(U/E)^\tau$ when:\\
$(K,\preceq)$ is an a.e.c. and $K_{<\mu}=\{N:(\exists M \in
K^U)(N=M/E)\}$, where $M/E$ is defined by the following way: Its
world is the set of equivalence classes of $E^M$, its vocabulary
is $\tau$ and it interprets the predicates and function symbols by
representatives of the equivalence classes.
\end{definition}

Now we are going to define approximations of cardinality $\mu$, by
the approximations of cardinality $<\mu$.
\begin{definition} \label{4.9}
\mbox{}
\begin{enumerate}
\item $K^{qt}=K^{qt,U}:=\{(\bar{M},\bar{J},f):\bar{M}= \langle
M_\alpha:\alpha<\mu \rangle,\bar{J}=\langle J_\alpha:\alpha<\mu
\rangle,f \in ^\mu\mu, \alpha<\mu \Rightarrow
(M_\alpha,M_{\alpha+1},J_\alpha) \in FR_2\}$ ($f$ plays a role in
the relation $\leq^{qt}$). \item $\leq^{qt}$ is a relation on
$K^{qt}$. $(M_0,J_0,f_0) \leq (M_1,J_1,f_1)$ iff there is a club
$E$ of $\mu$ such that for every $\delta \in E$ and $\alpha \leq
f^1(\delta)$ we have:
\begin{enumerate}
\item $f^1(\delta) \leq f^2(\delta)$. \item $M_{0,\delta+1}\leq
M_{1,\delta+1}$. \item
$(M_{0,\delta+\alpha},M_{0,\delta+\alpha+1},J_{0,\delta+\alpha})
\leq_2
(M_{1,\delta+\alpha},M_{1,\delta+\alpha+1},J_{1,\delta+\alpha})$.
\item $M_{1,\delta+\alpha} \bigcap \bigcup
\{M_{0,\epsilon}:\epsilon<\mu \}=M_{0,\delta+\alpha}$.
\end{enumerate}
\end{enumerate}
\end{definition}

\begin{definition} \label{4.11} We say that \emph{almost every}
$(\bar{M},\bar{J},f) \in K^{qt}$ satisfies the property $pr$ when:
There is a function $g:K^{qt} \rightarrow K^{qt}$ such that if
$\langle \bar{M}^\alpha:\bar{J}^\alpha,f^\alpha \rangle$ is an
$\leq^{qt}$-increasing continuous (in the sense which is defined
in \cite{sh838} and not here) and $sup\{\alpha \in
\delta:g((\bar{M}^\alpha,\bar{J}^\alpha,f^\alpha))=
(\bar{M}^{\alpha+1},\bar{J}^{\alpha+1},f^{\alpha+1} \allowbreak
)\}=\delta)$, then $(\bar{M}^\delta,\bar{J}^\delta,f^\delta) \in
pr$.
\end{definition}

\begin{definition} \label{4.12}
\mbox{}
\begin{enumerate}
\item Let $U$ be a nice construction frame. We say that $U$
satisfies the \emph{weak coding property} for $(K,\preceq)$ if
almost every $(\bar{M},\bar{J},f) \in K^{qt}$ satisfies the weak
coding property. \item We say that $(\bar{M},\bar{J},f) \in
K^{qt}$ satisfies the \emph{weak coding property} when: There are
$\alpha_0<\mu$ and $N_0,J_0$ such that $(M_{\alpha_0},N_0,J_0) \in
FR_1,\ N_0 \bigcap M=M_{\alpha_0}$ where $M:=\bigcup
\{M_\alpha:\alpha<\mu\}$, and there is a club $E$ of $\mu$ such
that for every $\alpha_1 \in E$ and every $N_1,J_1$, which satisfy
$(M_{\alpha_0},N_0,J_0) \leq _1(M_{\alpha_1},N_1,J_1) \wedge N_1
\bigcap M=M_{\alpha_1}$, there is $\alpha_2 \in (\alpha_1,\mu)$
and for $n=1,2$ there are $N_{2,n},J_{2,n}$ such that:
\begin{enumerate}
\item $(M_{\alpha_1},N_1,J_1) \leq_1
(M_{\alpha_2},N_{2,n},J_{2,n})$. \item $N_{2,1},N_{2,2}$ are
incomparable amalgamations of $M_{\alpha_2},N_1$ over
$M_{\alpha_1}$, i.e. there are no $N,f_1,f_2$ such that $f_n$ is
an embedding of $N_{2,n}$ into $N$ over $N_1 \bigcup
M_{\alpha_2}$.
\end{enumerate}
\end{enumerate}
\end{definition}

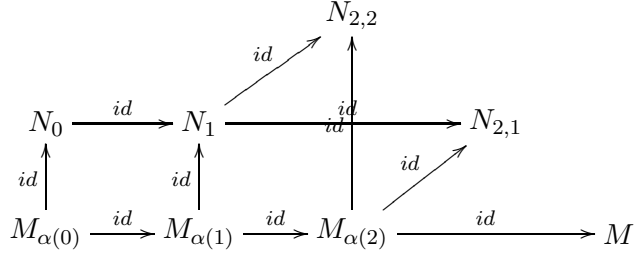
\begin{displaymath}
\xymatrix{
&&N_{2,2}\\
N_0 \ar[r]^{id} & N_1 \ar[rr]^{id} \ar[ru]^{id} && N_{2,1}\\
M_{\alpha(0)} \ar[r]^{id} \ar[u]^{id} & M_{\alpha(1)} \ar[r]^{id}
\ar[u]^{id} & M_{\alpha(2)} \ar[rr]^{id} \ar[uu]^{id} \ar[ru]^{id}
&& M
 }
\end{displaymath}

\begin{definition}
$\mu_{unif}(\mu^+,2^\mu):=Min \{|P|:P$ is a family of subsets of
$^{\mu^+}(2^{\mu})$
with union $^{\mu^+}(2^{\mu})$ and for each $A \in P$ there is a
function $c$ with domain $\bigcup \{^\alpha
(2^\mu):\alpha<\mu^+\}$ such that for each $f \in A$, the set
$\{\delta \in \mu^+:c(f \restriction \delta)=f(\delta)\}$ is not
stationary\}.
\end{definition}

\begin{proposition}\label{unifleq}
$\mu_{unif} \leq 2^{\mu^+}$.
\end{proposition}

\begin{proof}
Easy.
\end{proof}

\begin{remark}
$\mu_{unif}(\mu^+,2^\mu)$ is ``almost $2^{\mu^+}$'': If
$\beth_\omega \leq \mu$, then $\mu_{unif}(\mu^+,2^\mu) \allowbreak
=2^{\mu^+}$, and in any case it is not clear if
$\mu_{unif}(\mu^+,2^\mu)<2^{\mu^+}$ is consistent. There are
propositions which say that it is ``a big cardinal''.
\end{remark}

The following theorem is written in \cite{sh838}, and we bring it
without a proof.
\begin{theorem} \label{4.13}
Let $U$ be a nice construction frame which satisfies the weak
coding property for $(K,\preceq)$. Suppose the following set
theoretical assumptions:
\begin{enumerate}
\item $2^\theta=2^{<\mu}<2^\mu$. \item $2^\mu<2^{\mu^+}$. \item
The ideal $WdmId(\mu)$ is not saturated in $\mu^+$.
\end{enumerate}

Then $\mu_{unif}(\mu^+,2^\mu) \leq I(\mu^+,K)$, where $I(\mu^+,K)$
is the number of non-isomorphic models in $K_{\mu^+}$.
\end{theorem}

Now we are going to study a specific nice construction frame.
\relax From now $(K,\preceq)$ will denote the a.e.c. of
$\frak{s}$.
\begin{definition} \label{4.14}
Define $U=(\mu,(K^u,\preceq^u),FR_1,FR_2,\leq_1,\leq_2)$:
\begin{enumerate}
\item $\mu=\lambda^+$. \item The vocabulary of $K^U$ is
$\tau^U:=\tau\bigcup \{E\}$ where $E$ is a new predicate. \item
$K^U:=\{M:||M||=\lambda,\ M/E \in K_\lambda\}$. ($M/E$ is well
defined only if $E^M$ is a congruence relation on $|M|$, see
Definition \ref{4.8}. So if not, then $M$ does not belong to
$K^U$). \item $\preceq^U:=\{(M,N):M/E \preceq N/E \wedge M
\subseteq N \}$. \item $FR_n:=\{(M,N,J):M,N \in K^U,\ J \neq
\emptyset \Rightarrow (\exists a)[J=\{a\}\wedge (M/E,N/E,a/E) \in
K^{3,bs}]\}$. \item For $n=1,2$ the relation $\leq_n$ is defined
by the relation $\preceq_{bs}$ in the same way we defined $FR_n$.
\end{enumerate}
\end{definition}

\begin{proposition} \label{4.15}
Almost every $(\bar{M},\bar{J},f) \allowbreak \in K^{qt,U}$
satisfies: $\bigcup \{M_\alpha/E:\alpha<\lambda^+\}$ is a
saturated model in $\lambda^+$ over $\lambda$.
\end{proposition}

\begin{proof}
See \cite{sh838}.
\end{proof}

\begin{theorem} \label{4.16}
If $\bar{M}=\langle  M_\alpha:\alpha<\lambda^+ \rangle,\
\bar{a}=\langle a_\alpha:\alpha<\lambda^+ \rangle,\
(\bar{M},\bar{a},f) \in K^{qt}$ and $\bigcup
\{M_\alpha/E:\alpha<\lambda^+\}$ is saturated in $\lambda^+$ over
$\lambda$, then $(\bar{M},\bar{a},f)$ satisfies the weak coding
property.
\end{theorem}

\begin{proof}
For distinguishing between models in $K_\lambda$ to models in
$K^U$, we add to the names of models in $K_\lambda$ subscript $e$,
unless they are written in the form $M/E$. For example: $M_e,\
M_{2,e}$. Similarly for isomorphisms.

\begin{lemma}
\mbox{}
\begin{enumerate} \item Let $N_0 \in K^U,\
N_{1,e} \in K_\lambda$ be such that $N_0/E \preceq N_{1,e}$. Then
there is $N_1 \in K^U$ such that: \begin{enumerate} \item
$N_1/E=N_{1,e}$. \item $N_0 \preceq^U N_1$. \item $N_1$ is
embedded in every model which satisfies 1,2.
\end{enumerate} In this case we call $N_1$ the
canonical completion of $N_0,N_{1,e}$. There is exactly one such a
model up to isomorphism. Clearly every $[x] \in N_1-N_0$ is a
singleton. \item Suppose:
\begin{enumerate}
\item $N_0 \preceq ^U N_1,\ N_0 \preceq ^U N_2$. \item $g_e:N_1/E
\to N_2/E$ is an embedding over $N_0/E$. \item $N_1$ is the
canonical completion of $N_1/E,N_0$.
\end{enumerate} Then there is an embedding $g:N_1 \to N_2$
over $N_0$ such that $(\forall x \in N_1)(g(x) \in [g_e(x/E)])$.
\item Suppose for $n<3,\ N_n \in K^U$, $N_0/E \preceq N_n/E
\preceq N_{3,e} \in K_\lambda$ and $N_1 \bigcap N_2=N_0$. Then
there is $N_3 \in K^U$ such that $N_3/E=N_{3,e}$ and for $n=1,2 \;
N_n \preceq N_3$.
\end{enumerate}
\end{lemma}

\begin{proof}
\mbox{}
\begin{enumerate}
\item Trivial. \item Use the axiom of choice [For $x \in N_1-N_0\
g(x)$ choose an arbitrary element in $g_e([x])$]. \item Trivial.
\end{enumerate}
\end{proof}

Now we prove that $(\bar{M},\bar{a},f)$ satisfies the weak coding
property, by the following steps:\\
\step{Step a:} Denote $\alpha(0)=0$. $M_0/E \in K_\lambda$. So by
the categoricity in $K_\lambda$ and non-weak successfulness, there
are $N_{0,e} \in K_\lambda$ and $a \in N_{0,e}$ such that
$(M_0/E,N_{0,e},a) \in K^{3,bs}$ and every triple which is
$\preceq_{bs}$-bigger from it is not a uniqueness triple. Without
lose of generality $N_{0,e} \bigcap M/E=M_0/E$. Let $N_0 \in K^U$
be the model with world $N_{0,e}$, $E^{N_0}$ is the equality, and
$N_0/E=N_{0,e}$. $\lambda^+$ is of course a club of $\lambda^+$.
Let $\alpha(1) \in (\alpha(0),\mu)$, and let $N_1 \in K^U$ such
that $N_1 \bigcap M=M_{\alpha(1)},\ (M_0,N_0,a) \leq_n
(M_{\alpha(1)},N_1,a)$. We have to find $\alpha(2)$.\\
\step{Step b:} $(M_{\alpha(1)}/E,N_1/E,a)$ is not a uniqueness
triple. So for \ $n<2$ there are $M_{2,n,e},N_{2,n,e}^* \in
K_\lambda$ and an isomorphism $g_e:M_{2,0,e} \to M_{2,1,e}$ over
$M_{\alpha(1)}/E$ such that $(M_{\alpha(1)}/E,N_1/E,a)
\preceq_{bs} (M_{2,n,e},N_{2,n,e}^*,a)$ and there are no
$g_{0,e},g_{1,e},N_{3,e}$ such that $g_{n,e}:N_{2,n,e}^* \to
N_{3,e} \in K_{\lambda}$ an embedding over $N_1/E$ and $g_{1,e}
\circ g_e=g_{0,e}$. We choose new elements for
$N_{2,n,e}^*-(M_{\alpha(1)}/E)$, i.e. without loss of generality
$M/E \bigcap N_{2,n,e}^*=M_{\alpha(1)}/E$. By item 1 in the lemma
for $n<2$ there is a model $M_{2,n}$ which is canonical over
$M_{\alpha(1)},M_{2,n,e}$. By item 3 of the lemma for $n<2$ there
is a model $N_{2,n}^* \in K^U$ such that $M_{2,n} \preceq^U
N_{2,n}^*,\ N_1 \preceq N_{2,n}^*$ and $N_{2,n}^*/E=N_{2,n,e}^*$.
\begin{displaymath}
\xymatrix{N_0 \ar[r]^{id} & N_1 \ar[r]^{id} & N_{2,n,e}^*\\
M_0 \ar[r]^{id} \ar[u]^{id} & M_{\alpha(1)} \ar[r]^{id}
\ar[u]^{id} & M_{2,n,e} \ar[u]^{id} }
\end{displaymath}
\step{Step c:} $M/E$ is saturated in $\lambda^+$ over $\lambda$,
so by Lemma \ref{1.12} (the saturation = model homogeneity lemma),
there is an embedding $f_{0,e}:M_{2,0,e} \to M/E$ over
$M_{\alpha(1)}/E$. So by item b of the lemma over, there is an
embedding $f_0:M_{2,0} \to M$ over $M_{\alpha(1)}$. Define
$f_1=f_0 \circ g_e^{-1}$. Now for $n<2$ the function $f_n:M_{2,n}
\to M$ is an embedding.
\begin{displaymath}
\xymatrix{& N_{3,e}\\
M_{2,0,e} \ar[ru]^{g_{0,e}} \ar[rr]^{g_e} && M_{2,1,e}
\ar[lu]_{g_{1,e}}}
\end{displaymath}
\step{Step d:} For $n<2$ let $h_n$ be a function with domain
$N_{2,n}^*$ that extends $f_n$ by the identity. So
$h_n\restriction N_1$ is the identity.
\begin{displaymath}
\xymatrix{N_1 \ar[r]^{id} & h_n[N_{2,n}^*]\\
M_{\alpha(1)} \ar[r]^{id} \ar[u]^{id} & f_n[M_{2,n}] \ar[r]^{id}
\ar[u]^{id} & M }
\end{displaymath}
\step{Step e:} Define $\alpha(2):=\Min\{\alpha \in
\lambda^+:f_0[M_{2,0}] \preceq
M_{\alpha(2)}\}$.\\
\step{Step f:} For $n<2$ we can choose a model $N_{2,n} \in K^U$
such that $(f_n[M_{2,n}],\allowbreak h_n[N_{2,n}^*],a) \preceq_1
(M_{\alpha(2)},N_{2,n},a)$.
\begin{displaymath}
\xymatrix{N_1 \ar[r]^{id} & h_n[N_{2,n}^*] \ar[r]^{id} & N_{2,n}\\
M_{\alpha(1)} \ar[r]^{id} \ar[u]^{id} & f_0[M_{2,0}] \ar[r]^{id}
\ar[u]^{id} & M_{\alpha(2)} \ar[u]^{id} }
\end{displaymath}

By the transitivity of the relation $\leq_1$, we have
$(M_{\alpha(1)},N_1,a) \leq_1 (M_{\alpha(2)},\allowbreak N_{2,n},a)$.\\
\step{Step g:} $N_{2,0},N_{2,1}$ witness that $\alpha(2)$ is as
required [Toward contradiction assume that there are $N_{3,e} \in
K_\lambda$ and embeddings $g_{0,e},g_{1,e}$ such that
$g_n:N_{2,n}/E \to N_3$ is an embedding over
$M_{\alpha(2)}/E\bigcup N_1/E$ Define an isomorphism
$g_{n,e}^*:N_{2,n,e}^* \to N_{3,e}$ by
$g_{n,e}^*(x):=g_{n,e}([h_n(x)])$. This is an embedding over
$N_1/E$ and it includes $f_{n,e}$. This contradict the way we
chose $M_{2,n,e},N_{2,n,e}^*$ in step b]. Hence the triple
$(\bar{M},\bar{a},f)$ satisfies the weak coding property.
\end{proof}

\begin{corollary} \label{4.17}
${U}$ satisfies the weak coding property.
\end{corollary}

\begin{proof}
By \ref{4.15},\ref{4.16}.
\end{proof}

\begin{corollary} \label{4.18}
Let $\frak{s}$ be a semi-good $\lambda$-frame which is not weakly
successful in the sense of density. Then $I(\lambda^{+2},K) \geq
\mu_{unif}(\lambda^{+2},2^{\lambda^+})$.
\end{corollary}
\begin{proof}
By \ref{4.13},\ref{4.17}.
\end{proof}

\section{Non-forking amalgamation}
\begin{hypothesis} \label{5.0}
$\frak{s}$ is a weakly successful semi-good $\lambda$-frame with
conjugation, but we do not use local character in this section.
\end{hypothesis}

\subsection{The axioms of non-forking amalgamation}

\discussion{Introduction:} We want to find a relation of a
canonical amalgamation (see the discussion in the beginning of
section 3). In Definition \ref{5.1} we define the properties this
relation has to satisfy.

\begin{definition} \label{5.1}
Let $NF \subseteq \ ^4(K_\lambda)$ be a relation. We say
$\bigotimes_{NF}$ when the following axioms are satisfied:
\begin{enumerate}[(a)]
\item If $NF(M_0,M_1,M_2,M_3) \ then \ n\in \{1,2\}\Rightarrow
M_0\preceq M_n\preceq M_3 \ and \ M_1 \cap M_2=M_0$.  \item
Monotonicity: If $NF(M_0,M_1,M_2,M_3) \ and \ N_0=M_0, n<3
\Rightarrow N_n\preceq M_n\wedge N_0\preceq N_n\preceq N_3,
(\exists N^{*})[M_3\preceq N^{*}\wedge N_3\preceq N^{*}] \ then \
NF(N_0 \allowbreak ,N_1,N_2,N_3)$. \item Existence: For every
$N_0,N_1,N_2 \in K_\lambda$ if $l\in \{1,2\}\Rightarrow N_0
\preceq N_l$ and $N_1\bigcap N_2=N_0$ then there is $N_3$ such
that $NF(N_0,N_1,N_2,N_3).$ \item Weak uniqueness: Suppose for
$x=a,b$ $NF(N_0,N_1,N_2,N^{x}_3)$. Then there is a joint embedding
of $N^a_3,N^b_3$ over $N_1 \bigcup N_2$. \item Symmetry:
$NF(N_0,N_1,N_2,N_3) \Leftrightarrow NF(N_0,N_2,N_1,N_3)$. \item
Long transitivity: For $x=a,b$ let $\langle M_{x,i}:i\leq \alpha^*
\rangle$ an increasing continuous sequence of models in
$K_\lambda$. Suppose $i<\alpha^* \Rightarrow
NF(M_{a,i},M_{a,i+1},M_{b,i},\allowbreak M_{b,i+1})$. Then
$NF(M_{a,0},M_{a,\alpha^{*}},M_{b,0},M_{b,\alpha^{*}})$
\end{enumerate}
\end{definition}

We give another version of weak uniqueness:
\begin{proposition} \label{remark about
uniqueness} Suppose \begin{enumerate} \item $\bigotimes_{NF}$.
\item $NF(M_0,M_1,M_2,M_3)$ and $NF(M_0,M_1^*,M_2^*,M_3^*)$. \item
For $n=1,2$ there is an isomorphism $f_n:M_n \to M_n^*$ over
$M_0$.
\end{enumerate}
Then there are $M,f$ such that:
\begin{enumerate}
\item For $n<3$ $f\restriction M_n=f_n$. \item $M_3^* \preceq M$.
\item $f[M_3] \preceq M$.
\end{enumerate}
\end{proposition}

\begin{proof}
$M_1 \bigcap M_2=M_0$, so there is a function $g$ with domain
$M_3$ such that $f_1 \bigcup f_2 \subseteq g$. So $g[M_1]=M_1^*$
and $g[M_2]=M_2^*$. Hence $NF(M_0,M_1^*,M_2^*,\allowbreak g[M_3])$
and $NF(M_0,M_1^*,M_2^*,M_3^*)$. Therefore we can use the weak
uniqueness in Definition \ref{5.1}.
\end{proof}

Roughly speaking the following proposition says that finding a
relation $NF$ that satisfies clauses a,c,d of Definition \ref{5.1}
is equivalent to assigning to each triple $(M_0,M_1,M_2) \in
D:=\{(M_0,M_1,M_2):M_0,M_1,M_2 \in K_\lambda,\ M_0 \preceq M_1,M_0
\preceq M_2\}$ a disjoint amalgamation (see Definition
\ref{definition of disjoint amalgamation}) $(f_1,f_2,M_3)$ of
$M_1,M_2$ over $M_0$ up to $E_{M_0}$ (see Definition
\ref{equivalent amalgamations}.

\begin{proposition}\label{explanation to definition 5.1}
Let $NF$ be a relation that satisfies clauses a,c,d of Definition
\ref{5.1}. Then:
\begin{enumerate}
\item There is a function $G$ with domain
$D:=\{(M_0,M_1,M_2):M_0,M_1,M_2 \in K_\lambda,\ M_0 \preceq
M_1,M_0 \preceq M_2\}$ which assign to each triple $(M_0,M_1,M_2)$
an amalgamation $(f_1,f_2,M_3)$ of $M_1,M_2$ over $M_0$, such that
$NF(M_0,f_1[M_1],\allowbreak f_2[M_2],M_3)$ (in this item we do
not use clause d). \item If $G_1,G_2$ are two functions as in item
1 (with respect to $NF$), then for every $(M_0,M_1,M_2) \in D$,
$G_1((M_0,M_1,M_2))E_{M_0}G_2((M_0,M_1,M_2))$. \item If $G$ is a
function with domain $D:=\{(M_0,M_1,M_2):M_1,M_1,M_2 \in
K_\lambda,\ M_0 \preceq M_1,M_0 \allowbreak \preceq M_2\}$ which
assign to each triple $(M_0,M_1,M_2)$ a disjoint amalgamation,
then the relation $R:=\{(M_0,M_1,M_2,M_3):M_1 \bigcap M_2=M_0,\
G((M_0,M_1,M_2))E_{M_0}(id_{M_1},id_{M_2},M_3)$\} satisfies
clauses a,c,d of Definition \ref{5.1}.
\end{enumerate}
\end{proposition}

\begin{proof}
We leave to the reader.
\end{proof}

\begin{definition}\label{respect}
Suppose $\bigotimes_{NF}$. $NF$ is said to \emph{respect the
frame} $\frak{s}$ when: if $NF(M_0,M_1,M_2,M_3)$ and
$tp(a,M_0,M_1) \in S^{bs}(M_0)$ then $tp(a,M_2,M_3)$ does not fork
over $M_0$.
\end{definition}

\subsection{The relation $NF$}
First we define a relation $NF^*$ and then we define a relation
$NF$ as its monotonicity closure. Theorem \ref{5.13} asserts that
the relation $NF$ is the unique relation which satisfies
$\bigotimes_{NF}$ and respects the frame $\frak{s}$.

\begin{definition} \label{5.2}
Define a relation $NF^{*}=NF^{*}_\lambda$ on $^4(K_\lambda)$ by:
$NF^{*}(N_0, \allowbreak N_1,N_2 ,N_3)$ if there is
$\alpha^{*}<\lambda^+$ and for l=1,2 there are an increasing
continuous sequence $\langle  N_{l,i}:i\leq \alpha^{*} \rangle$
and a sequence $\langle d_i:i<\alpha^{*} \rangle$ such that:
\begin{displaymath}
\xymatrix{
N_2=N_{2,0} \ar[rr]^{id} &&N_{2,i} \ar[r]^{id} & N_{2,i+1} \ar[rr]^{id} && N_{2,\alpha^*}=N_3\\
N_0=N_{1,0} \ar[rr]^{id} \ar[u]^{id} &&N_{1,i} \ar[r]^{id}
\ar[u]^{id} & N_{1,i+1} \ar[u]^{id} \ar[rr]^{id} &&
N_{1,\alpha^*}=N_1 \ar[u]^{id} &&
 }
\end{displaymath}

\begin{enumerate}[(a)]
\item $n<3\Rightarrow N_0\preceq N_n \preceq N_3$. \item
$N_{1,0}=N_0, N_{1,\alpha^{*}}=N_1,
N_{2,0}=N_2,N_{2,\alpha^{*}}=N_3$. \item $i\leq \alpha^{*}
\Rightarrow N_{1,i} \preceq N_{2,i}$. \item $d_i \in
N_{1,i+1}-N_{1,i}$. 
\item $(N_{1,i},N_{1,i+1},d_i)\in K^{3,uq}$. \item
$tp(d_i,N_{2,i},N_{2,i+1})$ does not fork over $N_{1,i}$.
\end{enumerate}
In this case, $\langle N_{1,i}, d_i:i<\alpha^{*} \rangle ^\frown
\langle N_{1,\alpha^{*}} \rangle$ is said to be the first witness
for $NF^*(N_0,N_1,N_2,N_3)$, $d_i$ is said to be the $i$-th
element in the first witness for $NF^{*}$ and $\langle
N_{2,i}:i\leq \alpha^{*} \rangle$ is said to be the second witness
for $NF^*(N_0,N_1,N_2,N_3)$.
\end{definition}

\begin{definition} \label{5.3}\label{definition of NF}
$NF=NF_\lambda$ is the class of quadruples $(M_0,M_1,M_2,M_3)$
such that $M_0 \preceq M_1 \preceq M_3$, $M_0 \preceq M_2 \preceq
M_3$ and there are models $N_0,N_1,N_2,\allowbreak N_3$ such that:
$N_0=M_0, l<4\Rightarrow M_l\preceq N_l$
 and $NF^*(N_0,N_1,N_2,N_3)$.
\end{definition}

\begin{proposition}
The relations $NF^*,NF$ are closed under isomorphisms.
\end{proposition}

\begin{proof}
Trivial.
\end{proof}

\begin{proposition} \label{NF respects E}
Suppose for $x=a,b$ $(f_{x,1},f_{x,2},M_{x,3})$ is an amalgamation
of $M_1,M_2$ over $M_0$. If
$(f_{a,1},f_{a,2},M_{a,3})E_{M_0}(f_{b,1},f_{b,2},M_{b,3})$, then
$$NF(M_0,f_{a,1}[M_1],f_{a,2}[M_2],M_{a,3}) \Leftrightarrow
NF(M_0,f_{b,1}[M_1],f_{b,2}[M_2],M_{b,3})$$
\end{proposition}

\begin{proof}
Easy.
\end{proof}

\begin{proposition} \label{reduced}
\mbox{}
\begin{enumerate}
\item Every triple in $K^{3,uq}$ is reduced. \item If
$NF^*(N_0,N_1,N_2,N_3)$ then $N_1 \bigcap N_2=N_0$. \item If
$NF(N_0,N_1,N_2,N_3)$ then $N_1 \bigcap N_2=N_0$.
\end{enumerate}
\end{proposition}

\begin{proof}
(1) Suppose $(N_0,N_1,d) \preceq_{bs} (N_2,N_3,d), (N_0,N_1,d) \in
K^{3,uq}$. By Hypothesis \ref{5.0} and Proposition \ref{3.13}
(page \pageref{3.13}) there is a disjoint amalgamation of
$N_1,N_2$ over $N_0$, such that the type of $d$ does not fork over
$N_0$, so by the definition of uniqueness triple (definition 4.4),
$N_3$ is a disjoint amalgamation of $N_1,N_2$ over $N_0$.

(2) Let $x \in N_1 \bigcap N_2$. We will prove that $x \in N_0$.
Let $\langle N_{1,\alpha},d_\alpha:\alpha<\alpha^{*}
\rangle^\frown\langle N_{1,\alpha^{*}} \rangle, \langle
N_{2,\alpha}:\alpha \leq \alpha^{*} \rangle$ be witnesses for
$NF^{*}(N_0,N_1,N_2,N_3)$. Let $\alpha$ be the first ordinal such
that $x \in N_{1,\alpha}$. $\alpha$ is not a limit ordinal,
because a first witness for $NF^{*}$ is especially a continuous
sequence. we prove that $\alpha$ is not a successor ordinal, so we
conclude that $\alpha=0$. Suppose $\alpha=\beta+1$. By Definition
\ref{5.2}.e $(N_{1,\beta},N_{1,\beta+1},d_\beta) \in K^{3,uq}$. By
Definition \ref{5.2}.f $tp(d_\beta,N_{1,\beta},N_{1,\beta+1})$
does not fork over $N_{0,\beta}$. So by Proposition
\ref{reduced}.1 $N_{1,\beta+1} \bigcap N_{2,\beta}=N_{1,\beta}$.
But $x \in N_{1,\beta+1} \bigcap N_2 \subseteq N_{1,\beta+1}
\bigcap N_{2,\beta}$, so $x \in N_{1,\beta}$ in contradiction to
the assumption that $\alpha$ is the minimal ordinal with $x \in
N_{1,\alpha}$.

(3) By 2.
\end{proof}

\begin{theorem}[the existence theorem for $NF$]
\label{the existence theorem for NF} Suppose that for $n=1,2$ $N_0
\preceq N_n$ and $N_1 \bigcap N_2=N_0$.
\begin{enumerate}[(a)] \item  For some model $N_3 \in K_\lambda$ $NF(N_0,N_1,N_2,N_3)$. \item Moreover, if $N_1$ is decomposable over
$N_0$ by $K^{3,uq}$ then for some $N_3 \in K_\lambda$
$NF^{*}(N_0,N_1,N_2,N_3)$. \item Moreover, letting $a \in N_1-N_0$
we can choose $a$ as the first element in the first witness for
$NF^{*}$.
\end{enumerate}
\end{theorem}

\begin{proof}
\mbox{}
\begin{enumerate}[(a)]
\item By Theorem \ref{3.7}.3 (the extensions decomposition
theorem, page \pageref{3.7}), (and assumption \ref{5.0}), there is
a model $N^{*}_1$ with $N_1 \preceq N^{*}_1$ which is decomposable
over $N_0$, i.e. there is a sequence $\langle
N_{1,\alpha},d_\alpha:\alpha<\alpha^{*} \rangle^\frown \langle
N_{1,\alpha^{*}} \rangle$, such that: $N_0=N_{1,0},
(N_{1,\alpha},N_{1,\alpha+1},d_\alpha) \in K^{3,uq},N_1 \preceq
N_{1,\alpha^{*}}=N^{*}_1$. Therefore we can use item b. \item Let
$\langle N_{1,\alpha},d_\alpha:\alpha<\alpha^* \rangle ^\frown
\langle N_{1,\alpha^*} \rangle$ be an increasing continuous
sequence with $N_{1,0}=N_0$ and $N_{1,\alpha^*}=N_1$. By
Proposition \ref{3.3}.1
there is a sequence $\langle N_{2,\alpha}:\alpha \leq
\alpha^{*} \rangle$ which is a corresponding second witness for
$NF^{*}(N_0,N_{1,\alpha^{*}},N_2,N_{2,\alpha^{*}})$. \item By the
`more over' in Theorem \ref{3.7}.3 (the decomposing extensions
theorem, page \pageref{3.7}).
\end{enumerate}
\end{proof}

The following theorem is a private case of Theorem \ref{5.12},
i.e. the long transitivity theorem.
\begin{proposition}\label{private case of
transitivity} For $x=a,b$ let $\langle M_{x,\alpha}:\alpha \leq
\alpha^{*} \rangle$ be an increasing continuous sequence of
models. Suppose $\alpha<\alpha^{*} \Rightarrow
NF^{*}(M_{a,\alpha},M_{a,\alpha+1}, M_{b,\alpha},\allowbreak
M_{b,\alpha+1})$. Then $NF^*(M_{a,0},M_{a,\alpha^{*}},M_{b,0},
M_{b,\alpha^{*}})$.
\end{proposition}

\begin{proof}
Concatenate all the sequences together.
\end{proof}

\begin{proposition}[the monotonicity theorem] \label{5.4b}\label{monotonicity of NF}
\mbox{}
\begin{enumerate}
 \item If $NF^*(N_0,N_1,N_2,N_3)$ and $N_0 \preceq
M_2 \preceq N_2$, then $NF^*(N_0,N_1,M_2,\allowbreak N_3)$. \item
If $NF(M_0,M_1,M_2,M_3)$ then we can find $N_1,N_3$ such that
$NF^*(M_1,\allowbreak N_1,M_2,N_3)$ and $M_1 \preceq N_1 \preceq
N_3 \wedge M_3 \preceq N_3$. \item $NF^*(M_0,M_1,M_2,M_3) \wedge
M_3 \preceq M_3^* \Rightarrow NF(M_0,\allowbreak M_1,M_2,M_3^*)$.
\item The relation $NF$ satisfies monotonicity (in the sense of
Definition \ref{5.1}.b).
\end{enumerate}
\end{proposition}

\begin{proof}
\mbox{}\\
(1) Let $\langle N_{1,\alpha},d_\alpha:\alpha<\alpha ^* \rangle,\
\langle N_{2,\alpha}:\alpha<\alpha ^* \rangle$ be witnesses for
$NF^*(N_0,N_1,\allowbreak N_2,N_3)$. Then $\langle
N_{1,\alpha}:\alpha<\alpha ^* \rangle,\ \langle M_2 \rangle
^\frown \langle N_{2,\alpha}:0<\alpha<\alpha ^* \rangle$ are
witnesses for $NF^*(N_0,N_1,N_2,N_3)$ (notice that by Definition
\ref{2.1a}.3.b (monotonicity)
$tp(d_0,M_2,N_{2,1})$ does not fork over $N_0$).\\
(2) By the definition of $NF$ (Definition \ref{5.3}) and item 1.\\
(3)
\begin{displaymath}
\xymatrix{a \in M_1^* \ar[rr]^{f} && M_3^{**}\\
M_1 \ar[r]^{id} \ar[u]^{id} & M_3 \ar[r]^{id} & M_3^*
\ar[u]^{id}\\
M_0 \ar[r]^{id} \ar[u]^{id} & M_2 \ar[u]^{id} }
\end{displaymath}

Take $p \in S^{bs}(M_1)$, and take $M_1^*,a$ such that
$(M_1,M_1^*,a) \in p \bigcap K^{3,uq}$. By Definition
\ref{definition good frame}.1.3.f (on page \pageref{definition of
a good frame}) there is an amalgamation $(f,id_{M_3^*},M_3^{**})$
of $M_1^*,M_3^*$ over $M_1$ such that $tp(a,f[M_3^*],M_3^{**})$
does not fork over $M_1$. So $NF^*(M_1,f[M_1^*],M_3^*,\allowbreak
M_3^{**})$. Hence by item 1, $NF^*(M_1,f[M_1^*],M_3,M_3^{**})$.
Now by Proposition \ref{private case of transitivity}
$NF^*(\allowbreak M_0,M_1^*,M_2,M_3^{**})$. So the definition of
$NF$ (Definition \ref{5.3}),
$NF(M_0,M_1,\allowbreak M_2,M_3^*)$.\\
(4) Suppose $M_0^*=M_0,\ 0<n<3 \Rightarrow M_0^* \preceq M_n^*
\preceq M_3^*,\ M_n^* \preceq M_n,\ M_3^* \preceq M_3^{**},\ M_3
\preceq M_3^{**},\ NF(M_0,M_1,M_2,M_3)$.

\begin{displaymath}
\xymatrix{&& M_3^{**} \ar[r]^{f} & M_3^{***}\\
N_1 \ar[rrr]^{id} &&& N_3 \ar[u]^{id}\\
M_1 \ar[u]^{id} \ar[rr]^{id} && M_3 \ar[uu]^{id} \ar[ru]^{id}\\
M_1^* \ar[r]^{id} \ar[u]^{id} & M_3^* \ar[ruuu]^{id}\\
M_0 \ar[r]^{id} \ar[u]^{id} & M_2^* \ar[r]^{id} \ar[u]^{id} &
M_2=N_2 \ar[uu]^{id} }
\end{displaymath}

By item 2, for some $N_1,N_3$, $NF^*(M_0,N_1,M_2,N_3)$, $M_1
\preceq N_1 \preceq N_3$ and $M_3 \preceq N_3$. Take an
amalgamation $(f,id_{N_3},M_3^{***})$ of $M_3^{**}$ and $N_3$ over
$M_3$ (so over $M_1^* \bigcup M_2^*$). By item 3
$NF(M_0,N_1,M_2,M_3^{***})$. So by the definition of $NF$
(Definition \ref{5.3}), $NF(M_0,M_1^*,M_2^*,f[M_3^*])$. But the
relation $NF$ is closed under isomorphisms, so
$NF(M_0,M_1^*,M_2^*,M_3^*)$.
\end{proof}

\subsection{Weak Uniqueness}

We want to show that $NF$ satisfies weak uniqueness and long
transitivity. Proposition \ref{the transitivity of the weak
uniqueness} is a key point. To emphasize the exact hypotheses
involved in the proof, we extract from the axioms
$\bigotimes_{R}$, a smaller set $\bigotimes^-_{R}$.

\begin{definition} \label{5.1-}
Let $R \subseteq\ ^4(K_\lambda)$ be a relation. We say
$\bigotimes^-_{R}$ when:
\begin{enumerate}
\item If $R(M_0,M_1,M_2,M_3) \ then \ n\in \{1,2\}\Rightarrow
M_0\preceq M_n\preceq M_3$. \item Weak Uniqueness: Suppose for
$x=a,b$ $(f^x_1,f^x_2,N^x_3)$ is an amalgamation of $N_1$ and
$N_2$ over $N_0$ and $R(N_0,f^x_1[N_1],f^x_2[N_2],N^{x}_3)$. Then
$(f^a_1,f^a_2,N^a_3)E_{N_0}(f^b_1,f^b_2,N^b_3)$. \item If
$R(M_0,M_1,M_2,M_3)$ and $f:M_2 \to M_4$ is an embedding, then
there is an amalgamation $(g,id_{M_4},M_5)$ of $M_3,M_4$ over
$M_2$ such that $R(f[M_0],\allowbreak g[M_1],M_4,M_5)$.
\begin{displaymath}
\xymatrix{M_1 \ar[r]^{id} & M_3 \ar[r]^{g} & M_5\\
M_0 \ar[u]^{id} \ar[r]^{id} & M_2 \ar[u]^{id} \ar[r]^{f} & M_4
\ar[u]^{id}}
\end{displaymath}
\end{enumerate}
\end{definition}

\begin{definition}
$NF^{uq}:=\{(M_0,M_1,M_2,M_3)$:there is $a \in M_1-M_0$ such that
$(M_0,M_1,a) \in K^{3,uq}$ and $tp(a,M_2,M_3)$ does not fork over
$M_0$\}.
\end{definition}

\begin{proposition}\label{examples of bigotimes-}
\mbox{}
\begin{enumerate}
 \item $\bigotimes^-_{NF^{uq}}$.  \item For every relation $R$,
$\bigotimes_{R} \Rightarrow \bigotimes^-_{R}$.
\end{enumerate}
\end{proposition}

\begin{proof}
\mbox{}
\begin{enumerate}
\item By the definition of $K^{3,uq}$ (Definition \ref{definition
of uniqueness triples}), Definition \ref{definition good
frame}.3.f and Definition \ref{definition good frame}.1.d (to get
$M_5$). \item By axioms d,f in Definition \ref{5.1} and by
Proposition \ref{transitivity}.
\end{enumerate}
\end{proof}

\begin{proposition}[the transitivity of the weak uniqueness] \label{the transitivity
of the weak uniqueness} Suppose
\begin{enumerate} \item $\bigotimes^-_{R}$. \item $\alpha^* \leq \lambda^+$.
\item For every $\alpha<\alpha^*$
$N_{1,\alpha},N^a_{2,\alpha},N^b_{2,\alpha} \in K_\lambda$. \item
$\langle N_{1,\alpha}:\alpha \leq \alpha^{*} \rangle,\ \langle
N^a_{2,\alpha}:\alpha \leq \alpha^{*} \rangle,\ \langle
N^b_{2,\alpha}:\alpha \leq \alpha^{*} \rangle$ are increasing
continuous sequences. \item $N^a_{2,0}=N^b_{2,0}$. \item For every
$\alpha \leq \alpha^*$ $f^a_\alpha:N_{1,\alpha} \to
N^a_{2,\alpha}$ and $f^b_\alpha:N_{1,\alpha} \to N^b_{2,\alpha}$.
\item $(\alpha<\alpha^{*} \wedge x \in \{a,b\}) \Rightarrow
R(f^x_\alpha[N_{1,\alpha}],f^x_{\alpha+1}[N_{1,\alpha+1}],N^x_{2,\alpha},N^x_{2,\alpha+1})$.
\end{enumerate}
Then $(f^a_{\alpha^*},id_{N^a_{2,0}},N^a_{2,\alpha^*})E_{N_{1,0}}
(f^a_{\alpha^*},id_{N^a_{2,0}},N^b_{2,\alpha^*})$.
\end{proposition}

\begin{proof}
We choose $N_{2,\alpha},g_{a,\alpha},g_{b,\alpha}$ by induction on
$\alpha \leq \alpha^{*}$, such that for $x=a,b$ and $\alpha \leq
\alpha^*$ the following hold:
\begin{enumerate}[(i)]
\item $g_{x,\alpha}:N^x_{2,\alpha} \to N_{2,\alpha}$ is an
embedding such that $g_{a,\alpha} \circ f^a_\alpha=g_{b,\alpha}
\circ f_{b,\alpha}$. \item $N_{2,0}=N^x_{2,0}, g_{x,0}=identity$.
\item $\langle N_{2,\alpha}:\alpha \leq \alpha^{*} \rangle$ is an
increasing continuous sequence. \item $\langle g_{x,\alpha}:\alpha
\leq \alpha^{*} \rangle$ is an increasing continuous sequence.
\end{enumerate}

If we can construct this, then the following diagram commutes:

\begin{displaymath}
\xymatrix{& N^a_{2,\alpha^*} \ar[rr]^{g_{a,\alpha^*}} &&
N_{2,\alpha^*}\\
N_{1,\alpha} \ar[ru]^{f^a_{\alpha^*}} \ar[rr]^{f^b_{\alpha^*}} &&
N^b_{2,\alpha^*} \ar[ru]^{g_{b,\alpha^*}}\\
N_{1,0} \ar[u]^{id} \ar[r]^{id} & N_{2,0} \ar[uu]^{id}
\ar[ru]^{id} }
\end{displaymath}
[By clause (i) $g_{a,\alpha^*} \circ f^a_{\alpha^*}=g_{b,\alpha^*}
\circ f^b_{\alpha^*}$ and by clauses $(ii),(iv)$ $g_{x,\alpha^*}
\supseteq g_{x,0}=id_{N_{2,0}}$].

Therefore $(g_{a,\alpha^*},g_{b,\alpha^*},N_{2,\alpha^*})$
witnesses that $(f^a_{\alpha^*},id \allowbreak
_{N_{2,0}},N^a_{2,\alpha^*})E_{N_{1,0}}
(f^b_{\alpha^*},\allowbreak id_{N_{2,0}},\allowbreak
N^b_{2,\alpha^*})$.

Why can we construct this? For $\alpha=0$ only clause (ii) is
relevant. For $\alpha$ limit ordinal, take unions, and by the
smoothness, $g_{x,\alpha}$ is $\preceq$-embedding. What will we do
for $\alpha+1$? By clause 7 for $x=a,b$
$R(f^x_\alpha[N_{1,\alpha}],f^x_{\alpha+1}[N_{1,\alpha+1}],\allowbreak
N^x_{2,\alpha},N^x_{2,\alpha+1})$. By clause (i)
$g_{x,\alpha}[N^x_{2,\alpha}] \preceq N_{2,\alpha}$ and by clause
1 $\bigotimes^-_{R}$. So by Definition \ref{5.1-}.3 we can find
$g_x,N^x$ such that the following hold:
\begin{enumerate}
\item $g_x:N_{2,\alpha+1}^x \to N^x$ is an embedding. \item
$g_{x,\alpha} \subset g_x$. \item $R(g_x \circ
f^a_\alpha[N_{1,\alpha}],g_x \circ
f^a_{\alpha+1}[N_{1,\alpha+1}],N_{2,\alpha},N^x)$.
\end{enumerate}

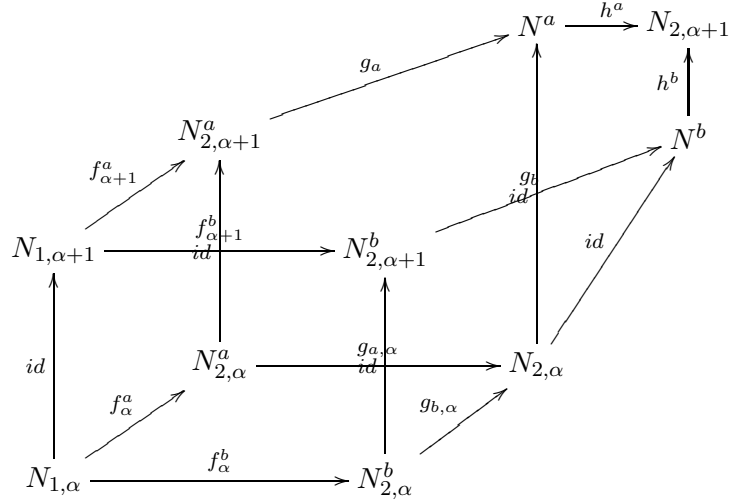
\begin{displaymath}
\xymatrix{&&& N^a \ar[r]^{h^a} & N_{2,\alpha+1}\\
&N^a_{2,\alpha+1} \ar[rru]^{g_a} &&& N^b \ar[u]^{h^b} \\
N_{1,\alpha+1} \ar[rr]^{f^b_{\alpha+1}} \ar[ru]^{f^a_{\alpha+1}} && N^b_{2,\alpha+1}  \ar[rru]^{g_b}\\
&N^a_{2,\alpha} \ar[rr]^{g_{a,\alpha}} \ar[uu]^{id} && N_{2,\alpha} \ar[ruu]^{id} \ar[uuu]^{id}\\
 N_{1,\alpha} \ar[rr]^{f^b_\alpha}  \ar[ru]^{f^a_\alpha} \ar[uu]^{id} && N^b_{2,\alpha} \ar[uu]^{id} \ar[ru]^{g_{b,\alpha}} }
\end{displaymath}

Hence by Definition \ref{5.1-}.2 $(g_a \restriction
f^a_{\alpha+1}[N_{1,\alpha+1}], \allowbreak id_{N_{2,\alpha}},
\allowbreak N^a)E_{f^a_\alpha[N_{1,\alpha}]}(g_b \restriction
f^b_{\alpha+1}[N_{1,\alpha+1}],id_{N_{2,\alpha}},N^b)$. So there
is a joint embedding $(h^a,h^b,N_{2,\alpha+1})$ of $N^a,N^b$ such
that for $x=a,b$ $id_{N_{2,\alpha}} \subseteq h^x$ and $h^a \circ
g_\alpha \circ f^a_{\alpha+1}=h^b \circ g_b \circ f^b_{\alpha+1}$.
Now we define $g_{x,\alpha+1}:=h^x \circ g_x$.
\end{proof}

The following proposition asserts that we have weak uniqueness
over first witness for $NF^*$.
\begin{proposition}
\label{5.5} If for $x=a,b$ $NF^*(N_0,\allowbreak N_1,N_2,N^x_3)$
and they have the same first witness, then there is a joint
embedding of $N^a_3,N^b_3$ over $N_1 \bigcup N_2$.
\end{proposition}

\begin{proof}
By Proposition \ref{examples of bigotimes-}.1,
$\bigotimes^-_{NF^{uq}}$. Hence it follows by Proposition \ref{the
transitivity of the weak uniqueness}.
\end{proof}

The following proposition is similar to weak uniqueness for
$NF^*$, but notice to the order of $N_1,N_2$ in the two
quadruples.
\begin{proposition}[the opposite uniqueness proposition]\label{opposite}
Suppose $NF^{*}(N_0,\allowbreak N_1,N_2,N^a_3)$ and
$NF^{*}(N_0,N_2,N_1,N^b_3)$. Then there is a joint embedding of
$N^a_3$ and $N^b_3$ over $N_1 \bigcup N_2$.
\end{proposition}

\begin{proof}
Let $\langle N^a_\alpha,d^a_\alpha:\alpha<\alpha^*
\rangle^\frown\langle N^a_{\alpha^*} \rangle$ be a first witness
for $NF^*(N_0,N_1,N_2,\allowbreak N^a_3)$ and let $\langle
N^b_\beta,d^b_\beta:\beta<\beta^* \rangle^\frown\langle
N^b_{\beta^*} \rangle$ be a first witness for
$NF^*(N_0,N_2,N_1,\allowbreak N^b_3)$. By Proposition \ref{3.4}
(page \pageref{3.4}), there is a rectangle
$\{M_{\alpha,\beta}:\alpha \leq \alpha^{*}, \beta \leq \beta^*\}$
such that:
\begin{enumerate}
\item $M_{\alpha,0}=N^a_\alpha$. \item $M_{0,\beta}=N^b_\beta$.
\item $tp(d^a_\alpha,M_{\alpha,\beta},M_{\alpha+1,\beta})$ does
not fork over $M_{\alpha,0}$. \item
$tp(d^b_\beta,M_{\alpha,\beta},M_{\alpha,\beta+1})$ does not fork
over $M_{0,\beta}$.
\end{enumerate}

\begin{displaymath}
\xymatrix{
&&N^a_3 \ar[r]^{f^a} &N^{a,*}_3 \ar[r]^{g^a} &N^*\\
N_1=N_{1,\alpha^*} \ar[rrr]^{id} \ar[rru]^{id} \ar[rrrrd]^{id} &&&M_{\alpha^*,\beta^*} \ar[r]^{id} \ar[u]^{id} &N^{b,*}_3 \ar[u]^{g^b} \\
d^a_\alpha \in N^a_{\alpha+1} \ar[r]^{id} \ar[u]^{id} &M_{\alpha+1,\beta} \ar[r]^{id} &M_{\alpha+1,\beta+1} \ar[ru]^{id} &&N^b_3 \ar[u]^{f^b}\\
N^a_{\alpha}=M_{\alpha,0} \ar[r]^{id} \ar[u]^{id} &M_{\alpha,\beta} \ar[r]^{id} \ar[u]^{id} &M_{\alpha,\beta+1} \ar[u]^{id}\\
N^a_1=M_{1,0} \ar[r]^{id} \ar[u]^{id} &M_{1,\beta} \ar[r]^{id} \ar[u]^{id} &M_{1,\beta+1} \ar[u]^{id}\\
N_0=M_{0,0} \ar[r]^{id} \ar[u]^{id} &N^b_\beta=M_{0,\beta}
\ar[r]^{id} \ar[u]^{id} &N^b_{\beta+1} \ar[r]^{id} \ar[u]^{id}
&N_2=M_{0,\beta^*} \ar[uuur]^{id} \ar[uuuu]^{id} \ar[uuuuul]^{id}}
\end{displaymath}

By clauses 1,3 $\langle d^a_\alpha,N^a_\alpha:\alpha<\alpha^a
\rangle$ is a first witness for $NF^{*}(N_0,N_1,N_2,\allowbreak
M_{\alpha^*,\beta^*})$. But by definition this is also a first
witness for $NF^{*}(N_0,N_1,\allowbreak N_2,N^a_3)$. So by
Proposition \ref{5.5}, there is a joint embedding
$(id_{M_{\alpha^*,\beta^*}},f^a,N^{a,*}_3)$ of
$M_{\alpha^*,\beta^*},N^a_3$ over $N_1 \bigcup N_2$. Similarly by
clauses 2,4 there is a joint embedding
$(id_{M_{\alpha^*,\beta^*}},f^b,N^{b,*}_3)$ of
$M_{\alpha^*,\beta^*},N^b_3$ over $N_1 \bigcup N_2$. Since
$(K_\lambda,\preceq \restriction K_\lambda)$ has amalgamation,
there is an amalgamation $(g^a,g^b,N_3)$ of $N^{a,*}_3,N^{b,*}_3$
over $M_{\alpha^*,\beta^*}$. $N_3$ is an amalgam of $N^a_3,N^b_3$
over $N_1 \bigcup N_2$.
\end{proof}

\begin{theorem}[the weak uniqueness theorem]
\label{the weak uniqueness theorem of NF} Suppose for $x=a,b$
$NF(M_0 \allowbreak ,M_1,\allowbreak M_2,M^x)$. Then there is a
joint embedding of $M^a,M^b$ over $M_1 \bigcup M_2$.
\end{theorem}

\begin{proof}
First note that since $M_1 \bigcap M_2=M_0$, the conclusion of the
theorem is equivalent to
$(id_{M_1},id_{M_2},M^a)E_{M_0}(id_{M_1},id_{M_2},M^b)$.

\case{Case a:} $NF^*(M_0,M_1,\allowbreak M_2,M^x)$ and $M_2$ is
decomposable over $M_0$. In this case, by Theorem \ref{the
existence theorem for NF}.b (the existence theorem for $NF$) there
is $M^c$ such that $NF^*(M_0,M_2,M_1,\allowbreak M^c)$. By
Proposition \ref{opposite} for $x=a,b$ $id_{M_1},id_{M_2},M^x)
E_{M_0} \allowbreak (id_{M_1},id_{M_2},M^c)$. But the relation
$E_{M_0}$ is an equivalence relation, so it is transitive.

\case{The general case:} Since $NF(M_0,M_1,M_2,M^a,)$ by
Proposition \ref{5.4b}.5 there are $N^a_1,N^{a,-}$ such that
$NF^*(M_0,N^a_1,M_2,N^{a,-})$ and $M_1 \preceq N^a_1 \preceq
N^{a,-} \wedge M^a \preceq N^{a,-}$. Similarly there are
$N^b_1,N^{b,-}$ such that $NF^*(M_0,N^b_1,M_2,\allowbreak N^b_3)$
and $M_1 \preceq N^b_1 \preceq N^{b,-} \wedge M^b \preceq
N^{b,-}$. By Theorem \ref{the extensions decomposition theorem}
(the extensions decomposition theorem) there is a model $M_2^+
\succeq M_2$ which is decomposable over $M_0$. Without loss of
generality for $x=a,b$ $M_2^+ \bigcap N^{x,-}=M_2$. So by Theorem
\ref{the extensions decomposition theorem}.3 (the extensions
decomposition theorem) there is $N^x \succeq N^{x,-}$ such that
$NF^*(M_0,N^x_1,\allowbreak M_2^+,N^x)$.

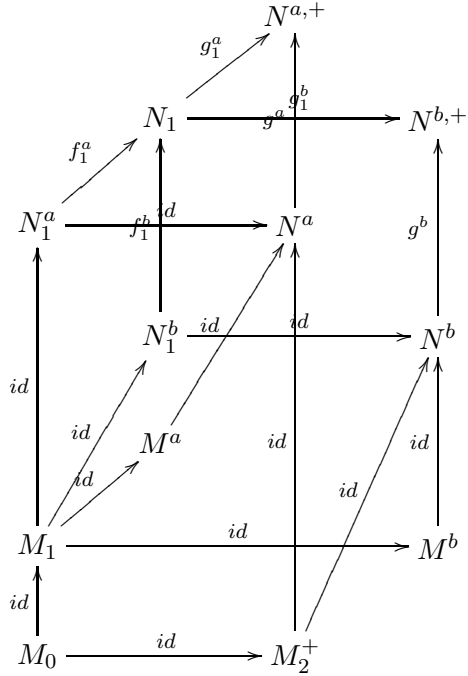
\begin{displaymath}
\xymatrix{&& N^{a,+}\\
& N_1 \ar[ru]^{g^a_1} \ar[rr]^{g^b_1} && N^{b,+}\\
N^a_1 \ar[ru]^{f^a_1} \ar[rr]^{id} && N^a \ar[uu]^{g^a}\\
& N^b_1 \ar[uu]^{f^b_1} \ar[rr]^{id} && N^b \ar[uu]^{g^b}\\
& M^a \ar[uur]^{id}\\
M_1 \ar[uuu]^{id} \ar[ruu]^{id} \ar[ru]^{id} \ar[rrr]^{id} &&& M^b \ar[uu]^{id}\\
M_0 \ar[u]^{id} \ar[rr]^{id} && M_2^+ \ar[uuuu]^{id}
\ar[ruuu]^{id}}
\end{displaymath}

By Proposition \ref{existence of decomposition over two models}
there is an amalgamation $(f^a_1,f^b_1,N_1)$ of $N^a_1,N^b_1$ over
$M_1$ such that $N_1$ is decomposable over $N^a_1$ and over
$N^b_1$. Hence for $x=a,b$ there is an amalgamation
$(g^x_1,g^x,N^{x,+})$ of $N_1,N^x$ over $N^x_1$ such that
$NF^*(g^x_1 \circ f^x_1[N^x_1],g^x_1[N_1],g^x[N^x],N^{x,+})$. So
for $x=a,b$ by Proposition \ref{5.4b}.8 (a private case of
transitivity), since $NF^*(M_0,N^x_1,M_2^+,N^x)$ and
$NF^*(N^x_1,N^x,N_1,N^{x,+})$ it follows that $NF^*(M_0
\allowbreak ,N_1,M_2^+,N^{x,+})$. So by case a $(g^a_1,g^a
\restriction M_2^+,N^{a,+})E_{M_0}(g^b_1,g^b \restriction
M_2^+,N^{b,+}$. Therefore $(g^a_1 \restriction M_1,g^a
\restriction M_2,N^{a,+})E_{M_0}(g^b_1 \restriction M_1,g^b
\restriction M_2,N^{b,+})$.

\end{proof}

\begin{proposition} \label{bigotimes-NF}
$\bigotimes^-_{NF}$.
\end{proposition}

\begin{proof}
We have to check clauses 1,2,3 of Definition \ref{5.1-} (on page
\pageref{5.1-}).

1. Trivial.

2. By Theorem \ref{the weak uniqueness theorem of NF}.

3. Suppose $NF(M_0,M_1,M_2,M_3)$ and $f:M_2 \to M_4$ is an
embedding. We have to find a model $M_5$ and an embedding $g:M_3
\to M_5$ over $M_2$ such that $NF(f[M_0],\allowbreak
g[M_1],M_4,M_5)$. By Theorem \ref{monotonicity of NF}.2 we can
find $N_1,N_3$ such that $NF^*(M_1,\allowbreak N_1,M_2,N_3)$ and
$M_1 \preceq N_1 \preceq N_3 \wedge M_3 \preceq N_3$. By Theorem
\ref{the existence theorem for NF}.b (the existence theorem for
$NF$, on page \pageref{the existence theorem for NF}) we can find
a model $M_5$ with $M_4 \preceq M_5$ and an embedding $h:N_3 \to
M_5$ such that $NF^*(M_0,M_4,N_1,M_5)$. Hence
$NF(M_0,M_1,M_4,M_5)$. Now we define $g:=h \restriction M_3$.
\end{proof}

\begin{theorem}[the symmetry theorem] \label{5.9}
$NF(N_0,N_1,N_2,N_3) \Leftrightarrow NF(N_0,\allowbreak
N_2,N_1,N_3)$.
\end{theorem}

\begin{proof}
By the monotonicity of NF, i.e. Propositon \ref{5.4b}.3, It is
sufficient to prove $NF^*(N_0,N_1,N_2,N_3) \Rightarrow
NF(N_0,N_2,N_1,N_3)$. Suppose $NF^*(N_0,N_1,\allowbreak N_2,N_3)$.
By Theorem \ref{the extensions decomposition theorem} (the
extensions decomposition theorem) there is $N_2^+ \succeq N_2$
which is decomposable over $N_0$. By Theorem \ref{the existence
theorem for NF}.b there is an amalgamation $(id_{N_1},f,N_3^+)$ of
$N_1,N_2^+$ over $N_2$ such that $NF^*(N_0,N_1,\allowbreak
f[N_2^+],N_3^+)$. So $N_1 \bigcap f[N_2^+]=N_0$. Hence by Theorem
\ref{the existence theorem for NF}.b, there is a model $N^*$ such
that $NF^*(N_0,f[N_2^+],\allowbreak N_1,N^*)$. By Proposition
\ref{opposite} (the opposite uniqueness proposition) there is a
joint embedding $id_{N_3^+},g,N^{**}$ of $N_3^+$ and $N^*$ over
$N_1 \bigcup f[N_2^+]$. Since $NF^*$ is closed under isomorphisms,
$NF^*(N_0,f[N_2^+],\allowbreak N_1,g[N^{*}])$. Now we have to use
the monotonicity of $NF$ twice. Since $N_0 \preceq N_2 \preceq
f[N_2^*]$ it follows that $NF^*(N_0,N_2,\allowbreak
N_1,g[N^{*}])$. Since $N_3 \preceq N_3^* \preceq N^{**} \succeq
g[N^*]$, it follows that $NF(N_0,N_2,N_1,N_3)$.
\end{proof}

\begin{theorem} \label{5.10}
$NF$ respects $\frak{s}$ (see Definition \ref{respect})
\end{theorem}

\begin{proof}
Suppose $NF(M_0,M_1,M_2,M_3),\ tp(a,M_0,M_1) \in S^{bs}(M_0)$. We
have to prove that $tp(a,M_2,M_3)$ does not fork over $M_0$.
Without loss of generality $NF^*(M_0,M_1,M_2,M_3)$ [Why? see the
Definition \ref{2.1a}.3.b (monotonicity)]. By the definition of
$NF^*,\ M_1$ is decomposable over $M_0$. By of NF, (Theorem
\ref{the existence theorem for NF}.c (the existence theorem for
$NF$), there is $M_3^*$ such that $NF^*(M_0,M_1,M_2,M_3^*)$ and
the first element in the first witness is a.

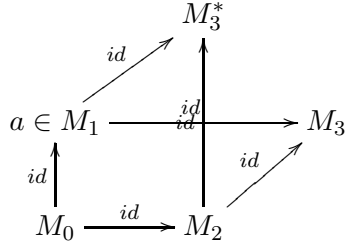
\begin{displaymath}
\xymatrix{&M_3^*\\ a \in M_1 \ar[rr]^{id} \ar[ru]^{id}&&M_3\\
M_0 \ar[r]^{id} \ar[u]^{id} & M_2 \ar[uu]^{id} \ar[ru]^{id} }
\end{displaymath}
By the definition of a first witness, $tp(a,M_2,M_3^*)$ does not
fork over $M_0$. By the weak uniqueness theorem (Theorem \ref{the
weak uniqueness theorem of NF}) there are $f,M_3^{**}$ such that
$M_3 \preceq M_3^{**},$ and $f:M_3^* \to M_3^{**}$ is an embedding
over $M_1 \bigcup M_2$. So
$tp(a,M_2,M_3)=tp(a,M_2,f[M_3^*])=tp(a,M_2,M_3^*)$ does not fork
over $M_0$.
\end{proof}
\subsection{Long transitivity}

\begin{proposition} \label{5.11}
Let $\langle M_\epsilon:\epsilon \leq \alpha^* \rangle$ be a
$\prec$-increasing continuous sequence of models in $K_\lambda$.
\begin{enumerate}[(a)]
\item There is an $\prec$-increasing continuous sequence of models
in $K_lambda$ $\langle N_\epsilon:\epsilon \leq \alpha^* \rangle$
such that: $N_0=M_0,\ M_\epsilon \preceq N_\epsilon$,
$NF(M_\epsilon,M_{\epsilon+1},N_\epsilon,N_{\epsilon+1})$ and
$N_{\epsilon+1}$ is decomposable over $N_\epsilon$ and over
$M_{\epsilon+1}$. \item Suppose $M^* \in K_\lambda$, $M^* \succ
M_0$ and $M^* \bigcap M_{\alpha^*}=M_0$. Then there is an
$\prec$-increasing continuous sequence of models in $K_\lambda$
$\langle N_\epsilon:\epsilon \leq \alpha^* \rangle$ such that:
$M^* \preceq N_0,\ M_\epsilon \preceq N_\epsilon$,
$NF(M_\epsilon,M_{\epsilon+1},N_\epsilon,N_{\epsilon+1})$, $N_0$
is decomposable over $M_0$ and $N_{\epsilon+1}$ is decomposable
over $N_\epsilon$ and over $M_{\epsilon+1}$.

\end{enumerate}
\end{proposition}

\begin{proof}
\mbox{} (a) We choose a pair $(N_\epsilon,f_\epsilon)$ by
induction on $\epsilon \leq \alpha$ such that:
\begin{enumerate}
\item $\langle N_\epsilon:\epsilon \leq \alpha \rangle$ is an
increasing continuous sequence of models in $K_\lambda$. \item
$f_\epsilon:M_\epsilon \to N_\epsilon$ is an embedding. \item
$f_0=id_{M_0}$. \item The sequence $\langle f_\epsilon:\epsilon
\leq \alpha \rangle$ is increasing and continuous. \item For
$\epsilon<\alpha^*$,
$NF(f_\epsilon[M_\epsilon],N_\epsilon,f_{\epsilon+1}[M_{\epsilon+1}],N_{\epsilon+1})$.
\item For $\epsilon<\alpha^*$, $N_{\epsilon+1}$ is decomposable
over $N_\epsilon$ and over $f_{\epsilon+1}[M_{\epsilon+1}]$.
\end{enumerate}
Why can we carry out this construction? For $\epsilon=0$ or limit
there is no problem. Suppose we chose $(N_\epsilon,f_\epsilon)$,
how will we choose $(N_{\epsilon+1},f_{\epsilon+1})$? By Theorem
\ref{the existence theorem for NF}.a we can find
$N_{\epsilon+1}^-$ and $f_{\epsilon+1}$ such that
$NF(f_\epsilon[M_\epsilon],N_\epsilon,f_{\epsilon+1}[M_{\epsilon+1}],\allowbreak
N_{\epsilon+1}^-)$. Now by Proposition \ref{existence of
decomposition over two models} we can find $N_{\epsilon+1}$ such
that $N_{\epsilon+1}^- \preceq N_{\epsilon+1}$ and
$N_{\epsilon+1}$ is decomposable over $N_\epsilon$ and over
$f_{\epsilon+1}[M_{\epsilon+1}]$. Therefore we can carry out this
construction.

Now, as in the proof of Proposition \ref{3.3}, without loss of
generality $f_\epsilon=id_{M_\epsilon}$ for every $\epsilon \leq
\alpha^*$ (because we can extend $f_{\alpha^*}^{-1}$ to a
bijection $g$ of $N_{\alpha^*}$ and take the sequence $\langle
g[N_\epsilon]:\epsilon \leq \alpha^* \rangle$).

(b) It demands a tiny change in the proof:
In the construction $M^* \preceq N_0$ and it is decomposable
over $M_0$.

\end{proof}

\begin{theorem} [the long transitivity theorem] \label{5.12}
For $x=a,b$ let $\langle  M_{x,\epsilon}:\epsilon \leq \alpha^*
\rangle$ be an $\prec$-increasing continuous sequence of models in
$K_\lambda$. Suppose $\epsilon<\alpha^* \Rightarrow
NF(M_{a,\epsilon},\allowbreak
M_{a,\epsilon+1},M_{b,\epsilon},M_{b,\epsilon+1})$. Then
$NF(M_{a,0},M_{a,\alpha^*},M_{b,0},M_{b,\alpha^*})$.
\end{theorem}

Similarly to the proof of Proposition \ref{transitivity} (the
transitivity proposition), we use the existence and weak
uniqueness theorems to prove the long transitivity. But here the
proof is more complicated, and it is divided to four cases, each
one is based on its previous and generalizes it.

\begin{proof}
\case{Case a:} $\epsilon<\alpha^* \Rightarrow
{NF^*}(M_{a,\epsilon},M_{a,\epsilon+1},M_{b,\epsilon},M_{b,\epsilon+1})$.
Concatenate all the sequences together.

In the other cases we are going to use the following claim:
\begin{claim}\label{for long transitivity}
It is enough to find $(N_{b,\epsilon},f_\epsilon)$ for $\epsilon
\leq \alpha^*$ such that:
\begin{enumerate}
\item $M_{b,0} \preceq N_{b,0}$. \item $\langle
N_{b,\epsilon}:\epsilon \leq \alpha^* \rangle$ is an increasing
continuous sequence of models in $K_\lambda$. \item $f_\epsilon$
is an embedding of $M_{a,\epsilon}$ to $N_{b,\epsilon}$. \item
$f_0=id_{M_{a,0}}$. \item $\langle f_\epsilon:\epsilon \leq
\alpha^* \rangle$ is an increasing continuous sequence. \item For
$\epsilon<\alpha^*$,
$NF(f_\epsilon[M_{a,\epsilon}],f_{\epsilon+1}[M_{a,\epsilon+1}],N_{b,\epsilon},N_{b,\epsilon+1})$.
\item
$NF(M_{a,0},f_{\alpha^*}[M_{a,\alpha^*}],N_{b,0},N_{b,\alpha^*})$.
\end{enumerate}
\end{claim}

\begin{proof}
Suppose we found $(N_{b,\epsilon},f_\epsilon)$ for $\epsilon \leq
\alpha^*$ such that clauses 1-7 are satisfied. By Proposition
\ref{bigotimes-NF}, $\bigotimes^-_{NF}$. Therefore by Proposition
\ref{the transitivity of the weak uniqueness} (the transitivity of
the uniqueness)
$(id_{M_{a,\alpha^*}},id_{M_{b,0}},M_{b,\alpha^*})E_{M_{a,0}}
\allowbreak (f^a_{\alpha^*},id_{M_{b,0}},N_{b,\alpha^*})$
[Substitute $\langle M_{a,\epsilon}:\epsilon \leq \alpha^*
\rangle,\ \langle M_{b,\epsilon}:\epsilon \leq \alpha^* \rangle,\
\langle N_{b,\epsilon}:\epsilon \leq \alpha^* \rangle,\ \langle
id_{M_{a,\epsilon}}:\epsilon \leq \alpha^* \rangle,\ \langle
f_\epsilon:\epsilon \leq \alpha^* \rangle$ in place of $\langle
N_{1,\alpha}:\alpha \leq \alpha^* \rangle,\ \langle
N^a_{2,\alpha}:\alpha \leq \alpha^* \rangle,\ \langle
N^b_{2,\alpha}:\alpha \leq \alpha^* \rangle,\ \langle
f^a_\alpha:\alpha \leq \alpha^* \rangle,\ \langle
f^b_\alpha:\alpha \leq \alpha^* \rangle$] . By clause 7
$NF(M_{a,0},M_{a,\alpha^*},N_{b,0},N_{b,\alpha^*})$. So by
Proposition \ref{NF respects E}
$NF(M_{a,0},M_{a,\alpha^*},M_{b,0},M_{b,\alpha^*})$.
\end{proof}

\case{Case b:} For every $\epsilon$, $M_{a,\epsilon+1}$ is
\emph{decomposable} over $M_{a,\epsilon}$. In this case we choose
$(N_{b,\epsilon},f_\epsilon)$ such that clauses 1-6 of Claim
\ref{for long transitivity} are satisfied: For $\epsilon=0$ we
define $N_{b,0}:=M_{b,0}$. In successor step we use Theorem
\ref{the existence theorem for NF}.a. For $\epsilon$ limit we
define $N_{b,\epsilon}:=\bigcup \{N_{b,\zeta}:\zeta<\epsilon\},\
f_\epsilon:=\bigcup \{f_\zeta:\zeta<\epsilon\}$. Now clause 7 is
satisfied by case a of the proof.

\case{Case c:} \emph{$\alpha^* \leq \omega$}. In this case we
apply Claim \ref{for long transitivity} with
$f_\epsilon=id_{M_{a,\epsilon}}$.
\begin{displaymath}
\xymatrix{ N_{b,0} \ar[r]^{id} & N_{b,1} \ar[r]^{id} & N_{b,2}
\ar[r]^{id} &
N_{b,\epsilon} \ar[r]^{id} & N_{b,\epsilon+1} \ar[r]^{id} & N_{b,\alpha^*}\\
M_{b,0} \ar[u]^{id}\\
M_{a,0} \ar[r]^{id} \ar[u]^{id} & N_{a,1} \ar[r]^{id} \ar[uu]^{id}
& N_{a,2} \ar[r]^{id} \ar[uu]^{id} &
N_{a,\epsilon} \ar[r]^{id}  \ar[uu]^{id} & N_{a,\epsilon+1} \ar[r]^{id}  \ar[uu]^{id} & N_{a,\alpha^*}  \ar[uu]^{id}\\
M_{a,0} \ar[r]^{id} \ar[u]^{id} & M_{a,1}  \ar[r]^{id} \ar[u]^{id}
& M_{a,2} \ar[r]^{id} \ar[u]^{id} & M_{a,\epsilon}  \ar[r]^{id}
\ar[u]^{id} & M_{a,\epsilon+1} \ar[r]^{id} \ar[u]^{id} &
M_{a,\alpha^*} \ar[u]^{id} }
\end{displaymath}

By Proposition \ref{5.11}.a, there is an increasing continuous
sequence $\langle N_{a,\epsilon}:\epsilon \leq \alpha^* \rangle$
such that: $N_{a,0}=M_{a,0},\ M_{a,\epsilon} \preceq
N_{a,\epsilon}$, $N_{a,\epsilon+1}$ is decomposable over
$N_{a,\epsilon}$ and over $M_{a,\epsilon+1}$ and
$\epsilon<\alpha^* \Rightarrow
NF(M_{a,\epsilon},M_{a,\epsilon+1},N_{a,\epsilon},N_{a,\epsilon+1})$.
Since $\alpha^* \leq \omega$, by Proposition \ref{5.11}.b, there
is an increasing continuous sequence $\langle
N_{b,\epsilon}:\epsilon \leq \alpha^* \rangle$ such that $N_{b,0}
\succ M_{b,0}$, for $\epsilon \leq \alpha^*$, $N_{b,\epsilon}$ is
decomposable over $N_{a,\epsilon}$ and
$NF^*(N_{a,\epsilon},N_{a,\epsilon+1},N_{b,\epsilon},N_{b,\epsilon+1})$.

Now it is enough to prove that $\langle
(N_{b,\epsilon},id_{M_{a,\epsilon}}):\epsilon \leq \alpha^*
\rangle$ satisfies clauses 1-7 of Claim \ref{for long
transitivity}. Clauses 1-5 are satisfied trivially. We check
clauses 6,7.

6. First assume $\epsilon>0$. As
$NF(M_{a,\epsilon},M_{a,\epsilon+1},N_{a,\epsilon},N_{a,\epsilon+1})$,
$NF(N_{a,\epsilon},\allowbreak
N_{a,\epsilon+1},N_{b,\epsilon},N_{b,\epsilon+1})$,
$N_{a,\epsilon}$ is decomposable over $M_{a,\epsilon}$ and
$N_{b,\epsilon}$ is decomposable over $N_{a,\epsilon}$, by case b
(for $\alpha^*=2$),
$NF(M_{a,\epsilon},M_{a,\epsilon+1},N_{b,\epsilon},N_{b,\epsilon+1})$.
Second assume $\epsilon=0$. As
$NF(N_{a,0},N_{a,1},N_{b,0},N_{b,1})$, $N_{a,0}=M_{a,0}$ and
$M_{a,1} \preceq N_{a,1}$, by the monotonicity of $NF$,
$NF(M_{a,0},M_{a,1},N_{b,0},N_{b,1})$

7. By case b, we have
$NF(N_{a,0},N_{a,\alpha^*},N_{b,0},N_{b,\alpha^*})$. By the
smoothness $M_{a,\alpha^*} \preceq N_{a,\alpha^*}$. So by the
monotonicity of $NF$,
$NF(M_{a,0},M_{a,\alpha^*},N_{b,0},N_{b,\alpha^*})$.

\case{The general case:} By the proof of case c. We have only one
problem: For $\epsilon$ limit it is not clear why does
$NF(M_{a,\epsilon},M_{a,\epsilon+1},N_{b,\epsilon},N_{b,\epsilon+1})$,
where we know
$NF(M_{a,\epsilon},M_{a,\epsilon+1},N_{a,\epsilon},N_{a,\epsilon+1})
\wedge
NF(N_{a,\epsilon},N_{a,\epsilon+1},N_{b,\epsilon},N_{b,\epsilon+1})$.
Here we cannot use case b, because we do not know if
$N_{b,\epsilon}$ is decomposable over $N_{a,\epsilon}$ and
$N_{a,\epsilon}$ is decomposable over $M_{a,\epsilon}$. But we can
use case c with $\alpha^*=2$.
\end{proof}

\begin{theorem} \label{5.13}
$NF=NF_\lambda$ is the unique relation which satisfies
$\bigotimes_{NF}$ and respects $\frak{s}$.
\end{theorem}

\begin{proof}
$NF$ satisfies $\bigotimes_{NF}$: Clause a is clear. Clause b (the
monotonicity) by Theorem \ref{5.4b}.4. Clause c (the existence) by
Theorem \ref{the existence theorem for NF}.a. Clause d (weak
uniqueness) by Theorem \ref{the weak uniqueness theorem of NF}.
Clause e (symmetry) by Theorem \ref{5.9}. Clause f (long
transitivity) by Theorem \ref{5.12}. By Theorem \ref{5.10} $NF$
respects $\frak{s}$.

Suppose the relation $R$ satisfies $\bigotimes_{R}$ and respects
$\frak{s}$. First we prove $NF(M_0,M_1,M_2,M_3) \Rightarrow
R(M_0,M_1,M_2,M_3)$.

\case{case a:} Take $a \in M_1-M_0$ with $(M_0,M_1,a) \in
K^{3,uq}$. As $NF$ respects $\frak{s}$, $tp(a,M_2,M_3)$ does not
fork over $M_0$. So as $R$ respects $\frak{s}$, by the definition
of unique triples (see Definition \ref{4.1} on page
\pageref{4.1}), $R(M_0,M_1,M_2,M_3)$.

\case{case b:} $NF^*(M_0,M_1,M_2,M_3)$. As $R$ satisfies long
transitivity, and by case a, $R(M_0,M_1,M_2,M_3)$.

\case{The general case:} Since $R$ satisfies monotonicity, by case
b, $R(M_0,M_1,M_2,\allowbreak M_3)$. So we have proved that the
relation $NF$ is included in the relation $R$.

\case{conversely}: Suppose $R(M_0,M_1,M_2,M_3)$. We have to prove
that $NF(M_0,M_1,M_2,M_3)$. As $\bigotimes_{R}$, $R$ satisfies
disjointness. So $M_1 \bigcap M_2=M_0$. By $\bigotimes_{NF}$, for
some model $M_4$ $NF(M_0,M_1,M_2,M_4)$. But by the first direction
of the proof $NF(M_0,M_1,\allowbreak M_2,M_4) \Rightarrow
R(M_0,M_1,M_2,M_4)$, so $R(M_0,M_1,M_2,M_4)$. As $\bigotimes_{R}$,
$R$ satisfies weak uniqueness, $R(M_0,M_1\allowbreak ,M_2,M_3)$
and $R(M_0,M_1,M_2,M_4)$, it follows that
$(id_{M_1},id_{M_2},\allowbreak M_3)E_{M_0}(id_{M_1}\allowbreak
,id_{M_2},M_4)$. Therefore by Proposition
 \ref{NF respects E} $NF(M_0,M_1,M_2,M_4)$ implies $NF(\allowbreak M_0,M_1,M_2,M_3)$, so $NF(M_0,M_1,M_2,M_3)$ as required.

\end{proof}

\section{A relation on $K_{\lambda^+}$ that is based on the relation $NF$}
Remember that we want to derive from $\frak{s}$ a good
$\lambda^+$-frame. So first we have to define an a.e.c. in
$\lambda^+$ with amalgamation.
Definition \ref{6.1} presents the
relation on models of this a.e.c. in $\lambda^+$.

\begin{hypothesis}
$\frak{s}$ is a weakly successful semi-good $\lambda$-frame with
conjugation.
\end{hypothesis}

\begin{definition} \label{5.14} Define a 4-place relation $\widehat{NF}$ on $K$ by\\ $\widehat{NF}(N_0,N_1,M_0,\allowbreak M_1)$ iff the following hold:
\begin{enumerate}
\item $n<2 \Rightarrow N_n \in K_\lambda,\ M_n \in K_{\lambda^+}$.
\item There is a pair of increasing continuous sequences $\langle
N_{0,\alpha}:\alpha<\lambda^+ \rangle,\ \langle
N_{1,\alpha}:\alpha<\lambda^+ \rangle$ such that for every
$\alpha,\ NF(N_{0,\alpha},N_{1,\alpha},N_{0,\alpha+1},\allowbreak
N_{1,\alpha+1})$ and for $n<2$, $N_{0,n}=N_n,\
M_n=\bigcup\{N_{n,\alpha}:\alpha<\lambda^+\}$.
\end{enumerate}
\end{definition}

\begin{theorem} [the $\widehat{NF}$-properties]
\label{5.15}\label{the widehat{NF}-properties} \mbox{}
\begin{enumerate}[(a)] \item Disjointness: If
$\widehat{NF}(N_0,N_1,M_0,M_1)$ then $N_1 \bigcap M_0=N_0$. \item
Monotonicity: Suppose $\widehat{NF}(N_0,N_1,M_0,M_1),\ N_0 \preceq
N^{*}_1 \preceq N_1,\ N_1^* \bigcup M_0 \allowbreak \subseteq
M^*_1 \preceq M_1$ and $M_1^* \in K_{\lambda^+}$. Then
$\widehat{NF}(N_0,N^{*}_1,M_0,M^*_1)$. \item Existence: Suppose
$n<2 \Rightarrow N_n \in K_\lambda,\ M_0 \in K_{\lambda^+},\ N_0
\preceq N_1,\ N_0 \preceq M_0,\ N_1 \bigcap M_0=N_0$. Then there
is a model $M_1$ such that $\widehat{NF}(N_0,N_1,\allowbreak
M_0,M_1)$. \item Weak Uniqueness: If $n<2 \Rightarrow
\widehat{NF}(N_0,N_1,M_0,M_{1,n})$, then there are $M,f_0,f_1$
such that $f_n$ is an embedding of $M_{1,n}$ into $M$ over $N_1
\bigcup M_0$. \item Respecting the frame: Suppose
$\widehat{NF}(N_0,N_1,M_0,M_1), tp(a,N_0,M_0) \in \allowbreak S^
{bs} \allowbreak (N_0)$. Then $tp(a,N_1,M_1)$ does not fork over
$N_0$.
\end{enumerate}
\end{theorem}
\begin{proof}
\mbox{} (a) Disjointness: Let $\langle
N_{0,\epsilon}:\epsilon<\lambda^+ \rangle,\ \langle
N_{1,\epsilon}:\epsilon<\lambda^+ \rangle$ be witnesses for
$\widehat{NF}(N_0,N_1,M_0,M_1)$. Especially $\epsilon<\lambda^+
\Rightarrow NF(N_{0,\epsilon},N_{1,\epsilon},\allowbreak
N_{0,\epsilon+1},N_{1,\epsilon+1})$. So by Theorem \ref{reduced}.3
$\epsilon<\lambda^+ \Rightarrow N_{1,\epsilon} \bigcap
N_{0,\epsilon+1}=N_{0,\epsilon}$. So by the proof of Proposition
\ref{reduced}.2 $N_1 \bigcap M_0=N_0$. Let $x \in N_1 \bigcap
M_0$. So there is $\epsilon<\lambda^+$ such that $x \in
N_{0,\epsilon}$. Denote $\epsilon:=\Min\{\epsilon<\lambda^+:x \in
N_{0,\epsilon}\}$. $\epsilon$ cannot be a limit ordinal as the
sequence $\langle N_{0,\epsilon}:\epsilon<\lambda^+ \rangle$ is
continuous. If $\epsilon=\zeta+1$ then $x \in N_{0,\zeta+1}
\bigcap N_1 \subseteq N_{0,\zeta+1} \bigcap
N_{1,\zeta}=N_{0,\zeta}$, in contradiction to the minimality of
$\epsilon$. So $\epsilon$ must be equal to 0.
Hence $x \in N_{0,0}=N_0$.\\
(b) Monotonicity: Let $\langle N_{0,\epsilon}:\epsilon<\lambda^+
\rangle,\ \langle N_{1,\epsilon}:\epsilon<\lambda^+ \rangle$ be
witnesses for $\widehat{NF}(N_0,N_1,M_0,M_1)$. Let $E$ be a club
of $\lambda^+$ such that $0 \notin E$ and $\epsilon \in E
\Rightarrow N_{1,\epsilon} \bigcap M_1^* \preceq N_{1,\epsilon}$
[Why do we have such a club? Let $E$ be a club of $\lambda^+$ such
that $0 \notin E$ and $\epsilon \in E \Rightarrow N_{1,\epsilon}
\bigcap M_1^* \preceq M_1^*$. By the assumption $M_1^* \preceq
M_1$. So $\epsilon \in E \Rightarrow N_{1,\epsilon} \bigcap M_1^*
\preceq M_1$. Now as $N_{1,\epsilon} \preceq M_1$, by axiom
\ref{1.1}.1.e $\epsilon \in E \Rightarrow N_{1,\epsilon} \bigcap
M_1^* \preceq N_{1,\epsilon}$]. We will prove that the sequences
$\langle N_0 \rangle ^\frown \langle N_{0,\epsilon}:\epsilon \in E
\rangle,\ \langle N_1^* \rangle^\frown \langle N_{1,\epsilon}
\bigcap M_1^*:\epsilon \in E \rangle$ witness that
$\widehat{NF}(N_0,N_1^*,M_0,M_1^*)$. First, they are increasing
[Why $\epsilon<\zeta \wedge \{\epsilon,\zeta \} \subseteq E
\Rightarrow N_{1,\epsilon} \bigcap M_1^* \preceq N_{1,\zeta}
\bigcap M_1^*$? By the properties of E, $N_{1,\epsilon} \bigcap
M_1^* \preceq N_{1,\epsilon}$. But $N_\epsilon \preceq N_\zeta$.
So $N_{1,\epsilon} \bigcap M_1^* \preceq N_{1,\zeta}$. In the
other side again by the properties of E, $N_{1,\epsilon} \bigcap
M_1^* \subseteq N_{1,\zeta} \bigcap M_1^* \preceq N_{1,\zeta}$. So
by axiom \ref{1.1}.1.e $N_{1,\epsilon} \bigcap M_1^* \preceq
N_{1,\zeta} \bigcap M_1^*$]. Second, we will prove that if
$\epsilon<\zeta,\ \{\epsilon,\zeta\} \subseteq E$ then
$NF(N_{0,\epsilon},N_{1,\epsilon} \bigcap
M_1^*,N_{0,\zeta},N_{1,\zeta} \bigcap M_1^*)$. Fix such
$\epsilon,\zeta$. By Theorem \ref{5.12}, (the long transitivity
theorem),
$NF(N_{0,\epsilon},N_{1,\epsilon},N_{0,\zeta},N_{1,\zeta})$. By
the properties of $E$ and axiom \ref{1.1}.1.e, $N_{0,\epsilon}
\preceq N_{1,\epsilon} \bigcap M_1^* \preceq N_{1,\epsilon},\
N_{0,\zeta} \allowbreak \bigcup (N_{1,\epsilon} \allowbreak
\bigcap M_1^*) \subseteq N_{1,\zeta} \bigcap M_1^* \preceq
N_{1,\zeta}$. Now by Theorem \ref{5.4b}.5 (the monotonicity of
NF), we have $NF(N_{0,\epsilon},N_{1,\epsilon}
\bigcap M_1^*,N_{0,\zeta},N_{1,\zeta} \bigcap M_1^*)$.\\
(c) Existence: By Proposition \ref{5.11}.b.\\
(d) Weak Uniqueness: Since $\bigotimes_{NF}$ holds, it follows by
Proposition \ref{examples of bigotimes-}.2 and Proposition
\ref{the transitivity of the weak uniqueness}. But we give another
proof using section 7: By Proposition \ref{7.4}.f, there is a
model $M_{1,n}^+$ such that $M_{1,n}\prec^+M_{1,n}^+$. By Theorem
\ref{7.13}.c, there is an isomorphism $f:M_{1,1}^+ \to M_{1,2}^+$
over $M_0 \bigcup N_1$. So $M_{1,2}^+,id_{M_{1,2}},f\restriction
M_{1,1}$ is a witness as
required.\\
(e) Let $\langle  N_{0,\epsilon}:\epsilon<\lambda^+ \rangle,\
\langle  N_{1,\epsilon}:\epsilon<\lambda^+ \rangle$ a witness for
$\widehat{NF}(N_0,N_1,M_0,M_1)$. There is $\epsilon$ such that $a
\in N_{0,\epsilon}$. By Definition \ref{5.14} (the definition of
$\widehat{NF}$), we have
$NF(N_0,N_1,N_{0,\epsilon},N_{1,\epsilon})$. So the proposition is
satisfied by Theorem \ref{5.10} (the relation $NF$ respects the
frame).
\end{proof}


\begin{definition} \label{6.1}
$M_0 \preceq^{NF} M_1$ when: there are $N_0,N_1$ such that
$\widehat{NF}(N_0,N_1,\allowbreak M_0,M_1)$.
\end{definition}

\begin{proposition} \label{6.2}\label{6.5}
$(K_{\lambda^+},\preceq^{NF})$ satisfies the following properties:
\begin{enumerate}[(a)]
\item Suppose $M_0 \preceq M_1,\ n<2 \Rightarrow M_n \in
K_{\lambda^+}$. For $n<2$ let $\langle
N_{n,\epsilon}:\epsilon<\lambda^+ \rangle$ be a representation of
$M_n$. Then $M_0 \preceq^{NF}M_1$ iff there is a club $E \subseteq
\lambda^+$ such that $(\epsilon<\zeta \wedge \{\epsilon,\zeta\}
\subseteq E) \Rightarrow
NF(N_{0,\epsilon},N_{0,\zeta},N_{1,\epsilon},N_{1,\zeta})$. \item
$\preceq^{NF}$ is a partial order. \item If $M_0 \preceq M_1
\preceq M_2$ and $M_0 \preceq^{NF} M_2$ then $M_0 \preceq^{NF}
M_1$. \item It satisfies axiom c of a.e.c. in $\lambda^+$, i.e.:
If $\delta \in \lambda^{+2}$ is a limit ordinal and $\langle
M_\alpha:\alpha<\delta \rangle$ is a $\preceq^{NF}$-increasing
continuous sequence, then $M_0 \preceq^{NF} \bigcup
\{M_\alpha:\alpha<\delta \}$ and obviously it is $\in
K_{\lambda^+}$. \item It has no $\preceq^{NF}$-maximal model.
\item If it satisfies smoothness (Definition \ref{1.1}.1.d), then
it is an a.e.c. in $\lambda^+$, (see Definition \ref{1.1}, page
\pageref{1.1}). \item LST for $\preceq^{NF}$: If $M_0 \preceq^{NF}
M_1,\ n<2 \Rightarrow (A_n \subseteq M_n \wedge |A_n| \leq
\lambda)$, then there are models $N_0,N_1 \in K_\lambda$ such
that: $\widehat{NF}(N_0,N_1,M_0,M_1)$ and $n<2 \Rightarrow
A_n\subseteq N_n$.
\end{enumerate}
\end{proposition}

\begin{proof}
\mbox{} (a) \case{One direction:} Let $E$ be such a club. So
$\langle N_{0,\epsilon}:\epsilon \in E \rangle,\ \langle
N_{1,\epsilon}:\epsilon \in E \rangle$ witness that $M_0
\preceq^{NF}M_1$.\\ \case{conversely:} Let $\langle
M_{0,\alpha}:\alpha<\lambda^+ \rangle,\ \langle
M_{1,\alpha}:\alpha<\lambda^+ \rangle$ be witnesses for $M_0
\preceq^{NF} M_1$. Let $E$ be a club such that $(n<2 \wedge
\epsilon \in E) \Rightarrow M_{n,\alpha}=N_{n,\alpha}$. Suppose
$\epsilon<\zeta \wedge \{\epsilon,\zeta\} \subseteq E$. We will
prove $NF(N_{0,\epsilon},N_{1,\epsilon},N_{0,\zeta},N_{1,\zeta})$,
i.e. $NF(M_{0,\epsilon},M_{1,\epsilon},M_{0,\zeta},M_{1,\zeta})$.
The sequences $\langle  M_{0,\alpha}:\epsilon \leq \alpha \leq
\zeta \rangle,\ \langle  M_{1,\alpha}:\epsilon \leq \alpha \leq
\zeta \rangle$ are increasing and continuous. So by Theorem
\ref{5.12} (the long transitivity theorem)
$NF(M_{0,\epsilon},M_{1,\epsilon},M_{0,\zeta},M_{1,\zeta})$.

(b) The reflexivity is obvious. The antisymmetry is satisfied by
the antisymmetry of the inclusion relation. The transitivity is
satisfied by item a, Theorem \ref{5.12} and the evidence that the
intersection of two clubs is a club.

(c) For $n<3$ let $\langle M_{n,\alpha}:\alpha<\lambda^+ \rangle$
be a representation of $M_n$ such that $\alpha<\lambda^+
\Rightarrow
NF(M_{0,\alpha},M_{0,\alpha+1},M_{2,\alpha},M_{2,\alpha+1})$. Let
E be a club of $\lambda^+$ such that $\alpha \in E \Rightarrow
M_{0,\alpha} \preceq M_{1,\alpha} \preceq M_{2,\alpha}$. By the
monotonicity of $NF$ $\alpha \in E \Rightarrow
NF(M_{0,\alpha},M_{0,\alpha+1},M_{1,\alpha},M_{1,\alpha+1})$. The
representations $\langle  M_{0,\alpha}:\alpha \in E \rangle,\
\langle  M_{1,\alpha}:\alpha \in E \rangle$ witness that $M_0
\preceq^{NF} M_1$.

(d) Without loss of generality $cf(\delta)=\delta$, so $\delta
\leq \lambda^+$. Denote $M_\delta:= \bigcup
\{M_\alpha:\alpha<\delta \}$. For $\alpha<\delta$ let $\langle
M_{\alpha,\epsilon}:\epsilon<\lambda^+ \rangle$ be a
representation of $M_n$. By item a for every $\alpha$ there is a
club $E_{\alpha,0} \subseteq \lambda^+$ such that
$(\epsilon<\zeta\wedge \{\epsilon,\zeta \} \subseteq E_{\alpha,0})
\Rightarrow
NF(M_{\alpha,\epsilon},M_{\alpha,\zeta},M_{\alpha+1,\epsilon},M_{\alpha+1,\zeta})$.
Let $\alpha$ be a limit ordinal. $\bigcup
\{M_{\alpha,\epsilon}:\epsilon<\lambda^+ \}=M_\alpha= \bigcup
\{M_\beta:\beta<\alpha \}= \bigcup \{ \bigcup \{
M_{\beta,\epsilon}:\epsilon<\lambda^+ \}:\beta<\alpha \} = \bigcup
\{\bigcup \{M_{\beta,\epsilon}:\beta<\alpha
\}:\epsilon<\lambda^+\}$. Every edge of this equivalences's
sequence is a limit of an $\subseteq$-increasing continuous
sequence of subsets of cardinality less than $\lambda$, and it is
equal to $M_\alpha$ [Why is the sequence in the right edge,
$\langle \bigcup \{M_{\beta,\epsilon}:\beta<\alpha
\}:\epsilon<\lambda^+ \rangle$ continuous? Let
$\epsilon<\lambda^+$ be a limit ordinal. Suppose $x \in \bigcup
\{M_{\beta,\epsilon}:\beta<\alpha \}$. Then there are
$\zeta,\beta$ such that $x \in M_{\beta,\zeta}$. So $x \in \bigcup
\{ M_{\beta,\zeta}:\beta<\alpha \}$]. So there is a club
$E_{\alpha,1} \subseteq \lambda^+$ such that $\epsilon \in
E_{\alpha,1} \Rightarrow M_{\alpha,\epsilon}=\bigcup \{
M_{\beta,\epsilon}:\beta<\alpha \}$. For $\alpha$ limit define
$E_\alpha:=E_{\alpha,0} \bigcap E_{\alpha,1}$, and for $\alpha$
not limit define $E_\alpha:=E_{\alpha,0}$.

\case{Case a:} $\delta<\lambda^+$. Define $E:=\bigcap
\{E_\alpha:\alpha<\delta \}$. If $\epsilon \in E$ then for
$\alpha<\delta$, $\epsilon \in E$, so
$NF(M_{\alpha,\epsilon},M_{\alpha,\\Min(E-(\epsilon+1))},
M_{\alpha+1,\epsilon},M_{\alpha+1,\Min(E-(\epsilon+1))})$. So be
Theorem \ref{5.12} (the long transitivity theorem), $\epsilon \in
E \Rightarrow NF(M_{0,\epsilon},\allowbreak
M_{0,\Min(E-(\epsilon+1))} ,\allowbreak
M_{\delta,\epsilon},M_{\delta,\Min(E-(\epsilon+1))})$. Hence $M_0
\preceq^{NF}M_1$.

\case{Case b:} \emph{$\delta=\lambda^+$}. Let $E:=\{\epsilon \in
E:\epsilon$ is a limit ordinal, $\alpha<\epsilon \Rightarrow
\epsilon \in E_\alpha \}$. Denote $N_\epsilon:= \bigcup
\{M_{\alpha,\epsilon}:\alpha<\epsilon \}$.

\begin{displaymath}
\xymatrix{M_0 \ar[r]^{id} &M_\alpha \ar[r]^{id} &M_\epsilon \ar[r]^{id} &M_\zeta \ar[r]^{id} &M_{\lambda^+}\\
M_{0,\zeta}  \ar[r]^{id} \ar[u]^{id} &M_{\alpha,\zeta}  \ar[r]^{id} \ar[u]^{id} &M_{\epsilon,\zeta}  \ar[r]^{id} \ar[u]^{id} &N_\zeta \ar[u]^{id}\\
M_{0,\epsilon} \ar[r]^{id} \ar[u]^{id} &M_{\alpha,\epsilon}  \ar[r]^{id} \ar[u]^{id}  & N_\epsilon \ar[u]^{id}\\
M_{0,\alpha} \ar[r]^{id} \ar[u]^{id} &N_\alpha \ar[u]^{id} \\
M_{0,0} \ar[u]^{id}  }
\end{displaymath}

\begin{claim} \label{*}
For every $\epsilon \in E$ the sequence $\langle
M_{\alpha,\epsilon}:\alpha<\epsilon \rangle ^\frown \langle
N_\epsilon \rangle$ is increasing and continuous (especially
$N_\epsilon \in K$),
\end{claim}

\begin{proof}
If $\epsilon \in E$ is limit, then $\alpha<\epsilon \Rightarrow
\epsilon \in E_{\alpha,1}$, so the sequence $\langle
M_{\alpha,\epsilon}:\alpha<\epsilon \rangle$ is continuous. So it
is sufficient to prove that $\alpha<\epsilon \Rightarrow
M_{\alpha,\epsilon} \preceq M_{\alpha,\epsilon+1}$. Suppose
$\alpha<\epsilon$. $\epsilon \in E$, so $\epsilon \in
E_{\alpha,0}$. Hence
$NF(M_{\alpha,\epsilon},M_{\alpha+1,\epsilon},
M_{\alpha,\Min(E-(\epsilon+1))},\allowbreak
M_{\alpha+1,\Min(E-(\epsilon+1))})$, and especially
$M_{\alpha,\epsilon} \preceq M_{\alpha+1,\epsilon}$.
\end{proof}

\begin{claim} \label{**}
The sequence $\langle N_\epsilon:\epsilon \in E \rangle$ is
$\preceq$-increasing.
\end{claim}

\begin{proof}
Suppose $\epsilon<\zeta,\
\{\epsilon,\zeta\} \subseteq E$. By (*), the sequences $\langle
M_{\alpha,\epsilon}:\alpha<\epsilon \rangle ^\frown \langle
N_\epsilon \rangle, \langle M_{\alpha,\zeta}:\alpha \leq \epsilon
\rangle$ are increasing and continuous. For every $\alpha \in
\epsilon$ the sequence $\langle M_{\alpha,\beta}:\beta<\lambda^+
\rangle$ is a representation of $M_\alpha$, and especially it is
$\preceq$-increasing. So $(\forall \alpha \in \epsilon)
M_{\alpha,\epsilon} \preceq M_{\alpha,\zeta}$. Hence by the
smoothness $N_\epsilon \preceq M_{\epsilon,\zeta}$. But by (*),
$M_{\epsilon,\zeta} \preceq N_\zeta$, so $N_\epsilon \preceq N_\zeta$.]\\
\end{proof}

\begin{claim} \label{***} The sequence $\langle N_\epsilon:\epsilon \in E \rangle$ is
continuous
\end{claim}

\begin{proof}
Suppose $\epsilon=sup(E \bigcap \epsilon)$. Let $x \in
N_\epsilon$. By the definition of $N_\epsilon$ there is
$\alpha<\epsilon$ such that $x \in M_{\alpha,\epsilon}$.
$\epsilon$ is limit and the sequence $\langle
M_{\alpha,\beta}:\beta \leq \epsilon \rangle$ is continuous. So
there is $\beta<\epsilon$ such that $x \in M_{\alpha,\beta}$.
$\epsilon=sup(E \bigcap \epsilon)$, so there is $\zeta \in
(\beta,\epsilon) \bigcap E$. $x \in M_{\alpha,\zeta}$ but by (*),
$M_{\alpha,\zeta} \subseteq N_\zeta$, so $x \in N_\zeta$].
\end{proof}

\begin{claim} \label{****}
$\bigcup \{N_\epsilon:\epsilon \in E \}=M_\delta$
\end{claim}

\begin{proof}
Clearly $\bigcup \{N_\epsilon:\epsilon \in E \} \subseteq
M_\delta$. The other inclusion: Let $x \in M_\delta$. Then there
is $\alpha<\delta$ such that $x \in M_\alpha$. So $(\exists
\alpha,\beta)x \in M_{\alpha,\beta}$. So as $sup(E)=\delta$, There
is $\zeta \in (\beta,\delta) \bigcap E$. So $x \in
M_{\alpha,\zeta}$ which by (*) is $\subseteq N_\zeta$. So $x \in
N_\zeta$].
\end{proof}

\begin{claim} \label{*****}
If $\epsilon<\zeta,\ \{\epsilon,\zeta \} \subseteq E$ then
$NF(M_{0,\epsilon},N_\epsilon,M_{0,\zeta},N_\zeta)$
\end{claim}

\begin{proof}
By the definition of E, $(\forall \alpha \in \epsilon)
\{\epsilon,\zeta \} \subseteq E_\alpha$. So $(\forall \alpha \in
\epsilon) NF(M_{\alpha,\epsilon},\allowbreak
M_{\alpha+1,\epsilon},M_{\alpha,\zeta},\allowbreak
M_{\alpha+1,\zeta})$. By (*), the sequences $\langle
M_{\alpha,\epsilon}:\alpha<\epsilon \rangle ^\frown \langle
N_\epsilon \rangle,\ \langle M_{\alpha,\zeta}:\alpha \leq \epsilon
\rangle$ are increasing and continuous. So by Theorem \ref{5.12}
(the long transitivity theorem),
$NF(M_{0,\epsilon},N_\epsilon,M_{0,\zeta},M_{\epsilon,\zeta})$.
But by Claim \ref{*} $M_{\epsilon,\zeta} \prec N_\zeta$, so
$NF(M_{0,\epsilon},N_\epsilon,M_{0,\zeta},\allowbreak N_\zeta)$].
\end{proof}

By Claims \ref{**},\ref{***},\ref{****}, the sequence $\langle
N_\epsilon:\epsilon<\delta \rangle$ is a representation of
$M_\delta$. The sequence $\langle
M_{0,\epsilon}:\epsilon<\lambda^+ \rangle$ is a representation of
$M_0$. Hence, by Claim \ref{*****} and item a, they witness that
$M_0 \preceq^{NF}M_\delta$.

(e) By Proposition \ref{the widehat{NF}-properties}.c. Derived
also by the existence proposition of the $\prec^+$-extension,
(Proposition \ref{7.4}.f), which we will prove later.

(f) We have actually proved it, (for example: axiom \ref{1.1}.1.e
by item c here and axiom \ref{1.1}.1.c. By item d here).

(g) Let $\langle N_{0,\epsilon}:\epsilon<\lambda^+ \rangle,\
\langle N_{1,\epsilon}:\epsilon<\lambda^+ \rangle$ be witnesses
for $M_0 \preceq^{NF}M_1$. By cardinality considerations there is
$\epsilon \in \lambda^+$ such that for $n<2$ we have $A_n
\subseteq N_{n,\epsilon}$. But for every $\epsilon<\lambda^+$,
$\widehat{NF}(N_{0,\epsilon},N_{1,\epsilon},M_0,M_1).$
\end{proof}

\section{$\prec^+$ and saturated models}
\begin{hypothesis}
$\frak{s}$ is a weakly successful semi-good $\lambda$-frame with
conjugation.
\end{hypothesis}

\begin{definition} \label{7.1}
$K^{sat}$ is the class of saturated models in $\lambda^+$ over
$\lambda$.
\end{definition}

Now we study the class $(K^{sat},\preceq^{NF}\restriction
K^{sat})$. Note that in the following theorem there is no any
set-theoretic hypothesis beyond $ZFC$.

\begin{theorem} \label{9.6}
If ($\frak{s}$ is a weakly successful semi-good $\lambda$-frame
with conjugation and) $(K^{sat},\preceq^{NF} \restriction
K^{sat})$ does not satisfy smoothness (see Definition
\ref{1.1}.1.d), then there are $2^{\lambda^{+2}}$ pairwise
non-isomorphic models in $K_{\lambda^{+2}}$.
\end{theorem}

\discussion{How can we prove this theorem?} First we find a
relation $\prec^+$ on $K_{\lambda^+}$ such that:
\begin{itemize}
 \item[(*)] For every model $M_0$ in
$K_{\lambda^+}$ there is a model $M_1$ such that $M_0\prec^+M_1$.
\item[(**)] If for $n=1,2$ $M_0\prec^+M_n$ then $M_1,M_2$ are
isomorphic over $M_0$. \item[(***)] If $\langle M_i:i \leq
\alpha^* \rangle$ is an increasing continuous sequence, and
$i<\alpha^* \Rightarrow M_i\prec^+M_{i+1}$ then
$M_0\prec^+M_{\alpha^*}$.
\end{itemize}

In section 7 we study the properties of $\prec^+$. Sections 8,9
are preparations for the proof of Theorem \ref{9.6}. A key theorem
is Theorem \ref{9.4}: Suppose that there is an increasing
continuous sequence $\langle M^*_\alpha:\alpha \leq \lambda+1
\rangle$ of models in $K^{sat}$ such that: $\alpha<\beta<\lambda^+
\Rightarrow M^*_\alpha \prec^+ M^*_\beta \wedge M^*_\alpha
\preceq^{NF}M_{\lambda^++1}$ and $M^*_{\lambda^+} \npreceq^{NF}
M^*_{\lambda^++1}$. \underline{Then} for every $S \in
S^{\lambda^{+2}}_{\lambda^+}:=\{S:S$ is a stationary subset of
$\lambda^{+2}$ and $(\forall \alpha \in S)cf(\alpha)=\lambda^+
\}$, there is a model $M^S$ in $K_{\lambda^{+2}}$ such that
$S(M^S)=S/D_{\lambda^{+2}}$. So there are $2^{\lambda^{+2}}$
pairwise non-isomorphic models in $K_{\lambda^{+2}}$.

Note that while $\prec^+$ is a priori defined on $K_{\lambda^+}$,
Proposition \ref{we can decrease E} shows that any $\prec^+$
extension is saturated in $\lambda^+$ over $\lambda$, so in
$K^{sat}$.

\begin{definition} \label{7.3}
$\prec^+$ is a 2-place relation on $K_{\lambda^+}$. For $M_0,M_1
\in K_{\lambda^+}$, we say $M_0\prec^+M_1$ iff: there are
increasing continuous sequences $\langle
N_{0,\alpha}:\alpha<\lambda^+ \rangle$, $\langle
N_{1,\alpha}:\alpha<\lambda^+ \rangle$, $\langle
N^{\oplus}_{1,\alpha}:\alpha<\lambda^+ \rangle$, and there is a
club $E$ of $\lambda^+$ such that:
\begin{enumerate}[(a)]
\item For $n=1,2$ $M_n=\bigcup \{N_{n,\alpha}:\alpha<\lambda^+\}$.
\item $\alpha \in E \Rightarrow N_{0,\alpha} \preceq N_{1,\alpha}
\preceq N^{\oplus}_{1,\alpha}$. \item If $\alpha<\beta$ and they
are in $E$, then
$NF(N_{0,\alpha},N^{\oplus}_{1,\alpha},N_{0,\beta},N_{1,\beta})$.
 \item For every $\alpha \in E$,
and every $p \in S^{bs}(N_{1,\alpha})$, there is an end-segment
$S$ of $\lambda^+$ such that for every $\beta \in S \bigcap E$ the
model $N^{\oplus}_{1,\beta}$ realizes the non-forking extension of
$p$ to $N_{1,\beta}$.
\end{enumerate}
In such a case $\langle N_{0,\alpha}:\alpha<\lambda^+ \rangle$,
$\langle N_{1,\alpha}:\alpha<\lambda^+ \rangle$, $\langle
N_{1,\alpha}^{\oplus}:\alpha<\lambda^+ \rangle$, $E$ are said to
be witnesses for $M_0 \prec^+ M_1$.
\end{definition}

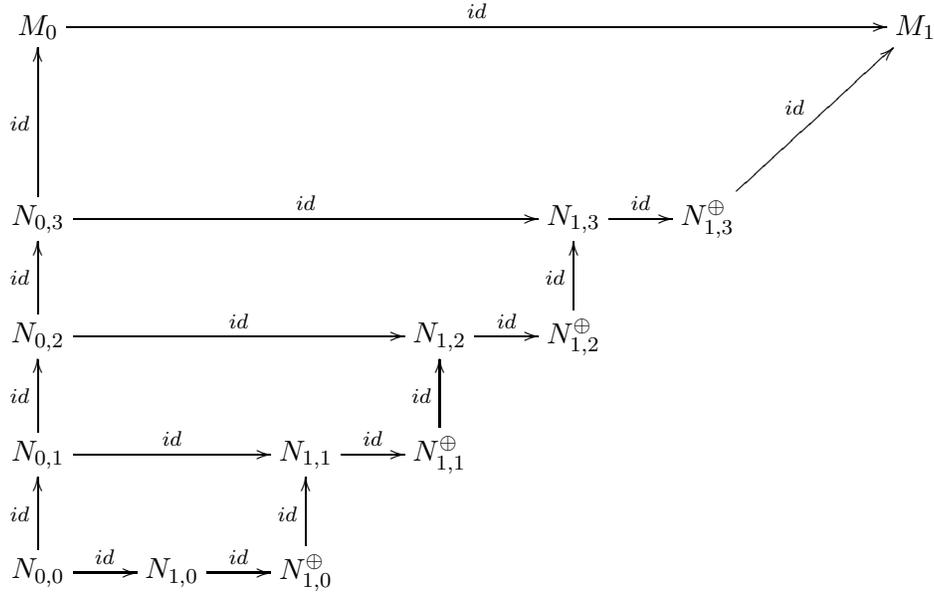
\begin{displaymath}
\xymatrix{
M_0 \ar[rrrrrrr]^{id} &&&&&&&M_1\\
\\
N_{0,3} \ar[rrrr]^{id} \ar[uu]^{id} &&&&N_{1,3} \ar[r]^{id}
&N_{1,3}^{\oplus}
 \ar[rruu]^{id} \\
N_{0,2} \ar[rrr]^{id} \ar[u]^{id} &&&N_{1,2} \ar[r]^{id}&N_{1,2}
^{\oplus} \ar[u]^{id}\\
N_{0,1} \ar[rr]^{id} \ar[u]^{id} && N_{1,1} \ar[r]^{id} & N_{1,1}^{\oplus} \ar[u]^{id} \\
N_{0,0} \ar[r]^{id} \ar[u]^{id} & N_{1,0} \ar[r]^{id} &
N^{\oplus}_{1,0} \ar[u]^{id}  }
\end{displaymath}

By the following proposition if $M_0 \prec^+ M_1$ then we can find
witnesses for it, with $E=\lambda^+$.
\begin{proposition}\label{wlog E=lambda^+}
If
\begin{enumerate} \item $\langle N_{0,\alpha}:\alpha<\lambda^+
\rangle$, $\langle N_{1,\alpha}:\alpha<\lambda^+ \rangle$,
$\langle N_{1,\alpha}^{\oplus}:\alpha<\lambda^+ \rangle$, $E$ are
witnesses for $M_0 \prec^+ M_1$. \item For $\alpha \in E$
$M_{0,otp(\alpha \bigcap E)}=N_{0,\alpha}$, $M_{1,otp(\alpha
\bigcap E)}=N_{1,\alpha}$, $M_{1,otp(\alpha \bigcap
E)}^{\oplus}=N_{1,\alpha}^{\oplus}$.
\end{enumerate}
then $\langle M_{0,\beta}:\beta<\lambda^+ \rangle$, $\langle
M_{1,\beta}:\beta<\lambda^+ \rangle$, $\langle
M_{1,\beta}^{\oplus}:\beta<\lambda^+ \rangle$, $\lambda^+$ are
witnesses for $M_0 \prec^+ M_1$.
\end{proposition}

\begin{proof}
Easy, so we prove Definition \ref{7.3}.c only. Suppose
$\gamma_0<\gamma_1$. We have to prove that
$NF(M_{0,\gamma_0},M_{1,\gamma_0}^{\oplus},M_{0,\gamma_1},M_{1,\gamma_1})$.
There is a unique ordinal $\alpha \in E$ with $otp(\alpha \bigcap
E)=\gamma_0$. So $M_{0,\gamma_0}=N_{0,\alpha} \wedge
M_{1,\gamma_0}^{\oplus}=N_{1,\alpha}^{\oplus}$. Similarly there is
a unique $\beta \in E$ such that $M_{0,\gamma_1}=N_{0,\beta}
\wedge M_{1,\gamma_1}=N_{1,\beta}$. Now by clause b in the
assumption
$NF(N_{0,\alpha},N_{1,\alpha}^{\oplus},N_{0,\beta},N_{1,\beta})$,
namely
$NF(M_{0,\gamma_0},M_{1,\gamma_0}^{\oplus},N_{0,\gamma_1},N_{1,\gamma_1})$.
\end{proof}

\begin{proposition} \label{we can decrease E}
If $\langle N_{0,\alpha}:\alpha<\lambda^+ \rangle$, $\langle
N_{1,\alpha}:\alpha<\lambda^+ \rangle$, $\langle
N_{1,\alpha}^{\oplus}:\alpha<\lambda^+ \rangle$, $E$ are witnesses
for $M_0 \prec^+ M_1$ and $E^-$ is a club of $\lambda^+$ with $E^-
\subseteq E$ then $\langle N_{0,\alpha}:\alpha<\lambda^+ \rangle$,
$\langle N_{1,\alpha}:\alpha<\lambda^+ \rangle$, $\langle
N_{1,\alpha}^{\oplus}:\alpha<\lambda^+ \rangle$, $E^-$ are
witnesses for $M_0 \prec^+ M_1$.
\end{proposition}

\begin{proof}
Trivial.
\end{proof}

\begin{proposition} \label{7.2}
Suppose:
\begin{enumerate}[(a)]
\item For $n=1,2$ $NF(M_{0,0},M_{0,1},M_{n,0},M_{n,1})$. \item
$M_{1,0} \preceq N_0,\ M_{2,0} \preceq N_0$. \item $N_0 \bigcap
M_{0,1}=M_{0,0}$.
\end{enumerate}

Then for some model $N_1$ with $NF(M_{0,0},M_{0,1},N_0,N_1)$ we
can assign to each $n \in \{1,2\}$ an embedding $f_n:M_{n,1} \to
N_1$ over $M_{0,1} \bigcup M_{n,0}$ such that
$NF(M_{n,0},f_n[M_{n,1}],N_0,N_1)$.
\end{proposition}

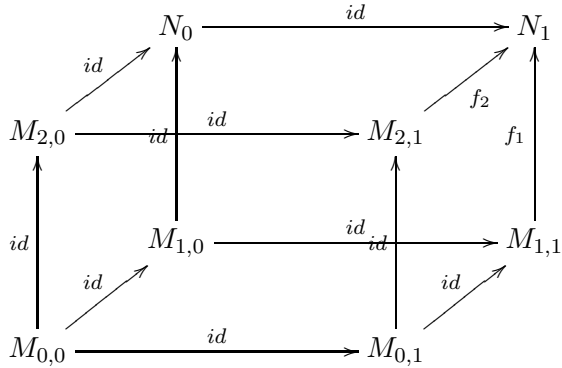
\begin{displaymath}
\xymatrix{
&N_0 \ar[rrr]^{id}&&&N_1\\
M_{2,0} \ar[ur]^{id} \ar[rrr]^{id} &&&M_{2,1} \ar[ur]_{f_2}\\
&M_{1,0} \ar[uu]^{id} \ar[rrr]^{id} &&& M_{1,1} \ar[uu]^{f_1}\\
M_{0,0} \ar[uu]^{id} \ar[ur]^{id} \ar[rrr]^{id} &&& M_{0,1}
\ar[uu]^{id} \ar[ur]^{id}}
\end{displaymath}

\begin{proof}
For each $n \in \{1,2\}$ by Theorem \ref{the existence theorem for
NF} (the existence theorem for $NF$), we can find an amalgamation
$(id_{N_0},g_n,N_{n,1})$ of $N_0,M_{n,1}$ over $M_{n,0}$ with
$NF(M_{n,0},N_0,g_n[M_{n,1}],N_{n,1})$. But
$NF(M_{0,0},M_{n,0},M_{0,1},M_{n,1})$. So by Theorem \ref{5.12}
(the long transitivity theorem)
$NF(M_{0,0},N_0,g_n[M_{0,1}],N_{n,1})$. By assumption c $N_0
\bigcap M_{0,1}=M_{0,0}$. So by Theorem \ref{the weak uniqueness
theorem of NF} (the weak uniqueness theorem) we can find
$h_1,h_2,N_1$ such that the following hold:
\begin{enumerate}
\item $h_n:N_{n,1} \to N_1$ is an embedding. \item $h_n
\restriction N_0=id_{N_0}$. \item $h_1 \circ g_1 \restriction
M_{0,1}=h_1 \circ g_2 \restriction M_{0,1}=id_{M_{0,1}}$.
\end{enumerate}
Now we define for $n=1,2$ $f_n:=h_n \circ g_n$. Why is $f_n$ over
$M_{0,1} \bigcup M_{n,0}$? By clause 3 $x \in M_{0,1} \Rightarrow
f_n(x)=x$. Let $x \in M_{n,0}$. Then $g_n(x)=x$. By assumption b
$M_{n,0} \subseteq N_0$, so $x \in N_0$. So by clause 2
$h_n(x)=x$. Hence $f_n(x)=h_n(g_n(x))=h_n(x)=x$.
\begin{claim}
$NF(M_{n,0},f_n[M_{n,1}],N_0,N_1)$.
\end{claim}

\begin{proof}
$NF(M_{n,0},N_0,g_n[M_{n,1}],N_{n,1})$. So by clauses 1,2
$NF(M_{n,0},N_0,\allowbreak f_n \allowbreak
[M_{n,1}],h_n[N_{n,1}])$. But $h_n[N_{n,1}] \preceq N_1$, so
$NF(M_{n,0},N_0,f_n[M_{n,1}],N_1)$.
\end{proof}

\begin{claim}
$NF(M_{0,0},M_{0,1},N_0,N_1)$.
\end{claim}

\begin{proof}
Since $NF(M_{1,0},M_{1,1},N_0,N_1)$, by Theorem \ref{5.12} (the
long transitivity theorem) it is enough to prove that
$NF(M_{0,0},M_{0,1},M_{1,0},f_1[M_{1,1}])$. But $f_n$ is over
$M_{0,1} \bigcup M_{1,0}$. Hence it follows by assumption a.
\end{proof}
\end{proof}

\begin{proposition}\label{before 7.4} \label{7.10 in 20.11.09}
\mbox{}
\begin{enumerate}[(a)]
\item If $M_0\prec^+ M_1$ then $M_0 \prec^{NF}M_1$. \item If
$M_0\prec^+M_1$ then $M_1 \in K^{sat}$. \item If $M_0 \preceq^{NF}
M_1\prec^+M_2$ then $M_0\prec^+M_2$. \item If $M_0 \prec^+ M_1
\prec^+ M_2$ then $M_0 \prec^+ M_2$.
\end{enumerate}
\end{proposition}

\begin{proof}
\mbox{}
\begin{enumerate}[(a)]
\item If $\langle N_{0,\alpha}:\alpha<\lambda^+ \rangle,\ \langle
N_{1,\alpha}:\alpha<\lambda^+ \rangle,\ \langle
N_{2,\alpha}:\alpha<\lambda^+ \rangle,\ E$ witness that $M_0
\prec^+ M_1$ then $\langle N_{0,\alpha}:\alpha \in E \rangle,\
\langle N_{1,\alpha}:\alpha \in E \rangle$ witness that
$\widehat{NF}(N_{0,0},N_{1,0},M_0,M_1)$. So $M_0 \preceq^{NF}M_1$.

\item By Theorem \ref{4.0}.2 (page \pageref{4.0}).

\item Easy.

\item By items a,c.
\end{enumerate}
\end{proof}

\begin{definition}
The \emph{$\prec^+$-game} is a game between two players. It lasts
$\lambda^+$ moves. In any move the players choose models in
$K_\lambda$ with the following rules:

The $0$ move: Player 1 chooses models $N_{0,0},N_{1,0} \in
K_\lambda$ with $N_{0,0} \preceq N_{1,0}$ and player 2 does not do
anything.

The $\alpha$ move where $\alpha$ is limit: Player 1 must choose
$N_{0,\alpha}:=\bigcup \{N_{0,\beta}:\beta<\alpha\}$ and Player 2
must choose $N_{1,\alpha}:=\bigcup \{N_{1,\beta}:\beta<\alpha\}$.

The $\alpha+1$ move: Player 1 chooses a model $N_{0,\alpha+1}$
such that the following hold:
\begin{enumerate}
\item $N_{0,\alpha} \preceq N_{0,\alpha+1}$. \item $N_{0,\alpha+1}
\bigcap N_{1,\alpha}=N_{0,\alpha}$.
\end{enumerate}

After player one chooses $N_{0,\alpha+1}$, player 2 has to choose
$N_{1,\alpha+1}$ such that the following hold:
\begin{enumerate}
\item $N_{1,\alpha} \preceq N_{1,\alpha+1}$. \item
$NF(N_{0,\alpha},N_{1,\alpha},N_{0,\alpha+1},N_{1,\alpha+1})$.
\end{enumerate}

In the end of the game, player 2 wins the game if $\bigcup
\{N_{0,\alpha}:\alpha<\lambda^+\}\prec^+ \bigcup
\{N_{1,\alpha}:\alpha<\lambda^+ \}$.

A \emph{strategy for player 2} is a function $F$ that assigns a
model $N_{1,\alpha+1}$ to each triple $(\alpha,\langle
N_{0,\beta}:\beta \leq \alpha+1 \rangle,\langle N_{1,\beta}:\beta
\leq \alpha \rangle)$ that satisfies the following conditions:
\begin{enumerate}
\item $\alpha<\lambda^+$. \item $\langle N_{0,\beta}:\beta \leq
\alpha+1 \rangle,\ \langle N_{1,\beta}:\beta \leq \alpha \rangle$
are increasing continuous sequences of models in $K_\lambda$.
\item
$NF(N_{0,\alpha},N_{1,\alpha},N_{0,\alpha+1},N_{1,\alpha+1})$ for
$\beta<\alpha$. \item $N_{0,\alpha+1} \bigcap
N_{1,\alpha}=N_{0,\alpha}$.
\end{enumerate}
 A \emph{winning strategy for player 2} is a
strategy for player 2, such that if player 2 acts by it, then he
wins the game, no matter what does player 1 do.

\end{definition}

\begin{proposition} \label{7.4}
\mbox{}
\begin{enumerate}[(a)]
 \item For every $M_0 \in
K_{\lambda^+}$ there is $M_1$ with $M_0\prec^+M_1$. \item If $M_0
\in K_{\lambda^+},\ n<2 \Rightarrow N_n \in K_\lambda,\ N_0\prec
M_0,\ N_0\prec N_1,\ N_1 \bigcap M_0=N_0$, then there is $M_1$
such that $M_0\prec^+M_1$ and $\widehat{NF}(N_0,N_1,M_0,M_1)$.
\item Player 2 has a winning strategy in the $\prec^+$-game.
\end{enumerate}
\end{proposition}

\begin{proof}
\mbox{}
(a) By c.\\
(b) By c.\\
(c) 
We describe a strategy: For $\alpha=0$ player 2 has nothing to do,
but he takes a paper and writes for himself: I define
$N_{1,0}^{temp}:=N_{1,0}$. For $\alpha$ limit player 2 chooses
$N_{1,\alpha}:=\bigcup \{N_{1,\beta}:\beta<\alpha \}$ and writes
for himself $N_{1,\alpha}^{temp}:=N_{1,\alpha}$. In the $\alpha+1$
move, he writes for himself 3 things:
\begin{enumerate}[(i)] \item A model $N_{1,\alpha+1}^{temp}$ with
$NF(N_{0,\alpha},N_{1,\alpha},N_{0,\alpha+1},N_{1,\alpha+1}^{temp})$.
By Theorem \ref{the existence theorem for NF}.a (on page
\pageref{the existence theorem for NF}) it is possible. \item A
sequence of types $\langle p_{\alpha,\beta}:\beta<\lambda^+
\rangle$ that includes $S^{bs}(N_{1,\alpha}^{temp})$.
\end{enumerate}
Now player 2 chooses a model $N_{1,\alpha+1}$ such that the
following hold:
\begin{enumerate}
\item $N_{1,\alpha+1}^{temp} \preceq N_{1,\alpha+1}$. \item For
each type in $p_{\gamma,\beta}$ with $\gamma<\alpha,\
\beta<\alpha$, $N_{1,\alpha+1}$ realizes the non-forking extension
of $p_{\gamma,\beta}$ over $N_{1,\alpha+1}^{temp}$.
\end{enumerate}
By Proposition \ref{1.8} (page \pageref{1.8}) it is possible.

Why shall player 2 win the game? Substitute the sequences $\langle
N_{0,\alpha}:\alpha<\lambda^+ \rangle,\ \langle
N_{1,\alpha}^{temp}:\alpha<\lambda^+ \rangle,\ \langle
N_{1,\alpha}:\alpha<\lambda^+ \rangle$ which appear here instead
of the sequences $\langle N_{0,\alpha}:\alpha<\lambda^+ \rangle,\
\langle N_{1,\alpha}:\alpha<\lambda^+ \rangle$, $\langle
N^{\oplus}_{1,\alpha}:\alpha<\lambda^+ \rangle$ in Definition
\ref{7.3}, and substitute $E=\lambda^+$.
\end{proof}

Roughly the following theorem says that:
\begin{enumerate}[(a)]
\item The $\prec^+$-extension is unique. \item Locality: Every
type over a model in $K_{\lambda^+}$ is determined by its
restrictions to submodels in $K_\lambda$. \item A preparation for
symmetry.
\end{enumerate}
\begin{theorem} \label{7.13}
Suppose for $n=1,2$ $M_0\prec^+M_n$ then:
\begin{enumerate}[(a)]
\item $M_1,M_2$ are isomorphic over $M_0$.\\
\item For every $a_1 \in M_1$, $a_2 \in M_2$ if for each $N \in
K_\lambda$ with $N \preceq M_0$ $tp(a_1,N,M_1)=tp(a_2,N,M_2)$ then
there is an isomorphism $f:M_1 \to M_2$ over $M_0$ with
$f(a_1)=a_2$. \item Let $N^* \in K_\lambda$, $N_0 \preceq N^*$. If
for $n=1,2$ $\widehat{NF}(N_0,N^*,M_0,M_n)$, then there is an
isomorphism $f:M_1 \to M_2$ over $M_0 \bigcup N^*$.
\end{enumerate}
\end{theorem}

\discussion{The plan of the proof:} We prove the three items at
once. The proof is similar to that of the uniqueness of the
saturated model in $\lambda^+$ over $\lambda$. Suppose $\langle
N_{0,\epsilon}:\epsilon<\lambda^+ \rangle,\ \langle
N_{1,\epsilon}:\epsilon<\lambda^+ \rangle,\ \langle
N_{1,\epsilon}^{\oplus}:\epsilon<\lambda^+ \rangle,\ \lambda^+$
witness that $M_0 \prec^+ M_1$. So $\langle
N_{0,\epsilon}:\epsilon<\lambda^+ \rangle$ is a representation of
$M_0$ and $\langle
N_{1,0},N_{1,0}^{\oplus},N_{1,1},N_{1,1}^{\oplus},...N_{1,\omega},N_{1,\omega}^{\oplus}...
\rangle$ is a representation of $M_1$. Suppose in addition that
$\langle N_{0,\epsilon}:\epsilon<\lambda^+ \rangle,\ \langle
N_{2,\epsilon}:\epsilon<\lambda^+ \rangle,\ \langle
N_{2,\epsilon}^{\oplus}:\epsilon<\lambda^+ \rangle,\ \lambda^+$
witness that $M_0 \prec^+ M_2$. We amalgamate $M_1,M_2$ over $M_0$
in $\lambda^+$ steps. In each step we amalgamate the corresponding
models in the representations of $M_1,M_2$ over the corresponding
model in the representation of $M_0$. Now if $(f_1,f_2,M_3)$ is an
amalgamation of $M_1,M_2$ over $M_0$ and $f_1,f_2$ are onto $M_3$,
then $f_2^{-1} \circ f_1$ is an isomorphism of $M_1$ into $M_2$
over $M_0$ as required. In odd steps we choose the amalgamations
such that in the end $f_1,f_2$ will be onto $M_3$, see requirement
8 below. In even steps we choose amalgamations with $NF$, see
requirement 4 below.

\begin{proof}
Roughly, the following claim says that one representation of $M_0$
can serve as a part of the witness to $M_0 \prec^+ M_1$ and $M_0
\prec^+ M_2$ together.
\begin{claim}
There are sequences $\langle N_{0,\epsilon}:\epsilon<\lambda^+
\rangle ,\ \langle N_{1,\epsilon}:\epsilon<\lambda^+ \rangle,\
\langle N^{\oplus}_{1,\epsilon}:\epsilon<\lambda^+ \rangle,\
\langle N_{2,\epsilon}:\epsilon<\lambda^+ \rangle,\ \langle
N^{\oplus}_{2,\epsilon}:\epsilon<\lambda^+ \rangle$ such that for
$n=1,2$ $\langle N_{0,\epsilon}:\epsilon<\lambda^+ \rangle,\
\langle N_{n,\epsilon}:\epsilon<\lambda^+ \rangle,\ E=\lambda^+,\
\langle N^{\oplus}_{n,\epsilon}:\epsilon<\lambda^+ \rangle$
witnesses that $M_0\prec^+M_n$ (so $\bigcup
\{N_{0,\epsilon}:\epsilon<\lambda^+\}=M_0$ and for $n=1,2$
$\bigcup \{N_{n,\epsilon}:\epsilon<\lambda^+\}=\bigcup
\{N_{n,\epsilon}^{\oplus}:\epsilon<\lambda^+ \}=M_n$).
\end{claim}

\begin{proof}
For $n=1,2$ we take witnesses $\langle
N_{0,n,\epsilon}^{temp}:\epsilon<\lambda^+ \rangle ,\ \langle
N_{n,\epsilon}^{temp}:\epsilon<\lambda^+ \rangle,\ \langle
N^{\oplus,temp}_{n,\epsilon}:\epsilon<\lambda^+ \rangle$, $E_n$
for $M_0 \prec^+ M_n$. Take a club $E$ of $\lambda^+$ such that $E
\subseteq E_1 \bigcap E_2$ and $\epsilon \in E \Rightarrow
N_{0,1,\epsilon}^{temp}=N_{0,2,\epsilon}^{temp}$. By Proposition
\ref{we can decrease E} for $n=1,2$ $\langle
N_{0,n,\epsilon}^{temp}:\epsilon<\lambda^+ \rangle ,\ \langle
N_{n,\epsilon}^{temp}:\epsilon<\lambda^+ \rangle,\ \langle
N^{\oplus,temp}_{n,\epsilon}:\epsilon<\lambda^+ \rangle$, $E$ are
witnesses for $M_0 \prec^+ M_n$. Define $N_{0,otp(\epsilon \bigcap
E)}:=N_{0,1, \epsilon)}^{temp}$. For $n=1,2$ and $\epsilon \in E$,
define $N_{n,otp(\epsilon \bigcap E)}:=N_{n,\epsilon}^{temp}$. By
Proposition \ref{wlog E=lambda^+} for $n=1,2$ $\langle
N_{0,\epsilon}:\epsilon<\lambda^+ \rangle,\ \langle
N_{n,\epsilon}:\epsilon<\lambda^+ \rangle,\ E=\lambda^+,\ \langle
N^{\oplus}_{n,\epsilon}:\epsilon<\lambda^+ \rangle$ witness that
$M_0\prec^+M_n$.
\end{proof}

For item b, we require in addition that $a_n \in N_{n,0}$ and
$tp(a_1,N_{0,0},N_{1,0})=tp(a_2,N_{0,0},N_{2,0})$. For item c, we
require in addition that $NF(N_0,N^*,N_{0,0},\allowbreak
N_{n,0})$.

Define by induction on $\epsilon \leq \lambda^+$ a triple
$(N_\epsilon, f_{1,\epsilon},f_{2,\epsilon})$ such that:
\begin{enumerate}
\item $\langle N_\epsilon:\epsilon \leq \lambda^+ \rangle$ is an
increasing continuous sequence of models in $K_\lambda$ and for
every $\epsilon<\lambda^+$ $N_{2\epsilon} \bigcap
M_0=N_{2\epsilon+1} \bigcap M_0=N_{0,\epsilon}$. \item For item c
we add: $f_{n,0}\restriction N^*$ is the identity. \item For item
b we add: $f_{1,0}(a_1)=f_{2,0}(a_2)$. \item $\epsilon<\lambda^+
\Rightarrow
NF(N_{0,\epsilon},N_{2\epsilon+1},N_{0,\epsilon+1},N_{2\epsilon+2})$.
\item For $n=1,2$ the sequence $\langle f_{n,\epsilon}:\epsilon
\leq \lambda^+ \rangle$ is increasing and continuous. \item For
$\epsilon<\lambda^+,\ f_{n,2\epsilon}$ is an embedding of
$N_{n,\epsilon}$ to $N_{2\epsilon}$ and $f_{n,2\epsilon+1}$ is an
embedding of $N^{\oplus}_{n,\epsilon}$ to $N_{2\epsilon+1}$. \item
$f_{n,2\epsilon}\restriction
N_{0,\epsilon}=f_{n,2\epsilon+1}\restriction N_{0,\epsilon}$ and
it is the identity on $N_{0,\epsilon}$. \item For every
$\epsilon<\lambda^+$ if for some $n \in \{1,2\}$
$(*)_{n,\epsilon}$ holds then for some $m \in \{1,2\}$
$(**)_{m,\epsilon}$ holds, where:

$(*)_{n,\epsilon}$ There is $p \in S^{bs}(N_{n,\epsilon})$ such
that $p$ is realized in $N^{\oplus}_{n,\epsilon}$ and
$f_{n,2\epsilon}(p)$ is realized in $N_{2\epsilon}$.

$(**)_{m,\epsilon},\ f_{m,2\epsilon+1}[N^{\oplus}_{m,\epsilon}]
\bigcap N_{2\epsilon} \neq f_{m,2\epsilon}[N_{m,\epsilon}]$.
\end{enumerate}
Note that requirement 4 is essentially a property of
$N_{2\epsilon+2}$ and $(**)_{m,\epsilon}$ is essentially a
property of $f_{m,2\epsilon+1}$.

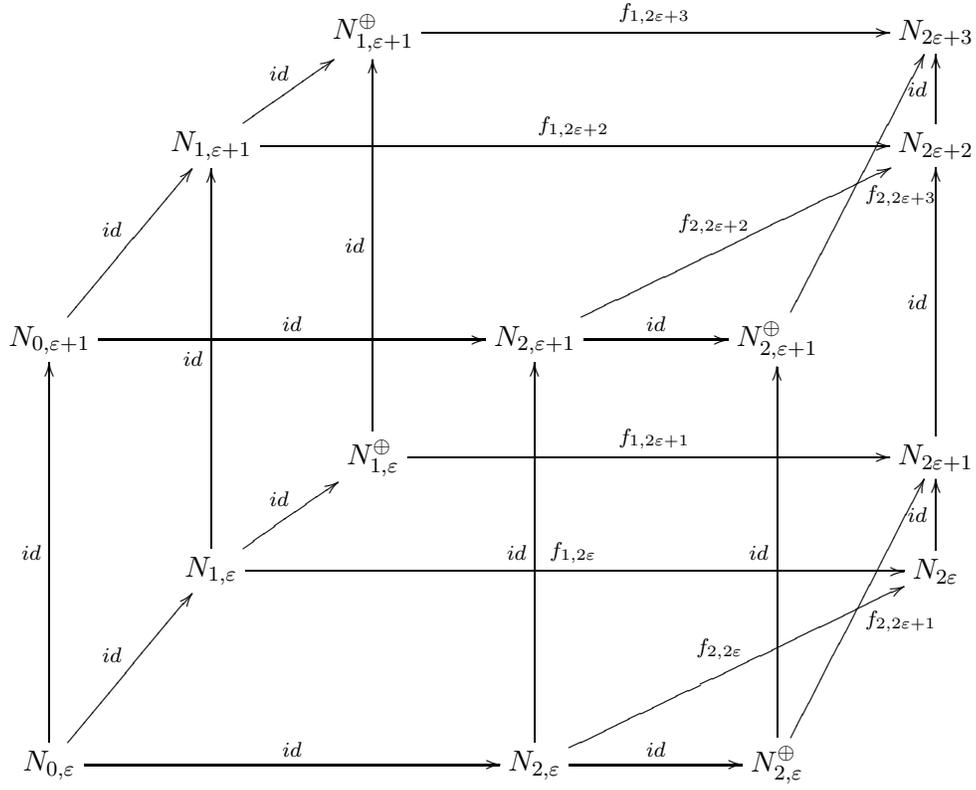
\begin{displaymath}
\xymatrix{&&N^{\oplus}_{1,\epsilon+1}
\ar[rrrr]^{f_{1,2\epsilon+3}}
&&&& N_{2\epsilon+3}\\
&N_{1,\epsilon+1} \ar[ur]^{id} \ar[rrrrr]^{f_{1,2\epsilon+2}}
&&&&&N_{2\epsilon+2} \ar[u]^{id} \\
\\
N_{0,\epsilon+1} \ar[uur]^{id} \ar[rrr]^{id}&&&N_{2,\epsilon+1}
\ar[rr]^{id} \ar[uurrr]^{f_{2,2\epsilon+2}}&&
 N^{\oplus}_{2,\epsilon+1} \ar[uuur]_{f_{2,2\epsilon+3}}\\
&&N^{\oplus}_{1,\epsilon} \ar[uuuu]^{id}
\ar[rrrr]^{f_{1,2\epsilon+1}}
&&&& N_{2\epsilon+1} \ar[uuu]^{id}\\
&N_{1,\epsilon} \ar[uuuu]^{id} \ar[ur]^{id} \ar[rrrrr]^{f_{1,2\epsilon}} &&&&& N_{2\epsilon} \ar[u]^{id}\\
\\
N_{0,\epsilon} \ar[uuuu]^{id} \ar[uur]^{id}
\ar[rrr]^{id}&&&N_{2,\epsilon} \ar[uuuu]^{id} \ar[rr]^{id}
\ar[ururr]^{f_{2,2\epsilon}} && N^{\oplus}_{2,\epsilon}
\ar[uuuu]^{id} \ar[uuur]_{f_{2,2\epsilon+1}}}
\end{displaymath}

\newpage
\case{Why can we carry out the construction?}
For $\epsilon=0$ let $(f_{1,0},f_{2,0},N_0)$ be an amalgamation of
$N_{1,0},N_{2,0}$ over $N_{0,0}$, such that $N_0 \bigcap
M_0=N_{0,0}$ (i.e. we choose new elements for $N_0-N_{0,0}$). In
the proof of item b, by the definition of the equality between
types without loss of generality $f_{1,0}(a_1)=f_{2,0}(a_2)$, so 3
is satisfied. In the proof of item c, by Theorem \ref{the weak
uniqueness theorem of NF} (the weak uniqueness theorem of NF),
there is a joint embedding $f_{1,0},f_{2,0},N_0$ of
$N_{1,0},N_{2,0}$ over $N_{0,0} \bigcup N^*$. So 2 is satisfied.

For \case{limit} $\epsilon$ define $N_\epsilon=\bigcup
\{N_\zeta:\zeta<\epsilon \},\ f_{n,\epsilon}=\bigcup
\{f_{n,\zeta}:\zeta<\epsilon \}$. 5 is satisfied. 1 is satisfied
by axiom \ref{1.1}.1.c. 6 is satisfied by the continuity of the
sequence $\langle N_{n,\epsilon}:\epsilon<\lambda^+ \rangle$, and
by the smoothness (Definition \ref{1.1}.1.d). Clearly 7 is
satisfied. 4,8 are not relevant in the limit case.

\case{The successor case:} How can we construct
$N_{2\epsilon+1},f_{1,2\epsilon+1},f_{2,2\epsilon+1}$ and
$N_{2\epsilon+2},\allowbreak f_{1,2\epsilon+2},f_{2,2\epsilon+2}$,
assuming we have constructed
$N_{2\epsilon},f_{1,2\epsilon},f_{2,2\epsilon}$?

\case{The construction of
$N_{2\epsilon+1},f_{1,2\epsilon+1},f_{2,2\epsilon+1}$:} Without
loss of generality for some $n \in {1,2}$, we have
$(*)_{n,\epsilon}$ [Otherwise requirement 8 is not relevant and we
can use the existence of an amalgamation in
$(K_\lambda,\preceq)$]. Fix $n^*$ with $(*)_{n^*,\epsilon}$. We
are going to find
$N_{2\epsilon+1},f_{n^*,2\epsilon+1},f_{3-n^*,2\epsilon+1}$ with
$(**)_{n^*,\epsilon}$, namely
$f_{n^*,2\epsilon+1}[N^{\oplus}_{n^*,\epsilon}] \bigcap
N_{2\epsilon} \neq f_{n^*,2\epsilon}[N_{n^*,\epsilon}]$. Let $p$
be a witness for $(*)_{n^*,\epsilon}$, so for some $a,b$
$tp(a,N_{n^*,\epsilon},N_{n^*,\epsilon}^{\oplus})=p,\
tp(b,f_{n^*,2\epsilon}[N_{n^*,\epsilon}],N_{2\epsilon})=f_{n^*,2\epsilon}(p)$.
So
$tp(f_{n^*,2\epsilon}(a),f_{n^*,2\epsilon}[N_{n^*,\epsilon}],f_{n^*,2\epsilon}[N_{n^*,\epsilon}^{\oplus}])=tp(b,f_{n^*,2\epsilon}[N_{n^*,\epsilon}],N_{2\epsilon})$.
Hence by the definition of equality of types, for some
$N_{2\epsilon+1}^{temp},\allowbreak f_{n^*,2\epsilon+1}^{temp}$
the following hold:
\begin{enumerate}
\item $N_{2\epsilon} \preceq N_{2\epsilon+1}^{temp}$. \item
$f_{n^*,2\epsilon+1}^{temp}:N_{n^*,\epsilon}^{\oplus} \to
N_{2\epsilon+1}^{temp}$ is an embedding. \item $f_{n^*,2\epsilon}
\subseteq f_{n^*,2\epsilon+1}^{temp}$  \item
$f_{n^*,2\epsilon+1}^{temp}(a)=b$.
\end{enumerate}

\begin{displaymath}
\xymatrix{& N_{3-n^*,\epsilon}^{\oplus}
\ar[rr]^{f_{3-n^*,2\epsilon+1}} && N_{2\epsilon+1} \\
&& a \in N_{n^*,\epsilon}^{\oplus} \ar[ru]^{f_{n^*,2\epsilon+1}}
\ar[r]_{f_{n^*,2\epsilon+1}^{temp}} &
N_{2\epsilon+1}^{temp} \ar[u]^{id}\\
& N_{3-n^*,\epsilon} \ar[ru]^{id} \ar[rr]^{f_{3-n^*,2\epsilon}}
\ar[uu]^{id} && N_{2\epsilon} \ni b \ar[u]^{id}\\
N_{0,\epsilon} \ar[ru]^{id} \ar[rr]^{id} && N_{n^*,\epsilon}
\ar[ru]^{f_{n^*,2\epsilon}} \ar[uu]^{id} }
\end{displaymath}

\begin{claim}\label{**_{n^*}}
$f_{n^*,2\epsilon+1}^{temp}[N^{\oplus}_{n^*,\epsilon}] \bigcap
N_{2\epsilon} \neq f_{n^*,2\epsilon}^{temp}[N_{n^*,\epsilon}]$.
\end{claim}

\begin{proof}
$b \in N_{2\epsilon}$. $p$ is a basic type so a non-algebraic one.
So $a \in N_{n^*,\epsilon}^{\oplus}-N_{n^*,\epsilon}$. Hence
$b=f_{n^*,2\epsilon+1}^{temp}(a) \in
f_{n^*,2\epsilon+1}^{temp}[N_{n^*,\epsilon}^{\oplus}]-f_{n^*,2\epsilon+1}^{temp}[N_{n^*,\epsilon}]$.
Therefore $b \in
f_{n^*,2\epsilon+1}^{temp}[N^{\oplus}_{n^*,\epsilon}] \bigcap
N_{2\epsilon}-f_{n^*,2\epsilon}^{temp}[N_{n^*,\epsilon}]$.
\end{proof}

 As
$(K_\lambda,\preceq)$ satisfies amalgamation, there are
$N_{2\epsilon+1},f_{3-n^*,2\epsilon+1}$ such that
$N_{2\epsilon+1}^{temp} \preceq N_{2\epsilon+1}$ and
$f_{3-n^*,2\epsilon+1}:N_{3-n^*,\epsilon}^{\oplus} \to
N_{2\epsilon+1}$ is an embedding that includes
$f_{3-n^*,2\epsilon}$. Now we define
$f_{n^*,2\epsilon+1}:N_{n^*,\epsilon}^{\oplus} \to
N_{2\epsilon+1}$ by
$f_{n^*,2\epsilon+1}(x)=f_{n^*,2\epsilon+1}^{temp}(x)$. By Claim
\ref{**_{n^*}} $(**)_{n^*}$ holds, so requirement 8 is satisfied.
As for $m=1,2$ the embedding $f_{m,2\epsilon+1}$ includes
$f_{m,2\epsilon}$, requirement 7 is satisfied. Without loss of
generality requirement 1 is satisfied. Requirement 4 is not
relevant in this case. Requirements 5,6 are satisfied.

\case{The construction of $N_{2\epsilon+2},f_{n,2\epsilon+2}$:} By
Proposition \ref{7.2}, there are $N_{2\epsilon+2},\allowbreak
f_{1,\allowbreak 2\epsilon+2},\allowbreak f_{2,2\epsilon+2}$ such
that: $NF(f_{n,2\epsilon+1}[N_{n,\epsilon}^{\oplus}],
f_{n,2\epsilon+2}[N_{n,\epsilon+1}],N_{2\epsilon+1},N_{2\epsilon+2})$,
and the reduction of $f_{n,2\epsilon+1}$ to $N_{0,\epsilon}$ is
the identity [Let $f_{n,2\epsilon+1}^+$ be an injection of
$N_{n,\epsilon+1},\ f_{n,2\epsilon+1} \subseteq
f_{n,2\epsilon+1}^+$, and the reduction of $f_{n,2\epsilon+1}^+$
to $N_{0,\epsilon+1}$ is the identity. Substitute the models
$N_{0,\epsilon},N_{0,\epsilon+1},f_{n,2\epsilon+1}[N_{n,\epsilon}^{\oplus}],
N_{2\epsilon+1},f_{2\epsilon+1}^+ \allowbreak
[N_{n,\epsilon+1}],N_{2\epsilon+2}$ which appear here, instead of
the models $M_{0,0},M_{0,1},M_{n,0},N_0,\allowbreak M_{n,1},N_1$
which appear in Proposition \ref{7.2} respectively. Assumption a
of Proposition \ref{7.2} (i.e. $NF(N_{0,\epsilon},\allowbreak
N_{0,\epsilon+1},f_{n,2\epsilon+1}[N_{n,\epsilon}^{\oplus}],
f_{n,2\epsilon+1}^+[N_{n,\epsilon+1}])$), is satisfied by
Definition \ref{7.3}.a (remember that $f_{n,2\epsilon+1}^+$ is an
isomorphism over $N_{0,\epsilon+1}$ and $NF$ respects
isomorphisms). Assumption b of Proposition \ref{7.2} is satisfied
by requirement 6 of the induction hypothesis. Assumption c of
Proposition \ref{7.2} is satisfied by requirement 4 of the
induction hypothesis]. Hence we can carry out the construction.

\case{Why is it sufficient?} By clause 7 for $n=1,2$
$f_{n,\lambda^+}:M_n \to N_{\lambda^+}$ is an embedding over
$M_0$.

\begin{claim}\label{sufficient 7.5.new}
$f_{1,\lambda^+}[M_1]=f_{2,\lambda^+}[M_2]=N_{\lambda^+}$.
\end{claim}

\begin{proof}
Toward a contradiction suppose there is $n \in \{1,2\}$ such that
$f_{n,\lambda^+}[M_n] \allowbreak \neq N_{\lambda^+}$. By Density
(Theorem \ref{2.10}.1), there is an element $b$ such that
$tp(b,\allowbreak f_{n,\lambda^+}[M_n],N_{\lambda^+})$ is basic.
$\langle f_{n,2\epsilon}[N_{n,\epsilon}]:\epsilon<\lambda^+
\rangle$ is a representation of $f_{n,\lambda^+}[M_n]$, so by
Definition \ref{2.9} there is $\epsilon<\lambda^+$ such that for
every $\zeta \in (\epsilon,\lambda^+)$ the type
$q_\zeta:=tp(b,f_{n,2\zeta}[N_{n,\zeta}],N_{\lambda^+})$ does not
fork over $f_{n,2\epsilon}[N_{n,\epsilon}]$. We choose this
$\epsilon$ such that $b \in N_{2\epsilon}$, (remember: $b \in
N_{\lambda^+}=\bigcup \{N_\epsilon:\epsilon<\lambda^+\}$). So
$q_\zeta$ is basic. Define $p_\zeta:=f_{n,2\zeta}^{-1}(q_\zeta)$.
So $p_\epsilon \in S^{bs}(N_{n,\epsilon})$. For every $\zeta \in
(\epsilon,\lambda^+),\ q_\zeta$ is the non-forking extension of
$q_\epsilon$, so $p_\zeta$ is the non-forking extension of
$p_\epsilon$. Hence by Definition \ref{7.3}, there is an end
segment $S^* \subseteq \lambda^+$ such that for $\zeta \in S^*$,
$p_\zeta$ is realized in $N^{\oplus}_{2\zeta}$. But
$q_\zeta=tp(b,f_{n,2\zeta}[N_{n,\zeta}],N_{2\zeta})$. So for every
$\zeta \in S^*$ we have $(*)_{n,\zeta}$ ($p_\zeta$ is a witness
for this). So by clause 8 there are $m \in \{1,2 \}$ and a
stationary set $S^{**} \subseteq S^*$ such that for every $\zeta
\in S^{**}$ we have $(**)_{m,\zeta}$, (there are no two
non-stationary subsets which their union is an end segment of
$\lambda^+$). The sequences $\langle N_{2\zeta}:\zeta \in S^{**}
\rangle,\ \langle N_{m,\zeta}:\zeta \in S^{**} \rangle$, $\langle
f_{m,2\zeta}:\zeta \in S^{**} \rangle$ are increasing and
continuous. But by $(**)_{m,\zeta}$, we have
$f_{m,2\zeta+1}[N^{\oplus}_{m,\zeta+1}] \bigcap N_{2\zeta} \neq
f_{m,2\zeta}[N_{m,\zeta}]$, in contradiction to Proposition
\ref{1.11}.
\end{proof}
By Claim \ref{sufficient 7.5.new} $f_{2,\lambda^+}^{-1} \circ
f_{1,\lambda^+}$ is an embedding of $M_1$ onto $M_2$ over $M_0$.
In the proof of item b we have to note that $f_{2,\lambda^+}^{-1}
\circ f_{1,\lambda^+}(a_1)=f_{2,0}^{-1} \circ f_{1,0}(a_1)=a_2$
(by clause 3). In the proof of item c we have to note that
$f_{2,\lambda^+}^{-1} \circ f_{1,\lambda^+} \restriction
N^*=f_{2,0}^{-1} \circ f_{1,0} \restriction N^*$ and by clause 3
it is the identity.
\end{proof}

\newpage
\begin{corollary} \label{7.7}
\mbox{}
\begin{enumerate}[(a)]
\item $(K_{\lambda^+},\preceq^{NF}\restriction K_{\lambda^+})$ has
amalgamation. So $(K^{sat},\preceq^{NF}\restriction K^{sat})$ has
amalgamation. \item Locality: Let $M_0,M_1,M_2$ be models in
$K_{\lambda^+}$, such that $M_0 \preceq M_1,\ M_0 \preceq M_2$.
Suppose there is $N_0 \in K_\lambda$ such that: $N_0\prec M_0$ and
for every $N \in K_\lambda$, $[N_0 \preceq N \preceq M_0]
\Rightarrow tp(a_1,N,M_1)\allowbreak =tp(a_2,N,M_2)$. Then
$tp(a_1,M_0,M_1)=tp(a_2,M_0,M_2)$. [The version we actually use:
Suppose there is $N_0 \in K_\lambda$ such that $tp(a_n,M_0,M_2)$
does not fork over $N_0$ and $tp(a_1,N_0,M_1)=tp(a_2,N_0,M_2)$.
Then $tp(a_1,M_0,M_1)=tp(a_2,M_0,M_2)$].
\end{enumerate}
\end{corollary}

\begin{proof}
\mbox{}
\begin{enumerate}[(a)]
\item Suppose for $n=1,2$ $M_0\prec^{NF}M_n$. By Proposition
\ref{7.4}.a, there is $M_n^+$ such that $M_n\prec^+M_n^+$. By
Proposition \ref{7.4}.a $M_0\prec^+M_n^+$. So by Theorem
\ref{7.13}.c (the uniqueness of the $\prec^+$-extension), there is
an isomorphism $f:M_1^+ \to M_2^+$ over $M_0$. Hence
$(f\restriction M_1,id_{M_2},M_2^+)$ is an amalgamation of
$M_1,M_2$ over $M_0$. Now Proposition \ref{before 7.4}.a. \item
Locality: By Proposition \ref{7.4}.a there is $M_n^+$ such that
$M_n\prec^+M_n^+$. By Theorem \ref{7.13}.b there is an isomorphism
$f:M_1^+ \to M_2^+$ over $M_0$, such that $f(a_1)=a_2$. So
$(f\restriction M_1,id_{M_2},M_2^+)$ witnesses that
$tp(a_1,M_0,M_1)=tp(a_2,M_0,M_2)$.
\end{enumerate}
\end{proof}

\begin{theorem} \label{7.8}
\begin{enumerate}[(a)]
\item $(K^{sat},\preceq^{NF} \restriction K^{sat})$ satisfies
axiom c in $\lambda^+$ (\ref{1.1}.2.c). \item If
$(K^{sat},\preceq^{NF} \restriction K^{sat})$ satisfies
smoothness, then it is an a.e.c. in $\lambda^+$. \item
$(K^{sat},\preceq^{NF} \restriction K^{sat})$ has the amalgamation
property.
\end{enumerate}
\end{theorem}

\begin{proof}
\mbox{} (a) Let $j<\lambda^{+2}$ and $\langle M_i:i<j \rangle$ be
an $\preceq^{NF}$-increasing continuous of models in $K^{sat}$.
Let $M_j$ be the union of this sequence. We prove that $M_j \in
K^{sat}$ by induction on $j$. Let $N$ be a model in $K_\lambda$
such that $N \prec M_j$.

\case{Case a:} $\lambda<cf(j)$. In this case for some $i<j$
$N\prec M_i$. Since $M_i$ is full over $N$, of course $M_j$ is.
Therefore $M_j \in K^{sat}$.

\case{Case b:} $cf(j) \leq \lambda$. Without loss of generality
$cf(j)=j$. So $|j|=j=cf(j) \leq \lambda$. Let $\langle
N_{i,\alpha}:\alpha \in \lambda^+ \rangle$ be a representation of
$M_i$. For every $i<j$ let $E_i$ be a club of $\lambda^+$ such
that for $\alpha \in E_i,\ NF(N_{i,\alpha},\allowbreak
N_{i+1,\alpha},\allowbreak N_{i,\alpha+1},N_{i+1,\alpha+1})$ and
if $i$ is a limit ordinal, then $N_{i,\alpha}=\bigcup
\{N_{\epsilon,\alpha}:\epsilon<i \}$. So $E:=\bigcap \{E_i:i<j \}$
is a club set of $\lambda^+$ (because $|j| \leq \lambda$). Define
$N_{j,\alpha}:=\bigcup \{N_{i,\alpha}:i<j\}$. $\langle
N_{j,\alpha}:\alpha \leq \lambda^+ \rangle$ is a representation of
$M_j$. Take $\alpha^* \in E$ such that $N \subseteq
N_{j,\alpha^*}$. By axiom \ref{1.1}.1.e $N \preceq
N_{j,\alpha^*}$, so it is sufficient to prove that $M_j$ is
saturated over $N_{j,\alpha^*}$. Let $q \in
S^{bs}(N_{j,\alpha^*})$. We will prove that $q$ is realized in
$M_j$. By the definition of $E$ the sequence $\langle
N_{i,\alpha^*}:i<j \rangle$ is increasing and continuous, so by
Definition \ref{2.1a}.3.c (the local character) there is an
ordinal $i<j$ such that $q$ does not fork over $N_{i,\alpha^*}$.
$M_i$ is saturated in $\lambda^+$ over $\lambda$, so there is $a
\in M_i$ such that $tp(a,N_{i,\alpha^*},M_i)=q\restriction
N_{i,\alpha^*}$. By Definition \ref{5.14} we have
$\widehat{NF}(N_{i,\alpha^*},N_{j,\alpha^*},M_i,M_j)$, so by
Theorem \ref{5.15}.e ($\widehat{NF}$ respects $\frak{s}$)
$tp(a,N_{j,\alpha^*},M_j)$ does not fork over $N_{i,\alpha^*}$.
Hence by Definition \ref{2.1a}.3.d (the uniqueness of the
non-forking extension) $tp(a,N_{j,\alpha^*},M_j)=q$.

(b) The first part of Axiom c of a.e.c. in $\lambda^+$ is item a
here. Axioms b,e and the second part of axiom c follows by
Proposition \ref{6.2}.f.

(c) By Corollary \ref{7.7}.a.
\end{proof}

\section{relative saturation} \discussion{Discussion:} This
section is, like its previous, a preparation for the proof of
Theorem \ref{9.6}. We study the relation $\preceq^{\otimes}$, a
kind of relative saturation. This relation is similar to `closure
of $\preceq^{NF}$ under smoothness' (see Proposition \ref{8.2}.b).
Theorem \ref{9.5} says that non-equality between the relations
$\preceq^{NF},\preceq^{\otimes}$ is equivalent to non-smoothness
and also to a strengthened version of non-smoothness.

\begin{hypothesis}
$\frak{s}$ is a weakly successful semi-good $\lambda$-frame with
conjugation.
\end{hypothesis}

\begin{definition} \label{8.1}
$\preceq^{\otimes}:=\{(M_0,M_1):M_0,M_1 \in K^{sat}$, $M_0 \prec
M_1$ and for every $N_0,N_1 \in K_\lambda$, if $M_0 \succeq N_0
\preceq N_1 \preceq M_1$ and $p \in S^{bs}(N_1)$ does not fork
over $N_0$, then for some element $d \in M_0$ $tp(d,N_1,M_1)=p$\}.
\end{definition}

\begin{proposition} \label{8.2}
\mbox{}
\begin{enumerate}[(a)]
\item If $M_0 \in K^{sat}$ and $M_0 \preceq^{NF} M_1$ then $M_0
\preceq^{\otimes} M_1$. \item If $\langle M_\epsilon:\epsilon \leq
\delta \rangle$ is an $\preceq^{NF}$-increasing continuous
sequence of models in $K^{sat}$ and for every $\epsilon \in
\delta,\ M_\epsilon \preceq^{NF}M_{\delta+1}$, then $M_\delta
\preceq^{\otimes} M_{\delta+1}$.
\end{enumerate}
\end{proposition}

\begin{proof}
\mbox{} (a) Suppose $M_0 \preceq^{NF} M_1$ and $M_0 \in K^{sat}$.
Let $N_0,N_1$ be models $K_\lambda$ with $M_0 \succeq N_0 \preceq
N_1 \preceq M_1$ and let $p$ be a type $S^{bs}(N_1)$ that does not
fork over $N_0$. We have to find an element $d \in M_0$ with
$tp(d,N_1,M_1)=p$.
By Proposition \ref{6.2}.g (LST for $\preceq^{NF}$) for some
$N_0^+,N_1^+ \in K_\lambda$ $N_0 \preceq N_0^+$, $N_1 \preceq
N_1^+$ and $\widehat{NF}(N_0^+,N_1^+,M_0,M_1)$. By axiom
\ref{1.1}.1.e $N_0 \preceq N_{0}^+$ and $N_1 \preceq N_1^+$. Let
$q$ be the non-forking extension of $p$ to $N_1^+$. Since $M_0 \in
K^{sat}$ for some $d \in M_0$ $tp(d,N_0^+,M_0)=q \restriction
N_0^+$. By Proposition \ref{s is transitive} $q$ does not fork
over $N_0$, so by Definition \ref{2.1a}.3.b (monotonicity) $q$
does not fork over $N_0^+$.
By Theorem \ref{5.15} $\widehat{NF}$ respects $\frak{s}$, so
$tp(d,N_1^+,M_1)$ does not fork over $N_1^+$. So by Definition
\ref{2.1a}.3.b (uniqueness) $tp(d,N_1^+,M_1)=q$. Therefore
$tp(d,N_1,M_1)=p$.

(b) Suppose $N_0,N_1 \in K_\lambda$, $M_\delta \succeq N_0 \preceq
N_1 \preceq M_{\delta+1}$ and $p \in S^{bs}(N_1)$ does not fork
over $N_0$. We have to find an element $d \in M_\delta$ that
realizes $p$. For every $\alpha \leq \delta+1$ there is a
representation $\langle N_{\alpha,\epsilon}:\epsilon<\lambda^+
\rangle$ of $M_\alpha$. without loss of generality
$cf(\delta)=\delta$.

\case{Case a:} $\delta=\lambda^+$. So for some $\alpha<\delta,\
N_0 \subseteq M_\alpha$ and we can use item a.

\case{Case b:} $\delta<\lambda^+$. For each $\alpha \in \delta$,
let $E_\alpha$  be a club of $\lambda^+$ such that for each
$\epsilon \in E_\alpha$:
$NF(N_{\alpha,\epsilon},N_{\alpha+1,\epsilon},\allowbreak
N_{\alpha,\epsilon+1},N_{\alpha+1,\epsilon+1})$ and if $\alpha$ is
limit then $N_{\alpha,\epsilon}=\bigcup
\{N_{\beta,\epsilon}:\beta<\alpha\}$. Let $E_\delta:=\{\alpha \in
\lambda^+:N_{\delta,\epsilon} \subseteq N_{\delta+1,\epsilon},\
N_{\delta,\epsilon}=\bigcup \{ N_{\alpha,\epsilon}:\alpha<\delta
\} \}$. Denote $E:=\bigcap \{E_\alpha:\alpha \leq \delta \}$. By
cardinality considerations there is $\epsilon \in E$ such that for
$n<2$ $N_n \subseteq N_{\delta+n,\epsilon}$, so by axiom
\ref{1.1}.1.e $N_n \preceq N_{\delta+n,\epsilon}$.
\begin{displaymath}
\xymatrix{d \in M_\alpha \ar[r]^{id} &M_\delta \ar[r]^{id} & M_{\delta+1}\\
N_{\alpha,\epsilon} \ar[r]^{id} \ar[u]^{id} & N_{\delta,\epsilon}
\ar[r]^{id} \ar[u]^{id} &
N_{\delta+1,\epsilon} \ar[u]^{id} & q\\
& N_0 \ar[r]^{id} \ar[u]^{id} & N_1 \ar[u]^{id} & p  }
\end{displaymath} Let $q \in
S^{bs}(N_{\delta+1,\epsilon})$ be the non-forking extension of
$p$. By Proposition \ref{s is transitive} (the transitivity
proposition), $q$ does not fork over $N_0$. By Definition
\ref{2.1a}.3.b (monotonicity) $q$ does not fork over
$N_{\delta,\epsilon}$, so $q\restriction N_{\delta,\epsilon}$ is
basic. As $\epsilon \in E$, the sequence $\langle
N_{\alpha,\epsilon}:\alpha\leq \delta \rangle$ is increasing and
continuous. So by Definition \ref{2.1a}.3.c (local character)
there is $\alpha<\delta$ such that $q\restriction
N_{\delta,\epsilon}$ does not fork over $N_{\alpha,\epsilon}$. So
by Proposition \ref{s is transitive} $q$ does not fork over
$N_{\alpha,\epsilon}$. Since $M_\alpha \preceq^{NF} M_{\delta+1}$
by item a for some $d \in M_\alpha$
$tp(d,N_{\delta+1,\epsilon},M_{\delta+1})=q$. So
$tp(d,N_1,M_{\delta+1})=p$.
\end{proof}

The following proposition is similar to the saturativity = model
homogeneity lemma.
\begin{proposition} \label{8.3}
Suppose
\begin{enumerate}
\item $M_0 \preceq^{\otimes}M_1$. \item For $n<3$ $N_n \in
K_\lambda$. \item $N_0 \preceq M_0$. \item $N_2 \succeq N_0
\preceq N_1 \preceq M_1$.
\end{enumerate}

Then for some $N_1^* \in K_\lambda$ and an embedding $f:N_2 \to
M_0$ the following hold:
\begin{enumerate}[(a)]
\item $f\restriction N_0=id_{N_0}$. \item
$NF(N_0,f[N_2],N_1,N_1^*)$. \item $N_1^* \preceq M_1$.
\end{enumerate}
\end{proposition}

\begin{displaymath}
\xymatrix{M_0 \ar[r]^{id} & M_1\\
f[N_2] \ar[r]^{id} \ar[u]^{id} & N_1^* \ar[u]^{id}\\
N_0 \ar[r]^{id} \ar[u]^{id} & N_1 \ar[u]^{id}  }
\end{displaymath}

\begin{proof}
We try to choose
$N_{0,\epsilon},N_{1,\epsilon},N_{2,\epsilon},f_\epsilon$ by
induction on $\epsilon<\lambda^+$ such that:
\begin{enumerate}[(1)]
\item For $n<3$ $\langle N_{n,\epsilon}:\epsilon<\lambda^+
\rangle$ is an increasing continuous of models in $K_\lambda$.
\item For $n<3$ $N_{n,0}=N_n,\ f_0=id_{N_0}$. \item For
$\epsilon<\lambda^+,\ N_{0,\epsilon} \preceq M_0 \wedge
N_{1,\epsilon} \preceq M_1$. \item $\langle
f_\epsilon:\epsilon<\lambda^+ \rangle$ is increasing and
continuous. \item $f_\epsilon:N_{0,\epsilon} \to N_{2,\epsilon}$
is an embedding over $N_0$. \item For every $\epsilon \in
\lambda^+$ there is $a_\epsilon$ such that
$(N_{0,\epsilon},N_{0,\epsilon+1},a_\epsilon)$ is a uniqueness
triple, $f_{\epsilon+1}(a_\epsilon) \in N_{2,\epsilon}$ and
$tp(a_\epsilon,N_{1,\epsilon},N_{1,\epsilon+1})$ does not fork
over $N_{0,\epsilon}$. \item $N_{0,\epsilon} \preceq
N_{1,\epsilon}$ (actually follows by 6).
\end{enumerate}

\begin{displaymath}
\xymatrix{&M_0 \ar[r]^{id} & M_1\\
N_{2,\epsilon+1} & N_{0,\epsilon+1} \ar[l]_{f_{\epsilon+1}}
\ar[r]^{id} \ar[u]^{id} &
N_{1,\epsilon+1} \ar[u]^{id}\\
N_{2,\epsilon} \ar[u]^{id} & N_{0,\epsilon} \ar[l]_{f_\epsilon} \ar[r]^{id} \ar[u]^{id} & N_{1,\epsilon} \ar[u]^{id}\\
& N_0 \ar[ul]^{id} \ar[u]^{id} \ar[r]^{id} & N_1 \ar[u]^{id} }
\end{displaymath}

By clauses 1,4,5 and particularly 6
and Proposition \ref{1.11} we cannot succeed. Where will we get
stuck? For $\epsilon=0$ or limit we will not get stuck. Suppose we
have defined
$N_{0,\epsilon},N_{1,\epsilon},N_{2,\epsilon},f_\epsilon$. By
clause 5, $f_\epsilon[N_{0,\epsilon}] \preceq N_{2,\epsilon}$.

\case{Case a:} $f_\epsilon[N_{0,\epsilon}] \neq N_{2,\epsilon}$.
In this case we can find
$N_{0,\epsilon+1},N_{1,\epsilon+1},N_{2,\epsilon+1},\allowbreak
f_{\epsilon+1}$ such that clauses 1-7 above hold [By the existence
of the basic types, there is $b \in
N_{2,\epsilon}-f_\epsilon[N_{0,\epsilon}]$ such that
$p:=tp(b,f_\epsilon[N_{0,\epsilon}],N_{2,\epsilon})$ is basic. Let
$q \in S^{bs}(N_{1,\epsilon})$ be the non-forking extension of
$f_\epsilon^{-1}(p)$. As $M_0 \preceq^{\otimes}M_1 \wedge (n<2
\Rightarrow N_{n,\epsilon} \preceq M_n^*) \wedge N_{0,\epsilon}
\preceq N_{1,\epsilon} \in K_\lambda$, there is $a \in M_0$ which
realizes $q$. So $tp(a,N_{0,\epsilon},M_0)=f_\epsilon^{-1}(p)$. As
$\frak{s}$ is weakly successful, we can find $N_{0,\epsilon+1}$
such that $(N_{0,\epsilon},N_{0,\epsilon+1},a) \in K^{3,uq}$. As
$M_0$ is saturated in $\lambda^+$ over $\lambda$, by Lemma
\ref{1.12} (the saturation = model homogeneity lemma), without
loss of generality $N_{0,\epsilon+1} \preceq M_0$. Denote $a$ as
$a_\epsilon$. Choose $N_{1,\epsilon+1} \preceq M_1$ such that
$N_{0,\epsilon+1} \bigcup N_{1,\epsilon} \subseteq
N_{1,\epsilon+1}$. By axiom \ref{1.1}.1.e $N_{0,\epsilon+1}
\preceq N_{1,\epsilon+1} \wedge N_{1,\epsilon} \preceq
N_{1,\epsilon+1}$. Now
$f_\epsilon(tp(a_\epsilon,N_{0,\epsilon},N_{0,\epsilon+1})=p$. So
there are $N_{2,\epsilon+1},f_{\epsilon+1}$ such that:
$N_{2,\epsilon} \preceq N_{2,\epsilon+1},\
f_{\epsilon+1}(a_\epsilon)=b,\ f_\epsilon \subseteq
f_{\epsilon+1}:N_{0,\epsilon+1} \to N_{2,\epsilon+1}$].

\case{Case b:} $f_\epsilon[N_{0,\epsilon}]=N_{2,\epsilon}$. Hence
$N_{1,\epsilon},f_\epsilon^{-1}\restriction N_2$ witness that our
proposition is true [By 6, Definition \ref{5.3} and Definition
\ref{5.2}, $\zeta<\epsilon \Rightarrow
NF(N_{0,\zeta},N_{0,\zeta+1},N_{1,\zeta},\allowbreak
N_{1,\zeta+1})$. So by Theorem \ref{5.12} (the long transitivity
theorem), $NF(N_0,N_{0,\epsilon},\allowbreak N_1,N_{1,\epsilon})$.
So by the monotonicity of NF, we have
$NF(N_0,f_\epsilon^{-1}[N_2],\allowbreak N_1,N_{1,\epsilon})$. So
clause b in the proposition is satisfied. Clauses a,c are
satisfied by 5,3 respectively].

Let $\epsilon+1$ be the first ordinal we will get stuck on . In
other words, suppose we have defined
$N_{0,\epsilon},N_{1,\epsilon},N_{2,\epsilon},f_\epsilon$ and we
cannot find
$N_{0,\epsilon+1},N_{1,\epsilon+1},N_{2,\epsilon+1},\allowbreak
f_{\epsilon+1}$ such that clauses 1-7 above hold. so case b holds
and the proposition is proved.
\end{proof}

\begin{proposition}\label{LST for pairs}
If $M_0 \preceq M_1,\ n<2 \Rightarrow (||M_n||)=\lambda^+ \wedge
A_n \subseteq M_n \wedge |A_n| \leq \lambda)$, then there are
models $N_0,N_1 \in K_\lambda$ such that: $n<2 \Rightarrow A_n
\subseteq N_n \preceq M_n$ and $N_1 \bigcap M_0=N_0$ (so of course
$N_0 \preceq N_1$).
\end{proposition}

\begin{proof}
For $n<2$ we will construct by induction on $m<\omega$ a model
$N_{n,m}$ such that $\langle N_{n,m}:m \leq \omega \rangle$ is
$\preceq$-increasing and continuous of models in $K_\lambda$, $A_n
\subseteq N_{n,0},\ N_{0,m} \subseteq N_{1,m},\ N_{1,m} \bigcap
M_0 \subseteq N_{0,m+1},\ N_{n,m} \preceq M_n$. This construction
is possible as $LST(K,\preceq) \leq \lambda$. Now $M_0 \bigcap
N_{1,\omega}=N_{0,\omega}$ [Why? If $x \in M_0 \bigcap
N_{1,\omega}$, then for some $m<\omega$ we have $x \in N_{1,m}
\bigcap M_0 \subseteq N_{0,m+1} \subseteq N_{0,\omega}$ and from
the other side, if $x \in N_{0,\omega}$ then for some $m<\omega$
we have $x \in N_{0,m} \subseteq N_{1,m}$, so $x \in M_0 \bigcap
N_{1,\omega}$].
\end{proof}

\begin{proposition} \label{8.4}
If $M_1^* \preceq^{\otimes} M_2^*$ then there is an increasing
continuous sequence of models in $K^{sat}$, $\langle
M_\epsilon:\epsilon \leq \lambda^++1 \rangle$ such that:
\begin{enumerate}[(a)]
\item $M_{\lambda^+}=M_1^*,\ M_{\lambda^++1}=M_2^*$. \item
$\epsilon<\lambda^+ \Rightarrow M_\epsilon\prec^+M_{\epsilon+1}$.
\item $\epsilon<\lambda^+ \Rightarrow M_\epsilon \preceq^{NF}
M_2^*$.
\end{enumerate}
\end{proposition}

\begin{proof}
By Proposition \ref{7.4}.c, there is a winning strategy for player
2 in the $\prec^+$-game. Let $F$ be such a winning strategy.
Enumerate $M_2^*$ by $\{a_\epsilon:\epsilon<\lambda^+\}$. We
construct $\langle N_{\alpha,\epsilon}:\epsilon \leq \alpha
\rangle,\ N_\alpha$ by induction on $\alpha$ such that the
following hold:
\begin{enumerate}
\item For each $\epsilon \leq \alpha$, $N_{\alpha,\epsilon} \in
K_\lambda$ and $N_{\alpha,\epsilon} \preceq M_1^*$. \item $\langle
N_{\alpha,\epsilon}:\epsilon \leq \alpha<\lambda^+ \rangle$ is
increasing continuous in the variables $\alpha,\epsilon$. \item
$\langle N_\alpha:\alpha<\lambda^+ \rangle$ is an increasing
continuous of models in $K_\lambda$. \item $N_{\alpha,\alpha}
\preceq N_\alpha \preceq M_2^*$. \item If $\alpha+1$ is odd then
for each $\epsilon \leq \alpha$, $N_{\alpha+1,\epsilon+1}$ is
isomorphic to $F(\langle N_{\beta,\epsilon}:\epsilon+1 \leq \beta
\leq \alpha+1 \rangle,\ \langle N_{\beta,\epsilon+1}:\epsilon+1
\leq \beta \leq \alpha \rangle)$ over $N_{\alpha,\epsilon+1}
\bigcup N_{\alpha+1,\epsilon}$. \item If $\alpha+1$ is odd then
$NF(N_{\alpha,\alpha},N_{\alpha},N_{\alpha+1,\alpha+1},N_{\alpha+1})$
\item $a_\alpha \in N_{2\alpha+2}$. \item $N_{2\alpha} \bigcap
M_1^* \subseteq N_{2\alpha,2\alpha}$. \item If $\alpha+1$ is odd
then $N_{\alpha+1,\alpha+1}=N_{\alpha+1,\alpha}$. \item If
$\alpha+1$ is odd then $N_{\alpha+1,0} \bigcap
N_{\alpha}=N_{\alpha,0},\ N_{\alpha+1,0} \neq N_{\alpha,0}$. \item
If $\alpha+1$ is even then for each $\epsilon \leq \alpha$
$N_{\alpha+1,\epsilon}=N_{\alpha,\epsilon}$.
\end{enumerate}

\begin{displaymath}
\xymatrix{ {M_\epsilon} \ar[r]^{id} & M_{\epsilon+1} \ar[r]^{id} &
M_\alpha \ar[r]^{id} & M_{\lambda^+}=M_1^* \ar[r]^{id} & M_{\lambda^++1}=M_2^*\\
N_{\alpha,\epsilon} \ar[r]^{id} \ar[u]^{id} &
N_{\alpha,\epsilon+1} \ar[r]^{id} \ar[u]^{id} & N_{\alpha,\alpha}
\ar[rr]^{id} \ar[u]^{id} &&
N_\alpha \ar[u]^{id} \\
N_{\epsilon+1,\epsilon} \ar[r]^{id} \ar[u]^{id} &
N_{\epsilon+1,\epsilon+1} \ar[rrr]^{id} \ar[u]^{id} &&&
N_{\epsilon+1} \ar[u]^{id}\\
N_{\epsilon,\epsilon} \ar[rrrr]^{id} \ar[u]^{id} &&&& N_\epsilon
\ar[u]^{id}  }
\end{displaymath}

[Explanation: $N_{\alpha,\alpha},\ N_\alpha$ are approximations
for $M_1^*,\ M_2^*$ respectively. $N_{\alpha,\epsilon}$ is an
approximation for $M_\epsilon$. When $\alpha+1$ is even, we
increase the approximations of $M_1^*,M_2^*$ such that in the end
we will have $M_2^* \subseteq \bigcup
\{N_\alpha:\alpha<\lambda^+\},\ M_1^*=\bigcup
\{N_{\alpha,\alpha}:\alpha<\lambda^+\}$ by 7,8 respectively. when
$\alpha+1$ is odd, we increase the approximations of $M_\epsilon$
(mainly by clause 10). Clause 11 says that in even step the
approximations to $M_\epsilon$ do not increase. Clause 5 insure
that in the end we will have $M_\epsilon\prec^+M_{\epsilon+1}$.
Clause 6 insure that in the end requirement c will be satisfied.
The point of the proof is, that we could not demand 6 for every
$\alpha$, (as otherwise we prove $M_1^* \preceq^{NF} M_2^*$, which
might be wrong). But we succeed to prove that
$NF(N_{\alpha,\epsilon},N_\alpha,N_{\alpha+1,\epsilon},N_{\alpha+1})$ so $M_\epsilon \preceq^{NF} M_2^*$].\\

\case{Why can we carry out the construction?} We construct by
induction on $\alpha$. For limit $\alpha$, by clauses 2,3 there is
no freedom. Clauses 1,4 are satisfied by the smoothness, clauses
5,6,7,9,10,11 are not relevant and clause 8 is satisfied. For
$\alpha=0$ we choose $N_0,N_{0,0}$ by Proposition \ref{LST for
pairs}. Suppose we have defined $\langle
N_{\alpha,\epsilon}:\epsilon \leq \alpha
\rangle,\ N_\alpha$. what will we do in step $\alpha+1$?\\
\case{Case a:} $\alpha+1$ is even. For $\epsilon \leq \alpha$
define $N_{\alpha+1,\epsilon}:=N_{\alpha,\epsilon}$. By
Proposition \ref{LST for pairs} there are $N_{\alpha+1},\
N_{\alpha+1,\alpha+1}$ as required, especially clauses 7,8 are
satisfied.\\ \case{Case b:} $\alpha+1$ is odd. Define
$N_{\alpha+1,\epsilon}^{temp}$ by induction on $\epsilon \leq
\alpha$ such that:
\begin{enumerate}
\item $\langle  N_{\alpha+1,\epsilon}^{temp}:\epsilon \leq \alpha
\rangle$ is an $\preceq$-increasing continuous sequence. \item
$N_{\alpha+1,\epsilon+1}^{temp}=F(\langle
N_{\beta,\epsilon}:\epsilon+1 \leq \beta \leq \alpha
\rangle^\frown \langle  N_{\alpha+1,\epsilon}^{temp} \rangle,\
\langle  N_{\beta,\epsilon+1}:\epsilon+1 \leq \beta<\alpha
\rangle)$. \item $N_{\alpha,0} \preceq N_{\alpha+1,0}^{temp}$.
\end{enumerate}
Now by Proposition \ref{8.3}, there are $N_{\alpha+1}$ and an
embedding $g:N_{\alpha+1,\alpha}^{temp} \to M_1^*$ over
$N_{\alpha,\alpha}$ such that we have
$NF(N_{\alpha,\alpha},N_\alpha,g[N_{\alpha+1,\alpha}^{temp}],N_{\alpha+1})$.
For every $\epsilon \leq \alpha$ define
$N_{\alpha+1,\epsilon}:=g[N_{\alpha+1,\epsilon}^{temp}]$. Now
define $N_{\alpha+1,\alpha+1}:=N_{\alpha+1,\alpha}$. So we can
carry out the construction.\\

\case{Why is it sufficient?} For $\epsilon<\lambda^+$ define
$M_\epsilon:=\bigcup \{N_{\alpha,\epsilon}:\epsilon \leq
\alpha<\lambda^+\}$. Define $M_{\lambda^+}:=\bigcup
\{M_\epsilon:\epsilon<\lambda^+\},\ M_{\lambda^++1}:=\bigcup
\{N_\alpha:\alpha<\lambda^+\}$. We will prove that the sequence
$\langle M_\epsilon:0<\epsilon<\lambda^++1 \rangle$ satisfies
requirements a,b,c:\\
(a) By 3,4,7 $M_{\lambda^++1}=M_2^*$. Why is
$M_{\lambda^+}=M_1^*$? By 1 $M_{\lambda^+} \subseteq M_1^*$. Let
$x \in M_1^*$. Then $x \in M_2^*=M_{\lambda^++1}$. So by the
definition of $M_{\lambda^++1}$ and 3, there is $\alpha$ such that
$x \in N_{2\alpha}$. So by 8 $x \in N_{2\alpha,2\alpha}$. But by
the definitions of $M_\epsilon,M_{\lambda^+},\ N_{2\alpha,2\alpha}
\subseteq M_{2\alpha}
\subseteq M_{\lambda^+}$.\\
(b) By 2,10 $|M_0|=\lambda^+$. By 2 and the smoothness, the
sequence $\langle M_\epsilon:\epsilon<\lambda^+ \rangle$ is
$\preceq$-increasing and continuous. So $|M_\epsilon|=\lambda^+$.
Does $\epsilon<\lambda^+ \Rightarrow M_\epsilon \in K^{sat}$? Not
exactly, but we can prove by induction on $\epsilon$ that
$0<\epsilon<\lambda^+ \Rightarrow (M_\epsilon \in K^{sat} \wedge
M_\epsilon\prec^+M_{\epsilon+1})$: For $\epsilon=0$ by 10. For
limit $\epsilon$ by Theorem \ref{7.8}.a. For $\epsilon$ successor
by 5
and Proposition \ref{before 7.4}.b. So requirement b is satisfied.\\
(c) The sequences $\langle  N_{\alpha,\epsilon}:\epsilon \leq
\alpha<\lambda^+ \rangle,\ \langle N_\alpha:\epsilon \leq
\alpha<\lambda^+ \rangle$ are representations of $M_\epsilon,\
M_{\lambda^++1}$ respectively. Let $\alpha \in \lambda^+$. We will
prove
$NF(N_{\alpha,\epsilon},N_\alpha,N_{\alpha+1,\epsilon},N_{\alpha+1})$.
If $\alpha+1$ is even, this is satisfied by clause 11. So let
$\alpha+1$ be odd. By 6 we have: (*) $NF(N_{\alpha,\alpha},N_\alpha,N_{\alpha+1,\alpha+1},N_{\alpha+1})$.\\
By 5 and Theorem \ref{5.12} (the transitivity of NF),
$NF(N_{\alpha,\epsilon},N_{\alpha,\alpha},N_{\alpha+1,\epsilon},\allowbreak
N_{\alpha+1,\alpha})$ [Why? By 5 (and Proposition \ref{7.4}.c),
$\forall \zeta \in [\epsilon,\alpha)
NF(N_{\alpha,\zeta},N_{\alpha,\zeta+1},\allowbreak
N_{\alpha+1,\zeta},\allowbreak N_{\alpha+1,\zeta+1})$. The
sequences $\langle N_{\alpha,\zeta}:\zeta \in [\epsilon,\alpha)
\rangle,\ \langle N_{\alpha+1,\zeta}:\zeta \in [\epsilon,\alpha)
\rangle$ are increasing and continuous. So by Theorem \ref{5.12}
(the long transitivity theorem),
$NF(N_{\alpha,\epsilon},N_{\alpha,\alpha},N_{\alpha+1,\epsilon},N_{\alpha+1,\alpha})$.
So by the monotonicity of NF, we have: (**)
$NF(N_{\alpha,\epsilon},N_{\alpha,\alpha},N_{\alpha+1,\epsilon},N_{\alpha+1,\alpha+1})$].
Now by (*),(**) and Theorem \ref{5.12}
$NF(N_{\alpha,\epsilon},N_{\alpha+1,\epsilon},N_\alpha,N_{\alpha+1})$.
Note that we use here freely Theorem \ref{5.9} (the symmetry
theorem of NF).
\end{proof}

\section{Non-smoothness implies non-structure}
\begin{hypothesis}
$\frak{s}$ is a weakly successful semi-good $\lambda$-frame with
conjugation.
\end{hypothesis}

\begin{definition} \label{9.1}
Let $\bar{M}=\langle  M_\alpha:\alpha<\alpha^* \rangle$ be an
increasing sequence of models in $K_{\lambda^+}$. We say that
$\bar{M}$ is \emph{$\preceq^{NF}$-increasing in the successor
ordinals} if $\beta<\gamma<\alpha^* \Rightarrow M_{\beta+1}
\preceq^{NF} M_{\gamma+1}$.
\end{definition}

\begin{definition} \label{9.2}
Let $\alpha \leq \lambda^{+2}$ and let $\bar{M}=\langle
M_\alpha:\alpha<\lambda^{+2} \rangle$ be an
$\preceq^{NF}$-increasing in the successor ordinals and continuous
sequence with union $M$. Define $S(\bar{M})=:\{\delta \in
\lambda^{+2}:\exists \alpha \in (\delta,\lambda^{+2})$ $M_\delta
\npreceq^{NF} M_\alpha\}$. Define
$S(M)=:S(\bar{M})/D_{\lambda^{+2}}$ where $D_{\lambda^{+2}}$ is
the clubs filter on $\lambda^{+2}$. (By Proposition \ref{is well
defined} $S(M)$ does not depend on the representation $\bar{M}$).
\end{definition}

\begin{proposition} \label{9.3}
Let $\bar{M}=\langle  M_\alpha:\alpha<\lambda^{+2} \rangle$ be a
$\preceq^{NF}$-increasing in the successor ordinals and continuous
sequence. Then:
\begin{enumerate}[(a)] \item  For each
$\alpha,\beta$ with $\alpha<\beta<\lambda^{+2}$, $M_\alpha
\preceq^{NF} M_{\alpha+1} \Leftrightarrow M_\alpha \preceq^{NF}
M_\beta$. \item $S(\bar{M})=\{\delta \in \lambda^{+2}:\forall
\alpha \in (\delta,\lambda^{+2})$ $M_\delta \npreceq^{NF}
M_\alpha\}$.
\end{enumerate}
\end{proposition}

\begin{proof}
\mbox{}
\begin{enumerate}[(a)]
\item Easy (by Proposition \ref{6.2}.c). \item By item a.
\end{enumerate}
\end{proof}

\begin{proposition}\label{is well defined}
Suppose:
\begin{enumerate}
\item The sequences $\bar{M}^1:=\langle
M_{\alpha,1}:\alpha<\lambda^{+2} \rangle,\ \bar{M}^2:=\langle
M_{\alpha,1}:\alpha<\lambda^{+2} \rangle$ are
$\preceq^{NF}$-increasing in the successor ordinals and
continuous. \item $M_1=\bigcup \{M_{\alpha,1}:\alpha<\lambda^{+2}
\}$ and $M_2=\bigcup \{M_{\alpha,2}:\alpha<\lambda^{+2} \}$. \item
$M_1,M_2$ are isomorphic.
\end{enumerate}

Then
$S(\bar{M}^1)/D_{\lambda^{+2}}=S(\bar{M}^2)/D_{\lambda^{+2}}$.
\end{proposition}

\begin{proof}
Let $f:M_1 \to M_2$ be an isomorphism. Define $E:=\{\alpha \in
\lambda^{+2}:f[M_{1,\alpha}]=M_{2,\alpha}\}$. So $S(\langle
M_{\alpha,1}:\alpha \in E \rangle)=S(\langle
f[M_{\alpha,1}]:\alpha \in E \rangle )=S(\langle
M_{\alpha,2}:\alpha \in E \rangle)$. By Proposition \ref{9.3}.b
$S(\langle M_{\alpha,1}:\alpha \in E \rangle)=S(\bar{M}^1) \bigcap
E$ and $S(\langle M_{\alpha,2}:\alpha \in E \rangle)=S(\bar{M}^2)
\bigcap E$. Hence $S(\bar{M}^1) \bigcap E=S(\bar{M}^2) \bigcap E$.
\end{proof}

\begin{proposition}\label{enough to assign M^S}
Assume that we can assign to each $S \in
S^{\lambda^{+2}}_{\lambda^+}:=\{S:S$ is a stationary subset of
$\lambda^{+2}$ and $(\forall \alpha \in S)cf(\alpha)=\lambda^+
\}$, a model $M^S \in K_{\lambda^{+2}}$ with
$S(M^S)=S/D_{\lambda^{+2}}$ (especially it is defined).

Then there are $2^{\lambda^{+2}}$ non-isomorphic models in
$K_{\lambda^{+2}}$.
\end{proposition}

\begin{proof}
Since $|S^{\lambda^{+2}}_{\lambda^+}|=2^{\lambda^{+2}}$ it follows
by Proposition \ref{is well defined}.
\end{proof}

The following theorem says that there is a kind of a witness for
non-$\preceq^{NF}$-smoothness, such that if it holds, then there
are $2^{\lambda^{+2}}$ non-isomorphic models in
$K_{\lambda^{+2}}$.

\begin{theorem} \label{9.4}
Suppose that there is an increasing continuous sequence $\langle
M^*_\alpha:\alpha \leq \lambda^++1 \rangle$ of models in $K^{sat}$
such that for each $\alpha,\beta$ with $\alpha<\beta<\lambda^+$ we
have $M^*_\alpha \prec^+ M^*_\beta \preceq^{NF}M^*_{\lambda^++1}$
but $M^*_{\lambda^+}\npreceq^{NF} M^*_{\lambda^++1}$.

Then  there are $2^{\lambda^{+2}}$ pairwise non-isomorphic models
in $K_{\lambda^{+2}}$.
\end{theorem}

\begin{proof}
By Proposition \ref{enough to assign M^S}, it is enough to assign
to each $S \in S^{\lambda^{+2}}_{\lambda^+}$ a model $M^S \in
K_{\lambda^{+2}}$ with $S(M^S)=S/D_{\lambda^{+2}}$. Let $S$ be a
stationary subset of $\lambda^{+2}$ such that $\alpha \in S
\Rightarrow cf(\alpha)=\lambda^+$. We will choose a model
$M_\beta$ by induction on $\beta<\lambda^{+2}$ such that:
\begin{enumerate}
\item $M_\beta \in K^{sat}$. \item The sequence $\langle
M_\beta:\beta<\lambda^{+2} \rangle$ is continuous. \item $\beta
\in \lambda^{+2}-S \Rightarrow M_\beta \prec^+ M_{\beta+1}$. \item
If $\beta \in S$ then $(M_\beta,M_{\beta+1}) \cong
(M^*_{\lambda^+},M^*_{\lambda^++1})$. \item For each
$\beta<\lambda^{+2}$ $M_\beta \preceq^{NF} M_{\beta+1}
\Leftrightarrow \beta \notin S$.
\end{enumerate}

Note that clause 5 is the crucial point and it actually follows by
clauses 3,4.

[Why is it possible to choose $M_\beta$? For $\beta=0$ we choose a
model $M_0 \in K^{sat}$. For limit ordinal $\beta$, define
$M_\beta=\bigcup\{M_\gamma:\gamma<\beta\}$. What will we do in the
$\beta+1$ step? Clause 5 follows by clauses 3,4. So it is enough
to find $M_{\beta+1}$ which satisfies clauses 3,4.

\case{case a:} $\beta \notin S$. In this case we choose
$M_{\beta+1}$ such that $M_\beta\prec^+M_{\beta+1}$ (see
Proposition \ref{7.4}.a).

\case{case b:} $\beta \in S$. Since $M_\beta,M_{\lambda^+}^*$ are
saturated in $\lambda^+$ over $\lambda$, they are isomorphic.
Hence we can find $M_{\beta+1}$ with clause 4]

Define $M^S:=\bigcup \{M_\alpha:\alpha<\lambda^{+2}\}$. It remains
to prove that $S(M^S)=S/D_{\lambda^{+2}}$ (especially $S(M^S)$ is
defined). But if $S(\langle M_\alpha:\alpha<\lambda^{+2} \rangle)$
is defined then by clause 5 $S(M^S)=S(\langle
M_\alpha:\alpha<\lambda^{+2}
\rangle)/D_{\lambda^{+2}}=S/D_{\lambda^{+2}}$. So it is enough to
prove that it is defined, namely to prove that for each
$\alpha,\beta$ with $\alpha<\beta<\lambda^{+2}$ we have
$M_{\alpha+1} \preceq^{NF} M_{\beta+1}$. But it is easier to prove
more:
\begin{claim}\label{(*)_beta}
For every $\beta \leq \lambda^+$ $(*)_\beta$: For each $\alpha$
with $\alpha<\beta$ the following hold:
\begin{enumerate}
\item $M_{\alpha+1} \preceq^{NF} M_{\beta+1}$. \item If $\beta
\notin S$ then $M_{\alpha+1} \prec^{+} M_{\beta+1}$.
\end{enumerate}
\end{claim}

\begin{proof}
$(*)_0$ is vacuous.

Why does $(*)_\beta \Rightarrow (*)_{\beta+1}$ hold? Fix
$\alpha<\beta+1$. We prove that $M_{\alpha+1} \prec^+
M_{\beta+2}$. By clause 3 $M_{\beta+1} \prec^{+} M_{\beta+2}$. So
if $\alpha=\beta$ then $M_{\alpha+1} \prec^+ M_{\beta+2}$. So
without loss of generality $\alpha<\beta$. By $(*)_\beta$
$M_{\alpha+1} \preceq^{NF} M_{\beta+1}$. But $M_{\beta+1} \prec^+
M_{\beta+2}$. So by Proposition \ref{7.10 in 20.11.09}.c
$M_{\alpha+1} \prec^+ M_{\beta+2}$. This establishes
$(*)_{\beta+1}$.

Assume that $\delta$ is a limit ordinal and $(*)_\beta$ holds for
each $\beta$ with $\beta<\delta$. We have to prove $(*)_\delta$.
Let $\langle \gamma(\epsilon):\epsilon<cf(\delta) \rangle$ be an
increasing continuous of ordinals with limit $\delta$, such that
for every $\epsilon,\ \gamma(\epsilon+1)$ is a successor of a
successor ordinal. Note that for every $\epsilon<\cf(\delta)$
$\gamma_\epsilon \notin S$, because
$cf(\gamma_\epsilon)<\cf(\delta) \leq \lambda^+$. Consider the
sequence $\langle M_{\gamma_\epsilon}:\epsilon<cf(\delta)
\rangle$.
\begin{claim}\label{9.4.1}
 $M_{\gamma_\epsilon} \prec^{+}
M_{\gamma_{\epsilon+1}}$ for each $\epsilon<cf(\delta)$.
\end{claim}

\begin{proof}
Since $\gamma_\epsilon \notin S$, by clause 3 $M_{\gamma_\epsilon}
\prec^+ M_{\gamma_\epsilon+1}$. If
$\gamma_{\epsilon+1}=\gamma_\epsilon+1$ then the claim is proved.
Assume $\gamma_{\epsilon+1}>\gamma_\epsilon+1$.
$\gamma_{\epsilon+1}=\zeta+1$ for some successor $\zeta$. $\zeta
\notin S$. So by $(*)_\zeta.2$, $M_{\gamma_\epsilon+1} \prec^+
M_{\zeta+1}=M_{\gamma_{\epsilon+1}}$. So $M_{\gamma_\epsilon}
\prec^+ M_{\gamma_\epsilon+1} \prec^+ M_{\gamma_{\epsilon+1}}$.
Hence by Proposition \ref{7.10 in 20.11.09}.d $M_{\gamma_\epsilon}
\prec^+ M_{\gamma_{\epsilon+1}}$.
\end{proof}

\begin{claim}\label{9.continuous}
The sequence $\langle M_{\gamma_\epsilon}:\epsilon<cf(\delta)
\rangle ^\frown \langle M_\delta \rangle$ is continuous.
\end{claim}

\begin{proof}
Take $\delta' \in \{\gamma_\epsilon:\epsilon<\cf(\delta)\} \bigcup
\{\delta\}$ and take $x \in M_{\delta'}$. We have to find
$\epsilon<cf(\delta)$ such that $\gamma_\epsilon<\delta'$ and $x
\in M_{\gamma_\epsilon}$. By clause 2 the sequence $\langle
M_\beta:\beta<\lambda^{+2} \rangle$ is continuous, so for some
$\beta<\delta'$ $x \in M_\beta$. The ordinals sequence $\langle
\gamma_\epsilon:\epsilon<cf(\delta) \rangle ^\frown \langle \delta
\rangle$ is increasing and continuous. Hence for some
$\epsilon<cf(\delta)$ with $\beta<\gamma_\epsilon<\delta'$. Since
$M_\beta \subseteq M_{\gamma_\epsilon}$, $x \in
M_{\gamma_\epsilon}$.
\end{proof}

\begin{claim} \label{9.11}
$M_{\gamma_\epsilon} \preceq^{NF} M_\delta$ for each
$\epsilon<\cf(\delta)$.
\end{claim}

\begin{proof}
By Proposition \ref{6.5}.d (and Claim \ref{9.4.1}, Claim
\ref{9.continuous} and Proposition \ref{7.10 in 20.11.09}.a).
\end{proof}

Now we return to the proof of $(*)_\delta$. Fix $\alpha<\delta$.

\begin{claim}\label{9.4.3} \label{9.12 in 11.09}
$M_{\alpha+1} \preceq^{NF} M_{\gamma_{\epsilon+1}}$ for some
$\epsilon<cf(\delta)$.
\end{claim}

\begin{proof}
Take $\epsilon<cf(\delta)$ with $\alpha+1<\gamma_{\epsilon+1}$.
$\gamma_{\epsilon+1}=\zeta+1$ for some $\zeta$. So by
$(*)_\zeta.1$ $M_{\alpha+1} \preceq^{NF} M_{\zeta+1}=
M_{\gamma_{\epsilon+1}}$.
\end{proof}

\case{Case a:} $\delta \notin S$. In this case by clause 4
$M_{\delta} \prec^+ M_{\delta+1}$. So by Proposition \ref{7.10 in
20.11.09}.c it is enough to prove that $M_{\alpha+1} \preceq^{NF}
M_\delta$. By Claim \ref{9.4.3} $M_{\alpha+1} \preceq^{NF}
M_{\gamma_{\epsilon+1}}$ for some $\epsilon$. By Claim \ref{9.11}
$M_{\gamma_{\epsilon+1}} \preceq^{NF} M_\delta$. So by Proposition
\ref{6.5}.b $M_{\alpha+1} \preceq^{NF} M_{\delta}$.

\case{Case b:} $\delta \in S$. In this case we have to prove that
$M_{\alpha+1} \preceq^{NF} M_{\delta+1}$. We choose $f_\alpha$ by
induction on $\alpha \leq \lambda^+$ such that:
\begin{enumerate}
\item For every $\alpha \leq \lambda^+$, $f_\alpha:M^{*}_\alpha
\to M_{\gamma_\alpha}$ is an isomorphism. \item $\langle
f_\alpha:\alpha \leq \lambda^+ \rangle$ is an increasing
continuous sequence of isomorphisms.
\end{enumerate}

There is no problem to carry out this induction [Why? We can
choose $f_0$ by Theorem \ref{1.13}, (the uniqueness of the
saturated model in $\lambda^+$ over $\lambda$). $M_\alpha^*
\prec^+ M_{\alpha+1}^*$. By Claim \ref{(*)_beta}
$M_{\gamma_\alpha} \prec^+ M_{\gamma_{\alpha+1}}$. So by Theorem
\ref{7.13}.a, for every $\alpha$, we can find $f_{\alpha+1}$. For
$\alpha$ limit take union].

Now by clause 4, $(M_\delta,M_{\delta+1}) \cong
(M_{\lambda^+}^*,M_{\lambda^++1}^*)$. So we can find an
isomorphism $f:M_{\lambda^++1} \to M_{\delta+1}$ that extends
$f_{\lambda^+}$. For every $\epsilon<\lambda^+$ $M_\epsilon^*
\preceq^{NF} M_{\lambda^++1}^*$, so
$M_{\gamma_\epsilon}=f[M_\epsilon^*] \preceq^{NF}
f[M_{\lambda^++1}^*]=M_{\delta+1}$. So $M_{\gamma_\epsilon}
\preceq M_{\delta+1}$ for each $\epsilon<cf(\delta)$. Hence
$M_{\gamma_{\epsilon+1}} \preceq^{NF} M_{\delta+1}$ for each
$\epsilon<cf(\delta)$. But by Claim \ref{9.12 in 11.09} for some
$\epsilon<cf(\delta)$ $M_{\alpha+1} \preceq^{NF}
M_{\gamma_{\epsilon+1}}$. Therefore by Proposition \ref{6.5}.b
$M_{\alpha+1} \preceq ^{NF} M_{\delta+1}$.
\end{proof}
\end{proof}

\begin{theorem} \label{9.5}
The following conditions are equivalent:
\begin{enumerate}[(a)]
\item $(K^{sat},\preceq^{NF} \restriction K^{sat})$ does not
satisfy smoothness. \item There are $M^*_1,M^*_2 \in K^{sat}$ such
that $M^*_1 \preceq^\otimes M^*_2$ but $M^*_1 \nprec^{NF} M^*_2$.
\item There is a sequence $\langle M_\epsilon:\epsilon \leq
\lambda^++1 \rangle$ of models in $K^{sat}$ such that for each
$\epsilon,\zeta$ with $\epsilon<\zeta \leq \lambda^++1$ we have
$\epsilon \neq \lambda^+ \Leftrightarrow M_\epsilon \prec^+M_\zeta
\Leftrightarrow M_\epsilon \preceq^{NF} M_\zeta$.
\end{enumerate}
\end{theorem}
\begin{proof} $c \Rightarrow a$
is clear. $b \Rightarrow c$ holds by Proposition \ref{8.4}. $a
\Rightarrow b$ holds by Proposition \ref{8.2}.b.
\end{proof}

Now we can prove Theorem \ref{9.6}, but first we remind it: If
$(K^{sat},\preceq^{NF} \restriction K^{sat})$ does not satisfy
smoothness, then there are $2^{\lambda^{+2}}$ pairwise
non-isomorphic models in $K_{\lambda^{+2}}$.

\begin{proof}
Condition a of Theorem \ref{9.5} is satisfied, so condition c is
satisfied too. Hence by Theorem \ref{9.4} we have the conclusion
of the theorem.
\end{proof}

\section{a good $\lambda^+$-frame}
\discussion{Discussion:} In Definitions \ref{preparation for
forking for big models}, \ref{forking for big models} and
\ref{basic for big models}
we expanded the definition of the
non-forking relation and basic types to models in $K_{>\lambda}$.
In Theorem \ref{2.10} we proved some axioms of a good frame for
this expansions. Here we are going to prove the
other axioms. So why are sections 3-9 needed? 
 In other words, what are the difficulties in proving that $S^+$ (defined below) is a good
$\lambda^+$-frame? The main problem is that amalgamation may not
hold in $(K_{\lambda^+},\preceq \restriction K_{\lambda^+})$. Now
we can overcome this problem by restricting the relation
$\preceq_{K_{\lambda^+}}$ to the relation $\preceq^{NF}$. But then
there is a problem with smoothness. We overcome this problem by
showing that non-smoothness is a non-structure property, see
section 9. For the non-structure theorem, we had to restrict to
the class of saturated models in $\lambda^+$ over $\lambda$. Now
the relation $\prec^+$ and the locality enable use to prove the
remaining axioms.

\begin{definition} \label{10.1}
Let $\frak{s}$ be a semi good frame. We say that $\frak{s}$ is
\emph{successful} when:
\begin{enumerate}
\item $\frak{s}$ is \emph{weakly successful} (i.e. we have
existence for $K^{3,uq}_\frak{s}$). \item $(K^{sat},\preceq^{NF}
\restriction K^{sat})$ satisfies smoothness.
\end{enumerate}
\end{definition}

\begin{hypothesis} \label{10.2}
$\frak{s}$ is a successful semi-good $\lambda$-frame with
conjugation.
\end{hypothesis}

We remind that the types in this paper are classes of triples
under some equivalence relation. But this relation depends on the
partial order, we define on the class of models.
For $M_0,M_1 \in K_{\lambda^+}$ when we write $tp(a,M,N)$ we mean
to the partial order $\preceq$. But when we want to consider the
partial order $\preceq^{NF}$ we have to write it explicitly.

\begin{definition}
For $M_0,M_1 \in K^{sat}$ and $a \in M_1-M_0$ we define
$$tp^+(a,M_0,M_1):=tp_{((K^{sat})^{up},(\preceq^{NF} \restriction
K^{sat})^{up})}(a,M_0,M_1).$$
\end{definition}
(About `$sat$' see Definition \ref{7.1} (page \pageref{7.1}) and
about `$up$' see Definition \ref{up} (page \pageref{up})).

\begin{proposition} \label{the type does not depend on the aec}
For every $M_0,M_1,M_2$ with $M_0 \preceq^{NF} M_1 \wedge M_0
\preceq^{NF} M_2$ and every $a_1,a_2$ with $a \in M_1-M_0 \wedge
a_2 \in M_2-M_0$:
$$tp^+(a_1,M_0,M_1)=tp^+(a_2,M_0,M_2) \Leftrightarrow tp(a_1,\allowbreak M_0,M_1)=tp(a_2,\allowbreak
M_0,M_2).$$
\end{proposition}

\begin{proof}
The first direction: Suppose
$tp^+(a_1,M_0,M_1)=tp^+(a_2,M_0,M_2)$. By Theorem \ref{7.8}.c
(page \pageref{7.8}) $(K^{sat},\preceq^{NF} \restriction K^{sat})$
has amalgamation. So there are $f_1,f_2,M_3$ such that: $M_0
\preceq^{NF}M_3,\ f_n:M_n \to M_3$ is a $\preceq^{NF}$-embedding
over $M_0$ and $f_1(a_1)=f_2(a_2)$. But $K^{sat} \subseteq K$, and
the relation $\preceq^{NF}$ is included in the relation $\preceq$
so the amalgamation $(f_1,f_2,M_3)$ witnesses that
$tp(a_1,M_0,M_1)=tp(a_2,M_0,M_2)$.

The second direction: Suppose $tp(a_1,M_0,M_1)=tp(a_2,M_0,M_2)$.
Take an amalgamation $(f_1,f_2,M_3)$ of $M_1,M_2$ over $M_0$ with
$f_1(a_1)=f_2(a_2)$. For each $N \in K_\lambda$ with $N \preceq
M_0$ $tp(f_1(a_1),N,f_1[M_1])=tp(f_2(a_2),N,f_2[M_2])$. So by
Theorem \ref{7.13}.b $tp^+(a_1,M_0,M_1)=tp^+(a_2,M_0,M_2)$.
\end{proof}

Although we defined restriction of types in Definition
\ref{definition of a type}.3 (on page \pageref{definition of a
type}), the following definition is needed:
\begin{definition}
For $p=tp^+(a,M_0,M_1)$ and $N \in K_\lambda$ with $N \preceq M_0$
we define $p \restriction N:=tp(a,N,M_1)$.
\end{definition}

The following definition is based on Definition \ref{2.9} (page
\pageref{2.9}).
\begin{definition} \label{10.3}
$\frak{s}^+:=((K^{sat})^{up},(\preceq^{NF} \restriction
K^{sat})^{up}, \frak{s}^{bs,+},\dnf^+)$, where:
\begin{enumerate}
 \item
For each $M \in K^{sat}$ we define
$S^{bs,+}(M):=\{tp^+(a,M,N):\{M,N\} \subseteq K^{sat},\ M
\preceq^{NF} N,\ tp(a,M,N) \in S^{bs}_{
>\lambda}\}$ \item $\dnf^+$ is defined by: $tp^+(a,M_1,M_2)$ does not fork over $M_0$
if $\{M_0,M_1,M_2\} \allowbreak \subseteq K^{sat},\ M_0
\preceq^{NF} M_1 \preceq^{NF} M_2$ and $tp(a,M_1,M_2)$ does not
fork over $M_0$.
\end{enumerate}
\end{definition}

\begin{proposition} \label{10.4}
\mbox{}
\begin{enumerate}[(a)]
\item $S^{bs}$ is well defined: It does not depend on the triple
$(M_0,M_1,a)$ that represents the type. \item $\dnf^+$ is well
defined: It does not depend on the triple $(M_0,M_1,a)$ that
represents the type.
\end{enumerate}
\end{proposition}

\begin{proof}
By Proposition \ref{the type does not depend on the aec}.
\end{proof}

\begin{proposition} \label{10.5}
Let $\frak{s}$ be a successful semi-good $\lambda$-frame with
conjugation.
 \begin{enumerate}
\item $(K^{sat},\preceq^{NF} \restriction K^{sat})$ satisfies
axiom c of a.e.c. in $\lambda^+$ (i.e. Definition \ref{1.1}.2.c).
\item $(K^{sat},\preceq^{NF} \restriction K^{sat})$ is an a.e.c.
in $\lambda^+$. \item $(K^{sat},\preceq^{NF} \restriction
K^{sat})$ satisfies the amalgamation property.
\end{enumerate}
\end{proposition}

\begin{proof}
By Theorem \ref{7.8} and hypothesis \ref{10.2}.
\end{proof}

\begin{theorem} \label{10.6}
Let $\frak{s}$ be a successful semi-good $\lambda$-frame with
conjugation. Then $\frak{s}^+$ is a good $\lambda^+$-frame.
\end{theorem}

\begin{proof}
By Proposition \ref{10.5} $(K^{sat},\preceq^{NF} \restriction
K^{sat})$ is an a.e.c. in $\lambda^+$ with amalgamation. So by
Fact \ref{1.14} (page \pageref{1.14})
$((K^{sat})^{up},(\preceq^{NF} \restriction K^{sat})^{up})$ is an
a.e.c. with LST number $\lambda^+$. By Theorem \ref{1.13} (page
\pageref{1.13}) $K^{sat}$ is categorical. So
$(K^{sat},\preceq^{NF} \restriction K^{sat})$ has joint embedding.
By Proposition \ref{7.4}.a (page \pageref{7.4}) and Proposition
\ref{7.10 in 20.11.09}.a there is no $\preceq^{NF}$-maximal model
in $K^{sat}$. What about the axioms of the basic types and the
non-forking relation? By Theorem \ref{2.10} the following axioms
are satisfied: Density, monotonicity, local character and
continuity.

\begin{proposition}
$\frak{s}^+$ satisfies basic stability.
\end{proposition}

\begin{proof}
Let $M \in K^{sat}$. $M \in K_{\lambda^+}$, so it has a
representation $\langle  N_\alpha:\alpha \in \lambda^+ \rangle$.
For $p \in S^{bs,+}(M)$ define $(\alpha_p,q_p)$ by: $\alpha_p$ is
the minimal ordinal in $\lambda^+$ such that $p$ does not fork
over $N_\alpha$. $q_p=:p\restriction N _{\alpha_p}$. For every
$\alpha \in \lambda^+$ we have $|S^{bs}(N_\alpha)| \leq
\lambda^+$, so $|{(\alpha_p,q_p):p \in S^{bs,+}(M)}| \leq
\lambda^+ \times \lambda^+ = \lambda^+$. So it is sufficient to
prove that the function $p \to (\alpha_p,q_p)$ is an injection.
For every $p_1,p_2 \in S^{bs,+}(M)$ if $\alpha_{p_1}=\alpha_{p_2}
\wedge q_{p_1}=q_{p_2}$ Then by Corollary \ref{7.7}.b (locality,
page \pageref{7.7}) $p_1=p_2$.
\end{proof}

\begin{proposition}\label{uniqueness in s^+}
\mbox{}
\begin{enumerate}
\item If \begin{enumerate} \item $N \in K_\lambda$ and $M \in
K_{\lambda^+}$. \item For $n=1,2$ $p_n \in S^{bs,+}(M)$ and does
not fork over $N$. \item $p_1 \restriction N=p_2 \restriction N$.
\end{enumerate}
Then $p_1=p_2$. \item $\frak{s}^+$ satisfies uniqueness.
\end{enumerate}
\end{proposition}

\begin{proof}
\mbox{}\\
1) By the proof of Corollary \ref{7.7}.b (locality, page
\pageref{7.7}). Remember that if $N_0 \preceq ^{NF} N_1 \prec^+
N_2$ then $N_0 \preceq N_2$
(By Proposition \ref{before 7.4}.c) .\\
2) Suppose $n<2 \Rightarrow M_n \in K^{sat},\ M_0 \preceq M_1,\
p,q \in S^{bs,+}(M_1),\ \allowbreak p\restriction
M_0=q\restriction M_0$ and $p,q$ does not fork over $M_0$. By the
definition of $\dnf^+$, there are $N_p,N_q \in K_\lambda$, such
that $N_p \preceq M_0,\ N_q \preceq M_0$, $p$ does not fork over
$N_p$ and $q$ does not fork over $N_q$. As $LST(K,\preceq) \leq
\lambda$, there is a model $N \in K_\lambda$ with $N_p \bigcup N_q
\subseteq N \preceq M_0$. By axiom \ref{1.1}.1.e $N_p \preceq N$
and $N_q \preceq N$. By Theorem \ref{2.10}(2) (monotonicity, page
\pageref{2.10}), $p,q$ does not fork over $N$. By the assumption
$p\restriction M_0=q\restriction M_0$, so $p\restriction N
=q\restriction N$. Hence by item 1, $p=q$.
\end{proof}

\begin{proposition}
$\frak{s}^+$ satisfies symmetry.
\end{proposition}

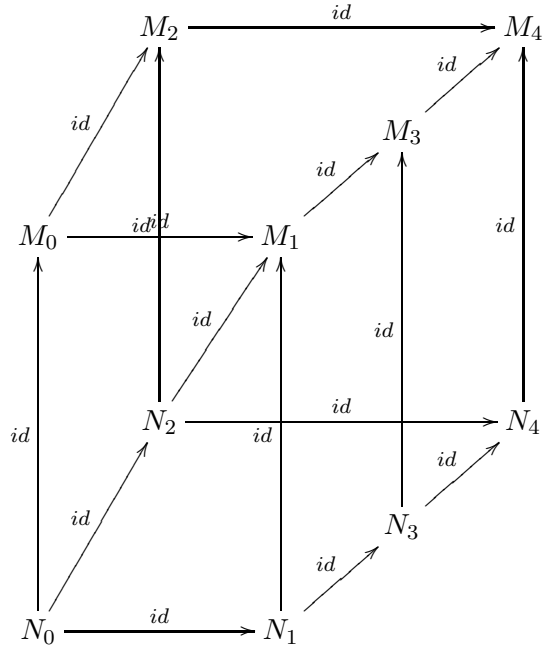
\begin{proof}
\begin{displaymath}
\xymatrix{
& M_2 \ar[rrr]^{id} &&& M_4\\
&&& M_3 \ar[ru]^{id}\\
M_0 \ar[uur]^{id} \ar[rr]^{id} && M_1 \ar[ru]^{id}\\ \\
& N_2 \ar[uur]^{id} \ar[rrr]^{id} \ar[uuuu]^{id} &&& N_4 \ar[uuuu]^{id}\\
&&& N_3 \ar[uuuu]^{id} \ar[ur]^{id}\\
N_0 \ar[uuuu]^{id} \ar[uur]^{id} \ar[rr]^{id} && N_1
\ar[uuuu]^{id} \ar[ur]^{id} }
\end{displaymath}

Suppose 1-5 where:\\
(1) $\{M_0,M_1,M_3\} \subseteq K^{sat}$.\\
(2) $M_0 \preceq^{NF}M_1 \preceq^{NF}M_3$.\\
(3) $tp(a_1,M_0,M_3) \in S^{bs,+}(M_0)$.\\
(4) $a_1 \in M_1$.\\
(5) $tp(a_2,M_1,M_3)$ does not fork over $M_0$.\\

\step{Step a:} We choose models $N_0,N_1,N_3 \in K_\lambda$
which satisfies 6-12 where:\\
(6) $n \in \{0,1,3\} \Rightarrow N_n \preceq M_n$ and $N_0 \preceq N_1 \preceq N_3$.\\
(7) $tp(a_2,M_1,M_3)$ does not fork over $N_0$.\\
(8) $tp(a_1,M_0,M_3)$ does not fork over $N_0$.\\
(9) $a_1 \in N_1$.\\
(10) $a_2 \in N_3$.\\
(11) $\widehat{NF}(N_0,N_1,M_0,M_1)$.\\
(12) $\widehat{NF}(N_1,N_3,M_1,M_3)$.\\
(Why is it possible? By 2, there are representations $\langle
N_{0,\alpha}:\alpha<\lambda^+ \rangle,\ \langle
N_{1,\alpha}:\alpha<\lambda^+ \rangle,\ \langle
N_{1,\alpha}^*:\alpha<\lambda^+ \rangle,\ \langle
N_{3,\alpha}:\alpha<\lambda^+ \rangle$ of $M_0,M_1,M_1,M_3$
respectively, such that: $\alpha<\lambda^+ \Rightarrow
NF(N_{0,\alpha},N_{1,\alpha},N_{0,\alpha+1},N_{1,\alpha+1}),\
NF(N_{1,\alpha}^*,N_{3,\alpha},\allowbreak
N_{1,\alpha+1}^*,N_{3,\alpha+1})$. Let $E$ be a club of
$\lambda^+$ such that $\alpha \in E \Rightarrow
N_{1,\alpha}=\N_{1,\alpha}^*$. Choose $\alpha \in E$ big enough
such that 7,8,9,10 will satisfied for $N_0=N_{0,\alpha}$
$N_1=N_{1,\alpha},\ N_3=N_{3,\alpha}$)\\

\step{Step b:} [We use the symmetry axiom] By 6,8 we
have:\\
(13) $tp(a_1,N_0,N_3) \in S^{bs}(N_0)$.\\
by 6,7 we have:\\
(14) $tp(a_2,N_1,N_3)$ does not fork over $N_0$.\\
Now by Definition \ref{2.1a}.3.e (symmetry) there are $N_2^*,N_4^*
\in
K_\lambda$ which satisfies 15-18:\\
(15) $N_0 \preceq N_2^* \preceq N_4^*$.\\
(16) $N_3 \preceq N_4^*$.\\
(17) $a_2 \in N_2^*$.\\
(18) $tp(a_1,N_2^*,N_4^*)$ does not fork over $N_0$.\\

\step{Step c:} [Move everything to $K^{sat}$]\\
We can choose $f$ which satisfies 19,20:\\
(19) $f$ is an injection, $dom(f)=N_4^*$ and $f \restriction N_3$ is the identity.\\
(20) $f[N_4^*] \bigcap M_3=N_3$.\\
Define $N_4:=f[N_4^*],\ N_2:=f[N_2^*]$. By the existence
proposition of the $\prec^+$-extensions (Proposition \ref{7.4}.c),
there is $M_4 \in K_\lambda$
which satisfies 21,22:\\
(21) $\widehat{NF}(N_3,N_4,M_3,M_4)$.\\
(22) $M_3\prec^+M_4$.\\
By 20 (mainly) we know:\\
(23) $N_2 \bigcap M_0=N_0$.\\
(Why? By 15 and the definitions of $f,N_2$, we have $N_0 \preceq
N_2$. By 6 $N_0 \preceq M_0$. Let $x \in N_2 \bigcap M_0$. By 2,15
$x \in N_4 \bigcap M_3$. So by 20 $x \in N_3$. So $x \in N_3
\bigcap M_1$. Hence by 12, $x \in N_1$. So $x \in N_1 \bigcap
M_0$. Hence by 11, we have $x \in N_0$). So by the existence
proposition of
$\widehat{NF}$ (Proposition \ref{5.15}.c), there is $M_2 \in K^{sat}$ such that:\\
(24) $\widehat{NF}(N_0,N_2,M_0,M_2)$.\\
Without loss of generality $N_4 \bigcap M_2=N_2$ as $M_0 \bigcap
N_4=N_0$. By Proposition \ref{7.4}.b there is $M_6 \in K^{sat}$
which satisfies 25,26:\\
(25) $M_2\prec^+M_6$.\\
(26) $\widehat{NF}(N_2,N_4,M_2,M_6)$.\\

\step{Step d:}
We will prove 27,28:\\
(27) $tp(a_1,M_2,M_6)$ does not fork over $N_0$.\\
(28) There is an isomorphism $g:M_6 \to M_4$ over $M_0
\bigcup N_2$.\\
Then we will conclude:\\
(29) $tp(a_1,g[M_2],M_4)$ does not fork over $M_0$. By 25,
Proposition \ref{7.4}.c=7.10 and 24 we have 30:\\
(30) $M_0\prec^+M_6$. \\
By 24,25 and Theorem \ref{5.15}.b (monotonicity, on page
\pageref{5.15}):
(31) ${NF}(N_0,N_2,M_0,M_6)$.\\
By 24,26,28 and the transitivity of the
relation $\widehat{NF}$ we have:\\
(32)${NF}(N_0,N_2,M_0,M_4)$.\\
By 2,22 and Proposition \ref{before 7.4}.c:\\
(33) $M_0\prec^+M_4$.

By 30-33 and Theorem \ref{7.13}.c, we know 28. By 26, and
Theorem \ref{5.15}.e (respecting the frame, page \pageref{5.15}):\\
(34) $tp(a_1,M_2,M_6)$ does not fork over $N_2$. By 18 (and
12,9,19):\\
(35) $tp(a_1,N_2,N_4)$ does not fork over $N_0$. By 26 $N_4
\preceq M_6$, so by Theorem \ref{2.10}(3) (the transitivity of the
non-forking relation), we have:\\
(27) $tp(a_1,M_2,M_6)$ does not fork over $N_0$.\\

\step{Step e:}\\
It remains to prove \\
(36) $a_2 \in g[M_2]$. By 28 , g is an isomorphism over $N_2$, so
it is sufficient to prove $a_2 \in N_2$. By 17 $a_2 \in N_2^*$. So
by
10,19 $a_2 \in N_2$.\\
\end{proof}

By the following proposition, $\frak{s}^+$ satisfies extension.
\begin{proposition}
\mbox{}
\begin{enumerate}
\item If $N \preceq M \in K^{sat},\ p \in S^{bs}(N),\ N \in
K_\lambda$, then there is $q \in S^{bs,+}(M)$ such that
$q\restriction N=p$ and $q$ does not fork over $N$. \item If
$\{M_0,M_1\} \subseteq K^{sat},\ M_0 \preceq^{NF} M_1,\ p \in
S^{bs,+}(M_0)$ than there is an extension of $p$ to
$S^{bs,+}(M_1)$.
\end{enumerate}
\end{proposition}

\begin{proof}
\mbox{}
\begin{enumerate}
\item Let $a,N_1$ be such that $tp(a,N,N_1)=p$. By Theorem
\ref{5.15}.c (page \pageref{5.15}) without loss of generality
there is a model $M_1$ such that $\widehat{NF}(N,N_1,\allowbreak
M,M_1)$. By Theorem \ref{5.15}.e $q:=tp(a,M,M_1)$ does not fork
over $N$. \item By the definition of $S^{bs,+}$, there is a model
$N \in K_\lambda$ such that $N \preceq M_0$ and $p$ does not fork
over N. By item (1), there is $q \in S^{bs,+}(M_1)$ which does not
fork over $N$, and $q\restriction N=p\restriction N $. $q$ does
not fork over $M_0$ as it does not fork over $N$. So it is
sufficient to prove that $q_0:=q \restriction M_0=p$. By Theorem
\ref{2.10}.2 (monotonicity), $q_0$ does not fork over $N$.
$q_0\restriction N=q\restriction N=p\restriction N$. Hence by
Corollary \ref{7.7}.b (locality) $p=q_0$.
\end{enumerate}
\end{proof}
This ends the proof of Theorem \ref{10.6}.
\end{proof}

\section{Conclusions}
\begin{theorem} \label{11.1}
Suppose:
\begin{enumerate}
\item $\frak{s}=(K,\preceq,S^{bs},\dnf)$ is a semi-good
$\lambda$-frame with conjugation. \item
$I(\lambda^{+2},K)<\mu_{unif}(\lambda^{+2},2^{\lambda^+})$. \item
$2^\lambda<2^{\lambda^+}<2^{\lambda^{+2}}$, and $WdmId(\lambda^+)$
is not saturated in $\lambda^{+2}$.
\end{enumerate}
Then
\begin{enumerate}
\item There is a good $\lambda^+$-frame
$\frak{s}^+=((K^{sat},\preceq^{NF}\restriction
K^{sat})^{up},S^{bs,+},\dnf^+)$, such that $K^{sat} \subseteq
K_{\lambda^+}$ and the relation $\preceq^{NF}\restriction K^{sat}$
is included in the relation $\preceq \restriction K^{sat}$. \item
$\frak{s}^+$ has the conjugation property. \item There is a model
in $K$ of cardinality $\lambda^{+2}$. \item There is a model in
$K$ of cardinality $\lambda^{+3}$.
\end{enumerate}
\end{theorem}

\begin{proof}
\mbox{} (1) By Corollary \ref{4.18} (page \pageref{4.18})
$\frak{s}$ is weakly successful in the density sense. $\frak{s}$
has conjugation, so by Proposition \ref{4.4} (page \pageref{4.4}),
$\frak{s}$ is weakly successful. By clause 2 of our assumption,
$I(\lambda^{+2},K)<\mu_{unif}(\lambda^{+2},2^{\lambda^+})$. But by
Proposition \ref{unifleq} $\mu_{unif}(\lambda^{+2},2^{\lambda^+})
\leq 2^{\lambda^{+2}}$. So $I(\lambda^{+2},K)<2^{\lambda^{+2}}$.
Hence by Theorem \ref{9.6} (page \pageref{9.6}),
$(K^{sat},\preceq^{NF} \restriction K^{sat})$ satisfies
smoothness, i.e. $\frak{s}$ is successful (Definition \ref{10.1}).
So Hypothesis \ref{10.2} is satisfied. Therefore by Theorem
\ref{10.6}, $\frak{s}^+$ is a good $\lambda^+$-frame. Obviously
$K^{sat} \subseteq K_{\lambda^+}$ and $\preceq^{NF}$ is included
in the relation $\preceq \restriction K_{\lambda^+}$.

(2) Why does $\frak{s}^+$ have conjugation? Suppose $M_0
\preceq^{NF}M_1,\ \{M_0,M_1\} \subseteq K^{sat}$ and $p \in
S^{bs,+}(M_1)$ does not fork over $M_0$. By the definition of
$\dnf^+$, there is $N \in K_\lambda$ such that $N \preceq M_0$ and
$p$ does not fork over $N$.

\begin{displaymath}
\xymatrix{p \restriction M_0 & f(p \restriction M_0)=p \\ M_0
\ar[r]^{id}_{f} & M_1
\\ N \ar[u]^{id} }
\end{displaymath}

By Theorem \ref{1.13}.a (the uniqueness of the saturated model),
there is an isomorphism $f:M_0 \to M_1$ over $N$. By Theorem
\ref{2.10}(2) (monotonicity), $p\restriction M_0$ does not fork
over $N$. So $f(p\restriction M_0)$ does not fork over $N$. But
also $p$ does not fork over $N$ and $f(p\restriction
M_0)\restriction N=(p\restriction M_0)\restriction N=p\restriction
N$ [Why do we have the first equality? There are $M_0^+,f^+,a$
such that $p \restriction M_0=tp(a,M_0,M_0^+)$ and $f \subseteq
f^+,\ dom(f^+)=M_0^+$. So $(p \restriction M_0) \restriction
N=tp(a,N,M_0^+)=tp(f^+(a),N,f^+[M_0^+])=tp(f^+(a),M_1,f^+[M_0^+])
\restriction N=f(p \restriction M_0) \restriction N$]. So by
\ref{uniqueness in s^+}(1), $f(p\restriction M _0)=p$.

(3) By Proposition \ref{3.3}.3 (page \pageref{3.3}).

(4) Substitute $\frak{s}^+$ instead of $\frak{s}$ in Proposition
\ref{3.3}.3.
\end{proof}

\begin{corollary} \label{11.2}
Suppose:
\begin{enumerate}
\item $n<\omega$. \item $\frak{s}=(K,\preceq,S^{bs},\dnf)$ is a
semi-good $\lambda$-frame with conjugation. \item $m<n \Rightarrow
I(\lambda^{+(2+m)},K)<\mu_{unif}(\lambda^{+(2+m)},2^{\lambda^{+(1+m)}})$.
\item For every $m<n$,
$2^\lambda<2^{\lambda^+}<2^{\lambda^{+2}}<...2^{\lambda^{+(1+n)}}$
and $WdmId(\lambda^{+1+m})$ is not saturated in
$\lambda^{+(2+m)}$.
\end{enumerate}
\emph{then} there is a good $\lambda^{+n}$-frame
$\frak{s}^n=:((K^n,\leq^n),S^{bs,+n},\dnf^{+n})$, such that:
\begin{enumerate}
\item $K^n_{\lambda^{+n}} \subseteq K_{\lambda^{+n}},\ \leq^n
\subseteq \preceq^k\restriction K^n$. \item $\frak{s}^n$ has
conjugation. \item There is a model in $K^n$ of cardinality
$\lambda^{+(2+n)}$.
\end{enumerate}
\end{corollary}
\begin{proof}
By induction on $n$, using Theorem \ref{11.1}.
\end{proof}

\begin{corollary} \label{11.2omega}
Suppose:
\begin{enumerate}
\item $\frak{s}=(K,\preceq,S^{bs},\dnf)$ is a semi-good
$\lambda$-frame with conjugation. \item $m<\omega \Rightarrow
I(\lambda^{+(2+m)},K)<\mu_{unif}(\lambda^{+(2+m)},2^{\lambda^{+(1+m)}})$.
\item $2^{\lambda^{+m}}<2^{\lambda^{+m+1}}$ and for every
$m<\omega$, $WdmId(\lambda^{+1+m})$ is not saturated in
$\lambda^{+(2+m)}$.
\end{enumerate}
\emph{Then} there is a model in $K^n$ of cardinality
$\lambda^{+n}$ for every $n<\omega$.
\end{corollary}
\begin{proof}
By Corollary \ref{11.2}.
\end{proof}

\section{Comparison to \cite{sh600}} A reader who knows
\cite{sh600}, might ask about the main problems in writing our
paper. As in \cite{sh600}, there is a wide use of brimmed
extensions (i.e. using stability), we had to find alternatives.

First the relation $NF$ is defined in \cite{sh600} using brimness,
so we found a natural definition (maybe an easier one) which is
equivalent to the definition in \cite{sh600}, but not using
brimness.

Another problem was proving conjugation (see Definition \ref{2.4},
page \pageref{2.4}). But in the main examples there is
conjugation, so it is reasonable to assume conjugation.

Another problem was to find a relation $\prec^+$ on $k^{nice}$
which satisfies the required properties (see the discussion before
Definition \ref{7.3}, page \pageref{7.3}). In \cite{sh600} it uses
essentially brimness. But as the needed relation is on models of
cardinality $\lambda^+$, We can find such a relation, using just
weak stability.





\subsection*{ Acknowledgment} We thank Boaz Tsaban and Alon Sitton for their useful suggestions and
comments. We thank the referee for doing an outstanding and
extremely conscientious job in improving the paper, for his useful
suggestions and for an example he pointed out.

\end{document}